\documentclass{amsart}
\usepackage{amssymb,amscd,graphics,axodraw,verbatim}

\numberwithin{equation}{section}

\newtheorem{prop}{Proposition}
\newtheorem{theorem}[prop]{Theorem}
\newtheorem{cor}[prop]{Corollary}
\newtheorem{lemma}[prop]{Lemma}
\newtheorem{conjecture}[prop]{Conjecture}

\theoremstyle{definition}
\newtheorem{example}[prop]{Example}
\newtheorem{remark}[prop]{Remark}

\numberwithin{prop}{section}

\newcommand{\A}{A_{2n-1}^{(1)}}
\newcommand{\AZ}{\mathbb{A}}
\newcommand{\At}{A_{2n}^{(2)}}
\newcommand{\Atd}{A_{2n}^{(2)\dagger}}
\newcommand{\BZ}{\mathbb{B}}
\newcommand{\CC}{\mathbb{C}}
\newcommand{\Cn}{C_n^{(1)}}
\newcommand{\Conf}{\mathrm{C}}
\newcommand{\Dt}{D_{n+1}^{(2)}}
\newcommand{\Hh}{\widehat{H}}
\newcommand{\Hom}{\mathrm{Hom}}
\newcommand{\Image}{\mathrm{Im}}
\newcommand{\Jb}{\overline{J}}
\newcommand{\Knuth}{\sim}
\newcommand{\LRT}{\mathrm{LR}}
\newcommand{\La}{\Lambda}
\newcommand{\Lab}{\overline{\La}}
\newcommand{\Lat}{\tilde{\La}}
\newcommand{\N}{\mathbb{N}}
\newcommand{\Path}{\mathcal{P}}
\newcommand{\Pb}{\overline{P}}
\newcommand{\Pcl}{P_{\mathrm{cl}}}
\newcommand{\Pfin}{\overline{P}}
\newcommand{\Plac}{\mathrm{Plac}}
\newcommand{\PSchen}{\mathbb{P}}
\newcommand{\Qb}{\overline{Q}}
\newcommand{\QQ}{\mathbb{Q}}
\newcommand{\QSchen}{\mathbb{Q}}
\newcommand{\RC}{\mathrm{RC}}
\newcommand{\RCtr}{\mathrm{RCtr}}
\newcommand{\RR}{\mathbb{R}}
\newcommand{\Rh}{\widehat{\R}}
\newcommand{\Rt}{\tilde{R}}
\newcommand{\R}{\sigma}
\newcommand{\Symm}{\mathfrak{S}}
\newcommand{\V}{V}
\newcommand{\Vh}{\widehat{\V}}
\newcommand{\Vho}{{\Vh}_{\mathrm{op}}}
\newcommand{\Vo}{\V_{\mathrm{op}}}
\newcommand{\Wb}{\overline{W}}
\newcommand{\X}{X}
\newcommand{\ZC}{Z}
\newcommand{\Z}{\mathbb{Z}}
\newcommand{\alphat}{\tilde{\alpha}}
\newcommand{\autoA}{{\psi}}
\newcommand{\autoAtwo}{{\psi_{A^{(2)}}}}
\newcommand{\autoC}{{\psi_C}}
\newcommand{\ba}[1]{\overline{#1}}
\newcommand{\bt}{\tilde{b}}
\newcommand{\cc}{cc}
\newcommand{\cd}{\vee}
\newcommand{\charge}{\mathrm{charge}}
\newcommand{\cl}{\mathrm{cl}}
\newcommand{\cont}{\mathrm{cont}}
\newcommand{\ed}{*}

\newcommand{\et}{\widetilde{e}}
\newcommand{\ev}{\mathrm{ev}}
\newcommand{\flipq}{\mathrm{comp}}
\newcommand{\floor}[1]{\lfloor #1 \rfloor}

\newcommand{\ft}{\widetilde{f}}

\newcommand{\gfin}{\overline{\gggg}}
\newcommand{\gggg}{\mathfrak{g}}
\newcommand{\gt}{\tilde{\gggg}}
\newcommand{\inner}[2]{\langle #1\,,\,#2\rangle}
\newcommand{\la}{\lambda}
\newcommand{\lb}{\overline{\ell}}
\newcommand{\len}[1]{\left|{#1}\right|}
\newcommand{\lev}{\mathrm{level}}
\newcommand{\lt}{\tilde{\ell}}
\newcommand{\mcbn}[1]{\multicolumn{1}{|c|}{#1}}
\newcommand{\mcb}{\multicolumn{1}{|c|}{}}

\newcommand{\mcrn}[1]{\multicolumn{1}{c|}{#1}}
\newcommand{\mcr}{\multicolumn{1}{c|}{}}
\newcommand{\mc}{\multicolumn{1}{c}{}}
\newcommand{\nub}{\overline{\nu}}
\newcommand{\pd}{\wedge}
\newcommand{\perm}{\mathrm{perm}}
\newcommand{\phib}{\overline{\phi}}
\newcommand{\phit}{\widetilde{\phi}}
\newcommand{\pos}{\mathrm{pos}}

\newcommand{\qbin}[2]{\genfrac{[}{]}{0pt}{}{#1}{#2}}
\newcommand{\rows}{\mathrm{rows}}
\newcommand{\shape}{\mathrm{shape}}
\newcommand{\slwt}{\mathrm{wt}_{\mathrm{sl}}}
\newcommand{\st}{\tilde{s}}
\newcommand{\wt}{\mathrm{wt}}
\newcommand{\wtembed}{\Psi}

\newcommand{\embedfin}{\overline{\Psi}}
\newcommand{\wtflip}{\psi_C}
\renewcommand{\emptyset}{\varnothing}

\begin{document}

\title{Virtual crystals and fermionic formulas of type $\Dt$,
$\At$, and $\Cn$}

\author[M.~Okado]{Masato Okado}
\address{Department of Informatics and Mathematical Science,
Graduate School of Engineering Science, Osaka University,
Toyonaka, Osaka 560-8531, Japan}
\email{okado@sigmath.es.osaka-u.ac.jp}

\author[A.~Schilling]{Anne Schilling}
\address{Department of Mathematics 2-279, M.I.T., Cambridge, MA 02139, U.S.A.
and\newline
Department of Mathematics, University of California, One Shields
Avenue, Davis, CA 95616-8633, U.S.A.}
\email{anne@math.mit.edu, anne@math.ucdavis.edu}

\author[M.~Shimozono]{Mark Shimozono${}^{\dagger}$}
\address{Department of Mathematics, 460 McBryde Hall, Virginia Tech,
Blacksburg, VA 24061-0123, U.S.A}
\email{mshimo@math.vt.edu}
\thanks{${}^{\dagger}$ Partially supported by NSF grant DMS-9800941.}

\subjclass{Primary 81R50 81R10 17B37; Secondary 05A30 82B23}

\keywords{Crystal bases, quantum affine Lie algebras,
fermionic formulas, rigged configurations, contragredient duality}

\begin{abstract}
We introduce ``virtual'' crystals of the affine types $\gggg=\Dt$,
$\At$ and $\Cn$ by naturally extending embeddings of crystals of
types $B_n$ and $C_n$ into crystals of type $A_{2n-1}$.
Conjecturally, these virtual crystals are the crystal bases of
finite dimensional $U_q'(\gggg)$-modules associated with multiples
of fundamental weights. We provide evidence and in some cases
proofs of this conjecture. Recently, fermionic formulas for the
one dimensional configuration sums associated with tensor products
of the finite dimensional $U_q'(\gggg)$-modules were conjectured
by Hatayama et al. We provide proofs of these conjectures in
specific cases by exploiting duality properties of crystals and
rigged configuration techniques. For type $\At$ we also conjecture
a new fermionic formula coming from a different labeling of the Dynkin
diagram.
\end{abstract}

\maketitle


\section{Introduction}

The quantized universal enveloping algebra $U_q(\gggg)$ associated
with a symmetrizable Kac--Moody Lie algebra $\gggg$ was introduced
independently by Drinfeld \cite{D} and Jimbo \cite{J} in their
study of two dimensional solvable lattice models in statistical
mechanics. The parameter $q$ corresponds to the temperature of the
underlying model. Kashiwara \cite{Ka1} showed that at zero
temperature or $q=0$ the representations of $U_q(\gggg)$ have
bases, which he coined crystal bases, with a beautiful
combinatorial structure and favorable properties such as
uniqueness and stability under tensor products.

\subsection{Finite affine crystals}
The underlying algebras $\gggg$ of affine crystals are
affine Kac--Moody algebras. There are two main categories of
affine crystals: (1) crystal bases of infinite dimensional integrable
highest weight $U_q(\gggg)$-modules, and (2) crystal bases
of finite dimensional $U'_q(\gggg)$-modules, where $U'_q(\gggg)$
is a subalgebra of $U_q(\gggg)$ without the degree operator. The
former category is well-understood. For instance it is known
\cite{Ka1h} that an irreducible integrable $U_q(\gggg)$-module has
a unique crystal basis. On the other hand, the latter is still far
from well-understood. It is not even known which finite
dimensional $U'_q(\gggg)$-modules have a crystal basis.

In \cite{HKOTY,HKOTT} it was conjectured that there exists a
family of finite dimensional $U'_q(\gggg)$-modules $\{W^{(r)}_s\}$
having crystal bases $\{B^{r,s}\}$ ($1\le r\le n,s\ge1$), where
$n+1$ is the number of vertices in the Dynkin diagram of $\gggg$.
If such crystals indeed exist, one can define one dimensional
configuration sums, which play an important role in the study of
phase transitions in two dimensional exactly solvable lattice
models. If $\gggg$ is of type $A^{(1)}_n$, the existence of the
crystal $B^{r,s}$ is settled in \cite{KKMMNN2}, and the one
dimensional configuration sums contain the Kostka polynomials
\cite{LLT,NY,SW,S3}, which turn up in the theory of symmetric
functions \cite{Mac}, combinatorics \cite{LS1}, the study of
subgroups of finite abelian groups \cite{But}, and
Kazhdan--Lusztig theory \cite{KL,L1}. In certain limits they are
branching functions of integrable highest weight modules
\cite{HKKOTY,NY}. For all affine types and $s=1$ it is shown in
\cite{Ka2} that $B^{r,1}$ exists (but not that it has the desired
decomposition as a classical crystal in all cases). For type
$\At$, the structure of $B^{r,1}$ has been described explicitly in
\cite{JMO}. The existence and structure of $B^{r,s}$ has been
proven for various other subcases \cite{KKM,KKMMNN2,Ko,Ko2,Y}.

\subsection{Virtual crystals}
Classical finite dimensional crystals of type $C_n$ and $B_n$ can
be embedded into crystals of type $A_{2n-1}$ \cite{B}. This
embedding has been extended to certain perfect affine crystals
\cite{Ka3}. Following this idea we consider ``virtual'' crystals
$\V^{r,s}$ for types $\Dt$, $\At$, $\Atd$ (the Dynkin diagram
$\At$ with different distinguished vertex $0$), and $\Cn$ by
naturally extending these embeddings to finite affine crystals. We
conjecture that these virtual crystals $\V^{r,s}$ are exactly the
conjectured crystals $B^{r,s}$ of \cite{HKOTY,HKOTT}. This is
proven for $s=1$ and types $\At$, $\Atd$, and $\Cn$. For type
$\Atd$ the work was already done by combining \cite{B} and
\cite{JMO}. As further evidence we show that the virtual crystals
have the expected decomposition into classical components for all
$s\ge 1$ in type $\Dt$. In all cases the virtual crystal
$\V^{r,s}$ is shown to have at least the desired classical
components. We give an explicit characterization of elements in
$\V^{r,1}$, which is used to establish the equality between the
one-dimensional sums (defined in terms of the virtual crystals
$\V^{r,1}$) and the fermionic formulas mentioned below.

\subsection{Fermionic formulas}

Crystal theory provides a general setting for the definition of
one dimensional configuration sums, denoted by $X$. Explicit
formulas $M$, called fermionic formulas, for $X$ have recently
been conjectured in \cite{HKOTY,HKOTT}. Fermionic formulas are
expressions for certain polynomials in $q$ (or $q$-series) which
are sums of products of $q$-binomial coefficients. These
expressions play a preponderant role in Bethe Ansatz studies of
spin chain models. At $q=1$ they give the number of solutions to
Bethe equations of the underlying model. Kirillov and Reshetikhin
\cite{KR} obtained fermionic expressions for the Kostka
polynomials from the Bethe Ansatz study of the $XXX$ model. In a
separate development Kedem and McCoy \cite{KM} discovered
fermionic expressions for the branching functions of the coset
pair $(A_3^{(1)})_1 \oplus (A_3^{(1)})_1 / (A_3^{(1)})_2$ based on
the Bethe Ansatz study of the 3-state Potts chain. This led to a
flurry of works on fermionic expressions for characters and
branching functions of conformal field theories (see \cite{HKOTT}
for references). At present many of these fermionic expressions
are understood as certain limits of $M$ or its level truncation.

Branching functions and one dimensional configuration sums are,
however, not the only places where fermionic formulas have
appeared. Kirillov and Reshetikhin \cite{KR1} conjectured that for
$\gggg$ of nontwisted affine type, the coefficients of the
decomposition of the representations of $U'_q(\gggg)$ naturally
associated with multiples of the fundamental weights into direct
sums of irreducible representations of $U_q(\gfin)$ are given by
the fermionic formulas at $q=1$. Here $\gfin$ is the simple Lie
algebra associated to the Kac--Moody algebra $\gggg$. Chari
\cite{Chari} proved a number of cases of this conjecture.
Recently, it was conjectured by Lusztig \cite{L} and Nakajima
\cite{Na} that quiver varieties are also related to fermionic
formulas.

The conjecture $X=M$ was proven in \cite{KSS} for type
$A_n^{(1)}$. The main technique used in the proof is a statistic
preserving bijection between crystals and rigged configurations.
Rigged configurations are combinatorial objects which label all
terms in the fermionic formula and also originate from the Bethe
Ansatz \cite{KKR,KR}. Here we prove the conjecture $X=M$ 
of \cite{HKOTY,HKOTT} for types $\Dt$, $\At$ and $\Cn$ for tensor 
products of crystals of the form $B^{r,1}$, assuming that 
$\V^{r,1}=B^{r,1}$ for type $\Dt$.
We also conjecture a new fermionic formula of type $\Atd$.

\subsection{Duality on rigged configurations}

The proof of $X=M$ is achieved by an explicit characterization of the
rigged configurations in the image of the embedding of the above
crystals into crystals of type $A_{2n-1}^{(1)}$ combined with the
bijection between type $A_{2n-1}^{(1)}$ crystals and rigged
configurations \cite{KSS}. An important tool that is used in
the proofs and the characterizations is the contragredient duality
on crystals and rigged configurations.

\subsection{Content of the paper}

Section \ref{sec:crystal} gives a review of crystal theory, in
particular the definition of crystal graphs, dual crystals, simple
crystals, the combinatorial $R$-matrix, the (intrinsic) energy
function and the definition of the one dimensional configuration
sums. In Section \ref{sec:type A} the combinatorial structure of
crystals of type $A_n^{(1)}$ is discussed in detail, and all
results and definitions of Section \ref{sec:crystal} are stated in
terms of operations on Young tableaux. Irreducible components of
tensor products of crystals of type $A_n^{(1)}$ can be
parametrized in terms of Littlewood--Richardson tableaux. These
are introduced in section \ref{sec:LR}. They are of importance in
section \ref{sec:RC} where we define rigged configurations and
recall the bijection between Littlewood--Richardson tableaux and
rigged configurations of \cite{KSS}. The relation between the
duality and this bijection is given in Theorem \ref{cdual RC}.
Virtual crystals $V^{r,s}$ for types $\Dt$, $\At$, $\Atd$, and
$\Cn$ are introduced in section \ref{sec:three types}. Their
explicit characterization for $s=1$ is given in section
\ref{subsec:char} in terms of classical virtual crystals
$\V(\la)$. Finally, in section \ref{sec:fermi} the fermionic
formulas for $s=1$ and types $\Dt$, $\At$ and $\Cn$ are proved.
The proof relies on the main Theorem \ref{thm rc prop} which gives
a characterization of the rigged configurations in the image of
the virtual crystals under the bijection of section \ref{sec:RC}.
For type $\Atd$ a similar characterization is given in 
Conjecture \ref{conj:rc Atd} which yields a new fermionic
formula given in Section \ref{sec:fermi Atd}.
Appendix \ref{app proof} is reserved for the proof of Theorem
\ref{thm rc prop}.

Finally, it worth mentioning that the techniques
of this paper have been applied in \cite{OSS} to obtain fermionic
formulas for level-restricted one dimensional configuration sums
of type $\Cn$.

\subsection*{Acknowledgments}

The self-duality property of rigged configurations of type
$C_n^{(1)}$ as stated in Theorem \ref{thm rc prop} was first
observed by Tim Baker. The authors thank him for this information.

\section{Crystal Theory} \label{sec:crystal}
In this section some results from affine crystal theory are
reviewed. In addition we show that the tensor category of simple
crystal bases of irreducible finite-dimensional modules over
quantized affine algebras, can be enhanced by including a grading
function which we call the intrinsic energy. This result
(Proposition \ref{pp:graded tensor}) is central to the definition
of the energy function of a virtual crystal, notions that are
introduced in Section \ref{sec:three types}. A prerequisite for
the explicit formula for the intrinsic energy of a tensor product,
is the discussion in Section \ref{subsec:iew} of isomorphisms and
energy functions indexed by reduced words.

\subsection{Crystals} \label{subsec:crystals}
Let $\gggg$ be a symmetrizable Kac-Moody algebra, $P$ the weight
lattice, $I$ the index set for the vertices of the Dynkin diagram
of $\gggg$, $\{\alpha_i\in P \mid i\in I \}$ the simple roots, and
$\{h_i\in P^*=\Hom_{\Z}(P,\Z)\mid i\in I \}$ the simple coroots.
Let $U_q(\gggg)$ be the quantized universal enveloping algebra of
$\gggg$ \cite{Ka}. A $U_q(\gggg)$-crystal is a nonempty set $B$
equipped with maps $\wt:B\rightarrow P$ and
$\et_i,\ft_i:B\rightarrow B\cup\{\emptyset\}$ for all $i\in I$,
satisfying
\begin{align}
\ft_i(b)=b' &\Leftrightarrow \et_i(b')=b
\text{ if $b,b'\in B$} \\
\wt(\ft_i(b))&=\wt(b)-\alpha_i \text{ if $\ft_i(b)\in B$} \\
\label{eq:string length}
\inner{h_i}{\wt(b)}&=\varphi_i(b)-\epsilon_i(b).
\end{align}
Here for $b \in B$
\begin{equation} \label{eq:eps phi def}
\begin{split}
\epsilon_i(b)&= \max\{n\ge0\mid \et_i^n(b)\not=\emptyset \} \\
\varphi_i(b) &= \max\{n\ge0\mid \ft_i^n(b)\not=\emptyset \}.
\end{split}
\end{equation}
(It is assumed that $\varphi_i(b),\epsilon_i(b)<\infty$ for all
$i\in I$ and $b\in B$.) A $U_q(\gggg)$-crystal $B$ can be viewed
as a directed edge-colored graph (the crystal graph) whose
vertices are the elements of $B$, with a directed edge from $b$ to
$b'$ labeled $i\in I$, if and only if $\ft_i(b)=b'$.

\subsection{Morphisms of crystals}
A morphism $\psi:B_1\rightarrow B_2$ of $U_q(\gggg)$-crystals is a
map $\psi:B_1\cup\{\emptyset\}\rightarrow B_2\cup\{\emptyset\}$
such that \begin{align} &\psi(\emptyset)=\emptyset. \\
&\text{If $\psi(b)\not=\emptyset$ for $b\in B_1$ then} \\
\notag &\wt(\psi(b))=\wt(b), \\
\notag &\epsilon_i(\psi(b))=\epsilon_i(b), \text{ and } \\
\notag &\varphi_i(\psi(b))=\varphi_i(b). \\
&\text{For $b\in B_1$ such that $\psi(b)\not=\emptyset$ and $\psi(\et_i(b))\not=\emptyset$,} \\
\notag &\psi(\et_i(b))=\et_i(\psi(b)). \\
&\text{For $b\in B_1$ such that $\psi(b)\not=\emptyset$ and $\psi(\ft_i(b))\not=\emptyset$,} \\
\notag &\psi(\ft_i(b))=\ft_i(\psi(b)).
\end{align}

\subsection{Tensor products of crystals} \label{sec:tensor
crystal} Let $B_1,B_2,\dotsc,B_L$ be $U_q(\gggg)$-crystals. The
Cartesian product $B_L\times \dotsm \times B_2 \times B_1$ has the
structure of a $U_q(\gggg)$-crystal using the so-called signature
rule \cite{KMOU}. The resulting crystal is denoted
$B=B_L\otimes\dots\otimes B_2\otimes B_1$ and its elements
$(b_L,\dotsc,b_1)$ are written $b_L\otimes \dotsm \otimes b_1$
where $b_j\in B_j$. The reader is warned that our convention is
opposite to that of Kashiwara \cite{Ka}. Fix $i\in I$ and
$b=b_L\otimes\dotsm\otimes b_1\in B$. The $i$-signature of $b$ is
the word consisting of the symbols $+$ and $-$ given by
\begin{equation*}
\underset{\text{$\varphi_i(b_L)$ times}}{\underbrace{-\dotsm-}}
\quad \underset{\text{$\epsilon_i(b_L)$
times}}{\underbrace{+\dotsm+}} \,\dotsm\,
\underset{\text{$\varphi_i(b_1)$ times}}{\underbrace{-\dotsm-}}
\quad \underset{\text{$\epsilon_i(b_1)$
times}}{\underbrace{+\dotsm+}} .
\end{equation*}
The reduced $i$-signature of $b$ is the subword of the
$i$-signature of $b$, given by the repeated removal of adjacent
symbols $+-$ (in that order); it has the form
\begin{equation*}
\underset{\text{$\varphi$ times}}{\underbrace{-\dotsm-}} \quad
\underset{\text{$\epsilon$ times}}{\underbrace{+\dotsm+}}.
\end{equation*}
If $\varphi=0$ then $\ft_i(b)=\emptyset$; otherwise
\begin{equation*}
\ft_i(b_L\otimes\dotsm\otimes b_1)= b_L\otimes \dotsm \otimes
b_{j+1} \otimes \ft_i(b_j)\otimes \dots \otimes b_1
\end{equation*}
where the rightmost symbol $-$ in the reduced $i$-signature of
$b$ comes from $b_j$. Similarly, if $\epsilon=0$ then
$\et_i(b)=\emptyset$; otherwise
\begin{equation*}
\et_i(b_L\otimes\dotsm\otimes b_1)= b_L\otimes \dotsm \otimes
b_{j+1} \otimes \et_i(b_j)\otimes \dots \otimes b_1
\end{equation*}
where the leftmost symbol $+$ in the reduced $i$-signature of $b$
comes from $b_j$. It is not hard to verify that this well-defines
the structure of a $U_q(\gggg)$-crystal with
$\varphi_i(b)=\varphi$ and $\epsilon_i(b)=\epsilon$ in the above
notation, with weight function
\begin{equation} \label{eq:tensor wt}
\wt(b_L\otimes\dotsm\otimes b_1)=\sum_{j=1}^L \wt(b_j).
\end{equation}
This tensor construction is easily seen to be associative. The
case of two tensor factors is given explicitly by
\begin{equation} \label{eq:f on two factors}
\ft_i(b_2\otimes b_1) = \begin{cases} \ft_i(b_2)\otimes b_1
& \text{if $\epsilon_i(b_2)\ge \varphi_i(b_1)$} \\
b_2\otimes \ft_i(b_1) & \text{if $\epsilon_i(b_2)<\varphi_i(b_1)$}
\end{cases}
\end{equation}
and
\begin{equation} \label{eq:e on two factors}
\et_i(b_2\otimes b_1) = \begin{cases} \et_i(b_2) \otimes b_1 &
\text{if $\epsilon_i(b_2)>\varphi_i(b_1)$} \\
b_2\otimes \et_i(b_1) & \text{if $\epsilon_i(b_2)\le
\varphi_i(b_1)$.}
\end{cases}
\end{equation}

\subsection{Crystals for subalgebras}
For a subset $J\subset I$ define the $J$-components of a
$U_q(\gggg)$-crystal $B$ to be the connected components of the
crystal graph using only edges labeled with elements of $J$. Say
that $b\in B$ is a $J$-highest weight vector if $\epsilon_i(b)=0$
for all $j\in J$. A highest weight vector is an $I$-highest weight
vector.

\subsection{Highest weight crystals}

Let $P^+ = \{\la\in P\mid \inner{h_i}{\la} \ge0 \text{ for all
$i\in I$}\}$. For $\la\in P^+$ let $B(\la)$ be the crystal basis of the
irreducible integrable $U_q(\gggg)$-module of highest weight $\la$.
By \cite{Ka} $B(\la)$ has the structure of a $U_q(\gggg)$-crystal
in the sense of section \ref{subsec:crystals} where $\et_i$ and
$\ft_i$ are the modified Chevalley generators.
$B(\la)$ is connected and has a unique
highest weight vector, denoted $u_\la$; it is the unique element
in $B(\la)$ of weight $\la$.

\subsection{Dual crystals}
\label{sec:gen dual} The notion of a dual crystal is given in
\cite[Section 7.4]{Ka}. Let $B$ be a $U_q(\gggg)$-crystal. Then
there is a $U_q(\gggg)$-crystal denoted $B^\cd$ obtained from $B$
by reversing arrows. That is, $B^\cd=\{b^\cd\mid b\in B\}$ with
\begin{equation} \label{eq:dual crystal}
\begin{split}
\wt(b^\cd)&=-\wt(b) \\
\epsilon_i(b^\cd)&=\varphi_i(b) \\
\varphi_i(b^\cd)&=\epsilon_i(b) \\
\et_i(b^\cd) &= (\ft_i(b))^\cd \\
\ft_i(b^\cd) &= (\et_i(b))^\cd.
\end{split}
\end{equation}

\begin{prop} \label{pp:dual tensor} \cite{Ka} There is an isomorphism
$(B_2\otimes B_1)^\cd \cong B_1^\cd \otimes B_2^\cd$ given by
$(b_2\otimes b_1)^\cd \mapsto b_1^\cd \otimes b_2^\cd$.
\end{prop}

\subsection{Extremal vectors}
\label{sec:extremal vector} Let $W$ be the Weyl group of $\gggg$,
$\{s_i\mid i\in I\}$ the simple reflections in $W$. Let $B$ be the
crystal graph of an integrable $U_q(\gggg)$-module. Say that $b\in
B$ is an extremal vector of weight $\la\in P$ provided that
$\wt(b)=\la$ and there exists a family of elements $\{b_w\mid w\in
W \}\subset B$ such that
\begin{enumerate}
\item $b_w=b$ for $w=e$.
\item If $\inner{h_i}{w\la}\ge0$ then $\et_i(b)=\emptyset$ and
$\ft_i^{\inner{h_i}{w\la}}(b_w)=b_{s_i w}$.
\item If $\inner{h_i}{w\la}\le0$ then $\ft_i(b)=\emptyset$ and
$\et_i^{\inner{h_i}{w\la}}(b_w)=b_{s_i w}$.
\end{enumerate}

\subsection{Affine case} \label{subsec:affine case}
We now fix notation for the affine case, which, for the most part,
follows \cite{Kac}. Let $\gggg$ be an affine Kac-Moody algebra
over $\CC$. Let $I$ be the index set for the Dynkin diagram, $P$
the weight lattice, $P^*=\Hom_{\Z}(P,\Z)$, $\inner{\cdot}{\cdot}$
the natural pairing $P^* \otimes P\rightarrow\Z$, $\{\alpha_i\mid
i\in I\}\subset P$ the simple roots, $\{h_i\mid i\in I\}\subset
P^*$ the simple coroots, $\{\La_i\mid i\in I\}\subset P$ the
fundamental weights, $c$ the canonical central element, and
$\delta\in P$ the null root. One has $P=\Z\delta\oplus
\bigoplus_{i\in I} \Z \La_i$.

Let the distinguished vertex in the Dynkin diagram be labeled
$0\in I$ as specified in \cite{Kac}. Let $\gfin\subset \gggg$ be
the simple Lie algebra whose Dynkin diagram is $J=I-\{0\}$ and
$\gggg'\subset\gggg$ the derived subalgebra. Denote the weight
lattices of $\gfin\subset\gggg'\subset\gggg$ by $\Pfin$, $\Pcl$,
and $P$ respectively. There are natural projections
\begin{equation*}
   P \overset{\cl}{\twoheadrightarrow}\Pcl\twoheadrightarrow\Pfin.
\end{equation*}
$\Pcl \cong P/\Z\delta$ and $\Pfin\cong
P/(\Z\delta\oplus\Z\La_0)$. We shall identify $\Pcl$ and $\Pfin$
with the sublattices of $P$ given by $\bigoplus_{i\in I} \Z\La_i$
and $\bigoplus_{i\in J} \Z\Lab_i $ where
$\Lab_i=\La_i-\inner{c}{\La_i}\La_0$. Let $\Pfin^+=\bigoplus_{i\in
J} \N \Lab_i$.

Denote by $U_q(\gfin)$, $U'_q(\gggg)$, and $U_q(\gggg)$ the
quantized universal enveloping algebras of $\gfin$, $\gggg'$, and
$\gggg$ respectively \cite{Ka}. Let $W$ and $\Wb$ be the Weyl
groups of $\gggg$ and $\gfin$; they are generated by the simple
reflections $\{s_i\mid i\in I\}$ and $\{s_i\mid i\in J\}$
respectively.

\subsection{Simple crystals}
\label{sec:simple crystal} Following \cite{AK}, say that a
$U'_q(\gggg)$-crystal $B$ is \textit{simple} if
\begin{enumerate}
\item $B$ is the crystal basis of a finite dimensional integrable
$U'_q(\gggg)$-module.
\item There is a weight $\la\in \Pfin^+$ such that $B$ has a
unique vector (denoted $u(B)$) of weight $\la$, and the weight of
any extremal vector of $B$ is contained in $\Wb\la$.
\end{enumerate}

In the definition of simple crystal in \cite{AK}, condition 1 is not
present. However we always want to assume both conditions,
so it is convenient to include condition 1 in the definition above.

\begin{theorem} \label{thm:simple} \cite{AK}
\begin{enumerate}
\item Simple crystals are connected.
\item The tensor product of simple crystals is simple.
\end{enumerate}
\end{theorem}

\begin{remark} \label{rem:simple} Suppose $\Psi:B\rightarrow B'$ is
an isomorphism of simple $U'_q(\gggg)$-crystals. Then $\Psi$ is
the only such isomorphism. Since $\Psi$ must preserve weight and
send extremal vectors to extremal vectors, it follows that
$\Psi(u(B))=u(B')$. Since $B$ is connected, the rest of the map
$\Psi$ is determined by the requirement that $\Psi$ be a
$U'_q(\gggg)$-crystal morphism.
\end{remark}

Suppose $B_j$ are simple $U'_q(\gggg)$-crystals. By Theorem
\ref{thm:simple} $B=B_L\otimes\dotsm\otimes B_1$ is simple, with
\begin{equation} \label{eq:extremal tensor}
 u(B)=u(B_L)\otimes\dotsm\otimes u(B_1).
\end{equation}

\subsection{Level and perfectness}
\label{sec:level perfect} For the $U'_q(\gggg)$-crystal $B$,
define $\epsilon,\varphi:B\rightarrow \Pcl$ by
\begin{equation}
\begin{split}
\epsilon(b)&=\sum_{i\in I} \epsilon_i(b) \La_i \\
\varphi(b) &= \sum_{i\in I} \varphi_i(b) \La_i.
\end{split}
\end{equation}
Define the level of $B$ by
\begin{equation}
\lev(B)=\min \{\inner{c}{\epsilon(b)}\mid b\in B\}.
\end{equation}
Say that $b\in B$ is minimal if $\inner{c}{\epsilon(b)}=\lev(B)$.

Let $B$ be simple. Say that $B$ is perfect of level $\ell$ if
$\ell=\lev(B)$ and $\epsilon$ and $\varphi$ are bijections from
the set of minimal elements of $B$ to
$(\Pcl^+)_\ell=\{\la\in\Pcl\mid \inner{c}{\la}=\ell \}$.

\subsection{Combinatorial $R$-matrix}
\label{sec:comb R def} The following two results essentially
appear in \cite{KKMMNN} in the case that $B_1=B_2=B_3$. They are a
consequence of the existence of a combinatorial analogue of the
$R$-matrix on the affinizations of finite $U'_q(\gggg)$-crystals
and the fact that it satisfies the Yang-Baxter equation.

\begin{theorem} \label{thm:two tensor} \cite[Section 4]{KKMMNN}
Suppose $B_1$ and $B_2$ are simple $U'_q(\gggg)$-crystals. Then
\begin{enumerate}
\item There is a unique isomorphism of $U'_q(\gggg)$-crystals
$\R=\R_{B_2,B_1}:B_2\otimes B_1\rightarrow B_1\otimes B_2$.
\item There is a function $H=H_{B_2,B_1}:B_2\otimes
B_1\rightarrow\Z$, unique up to global additive constant, such
that $H$ is constant on $J$-components and, for all $b_2\in B_2$
and $b_1\in B_1$, if $\R(b_2\otimes b_1)=b_1'\otimes b_2'$, then
\begin{equation} \label{eq:local energy}
  H(\et_0(b_2\otimes b_1))=
  H(b_2\otimes b_1)+
  \begin{cases}
    -1 & \text{if $\epsilon_0(b_2)>\varphi_0(b_1)$ and
    $\epsilon_0(b_1')>\varphi_0(b_2')$} \\
    1 & \text{if $\epsilon_0(b_2)\le\varphi_0(b_1)$ and
    $\epsilon_0(b_1')\le \varphi_0(b_2')$} \\
    0 & \text{otherwise.}
  \end{cases}
\end{equation}
\end{enumerate}
\end{theorem}

We shall call the maps $\R$ and $H$ the local isomorphism and
local energy function on $B_2\otimes B_1$. The pair $(\R,H)$ is
called the combinatorial $R$-matrix.

By Remark \ref{rem:simple} and \eqref{eq:extremal tensor},
\begin{equation} \label{eq:extremal R}
  \R(u(B_2)\otimes u(B_1)) = u(B_1)\otimes u(B_2).
\end{equation}
It is convenient to normalize the local energy function $H$ by
requiring that
\begin{equation} \label{eq:extremal H}
  H(u(B_2)\otimes u(B_1)) = 0.
\end{equation}
With this convention it follows by definition that
\begin{equation} \label{eq:HR=H}
H_{B_1,B_2} \circ \R_{B_2,B_1} = H_{B_2,B_1}
\end{equation}
as operators on $B_2\otimes B_1$.

The following observation, which follows from the uniqueness of
the local isomorphism and energy function, is obvious but useful.

\begin{prop}\label{pp:local R iso}
Suppose $B_1$, $B_2$, and $B_1'$ are simple $U'_q(\gggg)$-crystals
and there is a $U'_q(\gggg)$-crystal isomorphism
$\Psi:B_1\rightarrow B_1'$. Then
\begin{align}
\label{eq:R iso} \R_{B_2,B_1} &= (\Psi^{-1}\otimes 1_{B_2})\circ
\R_{B_2,B_1'}\circ (1_{B_2}\otimes \Psi) \\
\label{eq:H iso} H_{B_2,B_1} &= H_{B_2,B_1'}\circ (1_{B_2}\otimes
\Psi).
\end{align}
\end{prop}

Suppose $B_j$ is a simple $U'_q(\gggg)$-crystal for $1\le j\le L$.
By abuse of notation let $\R_k$ and $H_k$ denote the local
isomorphism and local energy function acting in the $k$-th and
$(k+1)$-st tensor positions (from the right).

\begin{theorem} \label{thm:three tensor} \cite{KKMMNN}
There is a unique $U'_q(\gggg)$-crystal isomorphism $B_3\otimes
B_2\otimes B_1\rightarrow B_1\otimes B_2\otimes B_3$ given by
either side of the Yang-Baxter equation
\begin{equation} \label{eq:iso YB}
\R_1 \circ \R_2\circ \R_1 = \R_2 \circ \R_1 \circ \R_2.
\end{equation}
Moreover
\begin{align}
\label{eq:energy YB1} H_1+H_2\R_1 &= H_1\R_2+H_2\R_1\R_2 \\
\label{eq:energy YB2} H_2+H_1\R_2 &= H_2\R_1+H_1\R_2\R_1.
\end{align}
\end{theorem}
Observe that \eqref{eq:energy YB1} and \eqref{eq:energy YB2} are
equivalent: the latter is obtained from the former by composing on
the right by $\R_1\R_2$, using \eqref{eq:iso YB}, and the
identities $H_1\R_1=H_1$ and $H_2\R_2=H_2$ which hold by
\eqref{eq:HR=H}.

\subsection{Isomorphisms and energy functions indexed by reduced
words} \label{subsec:iew}
Let $B_j$ be a simple $U'_q(\gggg)$-module for $1\le j\le L$
and let $B=B_L\otimes\dotsm\otimes B_1$. Let
$a=(a_1,a_2,\dotsm,a_p)$ be a sequence of indices in the set
$\{1,2,\dots,L-1\}$. Let $s_j$ denote the permutation that
exchanges the $j$-th and $(j+1)$-st positions. Define $\perm(a)$,
the permutation of the set $\{1,2,\dots,L\}$ associated with the
sequence $a$, by
\begin{equation}\label{eq:word perm}
\perm(a)=s_{a_1}s_{a_2}\dotsm s_{a_p}.
\end{equation}
Denote by $\R_a:B\rightarrow \perm(a) B$ the
$U'_q(\gggg)$-crystal isomorphism given by
\begin{equation} \label{eq:R word def}
\R_a = \R_{a_1}\R_{a_2}\dotsm\R_{a_p}
\end{equation}
and the energy function $E_a:B\rightarrow\Z$ by
\begin{equation} \label{eq:H word def}
  E_a = \sum_{j=1}^p H_{a_j} \R_{a_{j+1}}\R_{a_{j+2}}\dotsm\R_{a_p}.
\end{equation}

\begin{prop} \label{rem:word iso} With the above notation, let
$B'$ be a tensor product obtained from $B$ by reordering the
factors. Then there is a unique $U'_q(\gggg)$-crystal isomorphism
$\R:B\rightarrow B'$, and $\R=\R_a$ for any sequence $a$ such that
$B'=\perm(a)B$.
\end{prop}
\begin{proof} By Theorem \ref{thm:simple}, $B$ and $B'$ are simple
since they are tensor products of simple crystals. Let $a$ be any
sequence such that $B'=\perm(a)B$. Then $\R_a:B\rightarrow B'$ is
an isomorphism of simple $U'_q(\gggg)$-crystals. But there is a
unique $U'_q(\gggg)$-crystal isomorphism $B\rightarrow B'$ by
Remark \ref{rem:simple}.
\end{proof}

\begin{remark} \label{rem:exchange energy}
If $a'$ is obtained from $a$ by a sequence of exchanges of the
form $(\dotsc,i,j,\dotsc)\rightarrow(\dotsc,j,i,\dotsc)$ where
$|i-j|>1$, then it is easy to see that $E_a=E_{a'}$.
\end{remark}

\begin{remark} In general $E_{(1,2,1)}\not=E_{(2,1,2)}$. Take
$B_1=B_2=B_3=B^{1,1}$ in $A^{(1)}_1$ with $b=1\otimes 2\otimes 1$
(see sections \ref{sec:crystal A} and \ref{sec:A perf}). Then
$E_{(1,2,1)}(b)=1$ and $E_{(2,1,2)}=2$.
\end{remark}

For later use we define a few specific reduced words. Let
\begin{equation} \label{eq:word rev}
\begin{split}
  a^{(1)}&=\emptyset \\
  a^{(k)}&=(1,2,\dots,k-1,a^{(k-1)}) \qquad\text{for $k>1$.}
\end{split}
\end{equation}
Note that $\perm(a^{(L)})$ is the permutation that reverses the
numbers $\{1,2,\dots,L\}$. For $\ell$ and $m$ positive integers
summing to $L$, define
\begin{equation} \label{eq:word weave}
\begin{split}
  a^{(\ell,1)} &= (1,2,\dotsc,\ell) \\
  a^{(\ell,m)} &= (m,m+1,\dotsc,\ell+m-1,a^{(\ell,m-1)})
  \qquad\text{for $m>1$.}
\end{split}
\end{equation}
Observe that $\perm(a^{(\ell,m)})$ is the shortest permutation
that exchanges the numbers $\{1,2,\dots,\ell\}$ past the numbers
$\{\ell+1,\ell+2,\dots,\ell+m\}$.

\subsection{Tensoring tensor products}
Let $\ell$ and $m$ be positive integers. Let $B_j$ be a simple
$U'_q(\gggg)$-crystal for $1\le j\le \ell+m$. Then
$B'_2=B_{\ell+m}\otimes\dots\otimes B_{\ell+2}\otimes B_{\ell+1}$
and $B'_1=B_\ell\otimes\dots\otimes B_1$ are simple crystals by
Theorem \ref{thm:simple}. We compute the local isomorphism
$\R_{B'_2,B'_1}$ and energy function $H_{B'_2,B'_1}$ explicitly.

\begin{prop} \label{pp:weave} With the above notation,
\begin{align}
\label{eq:weave R} \R_{B'_2,B'_1} &= \R_{a^{(\ell,m)}} \\
\label{eq:weave H} H_{B'_2,B'_1} &= E_{a^{(\ell,m)}}.
\end{align}
\end{prop}

Proposition \ref{pp:weave} is proved by induction. The first
nontrivial cases are for three tensor factors.

\begin{prop} \label{pp:RH 3}
\begin{align}
\label{eq:R 2 1} \R_{B_3,B_2\otimes B_1} &= \R_1\R_2 \\
\label{eq:H 2 1} H_{B_3,(B_2\otimes B_1)} &= H_2 + H_1\R_2 \\
\label{eq:R 1 2} \R_{(B_3\otimes B_2),B_1} &= \R_2 \R_1 \\
\label{eq:H 1 2} H_{(B_3\otimes B_2),B_1} &= H_1+H_2\R_1.
\end{align}
\end{prop}
\begin{proof} We will prove the case $\ell=2$ and $m=1$ as
the case $\ell=1$ and $m=2$ is similar. We have $\R=\R_1\R_2$. It
must be shown that $H_{B_3,B_2\otimes B_1}=H_2+H_1\R_2$. Let us
fix $b\in B_3\otimes B_2\otimes B_1$. For any isomorphism $\Psi$
that reorders $B_3\otimes B_2\otimes B_1$ by a composition of
local isomorphisms, define $\pos(\Psi)\in\{1,2,3\}$ to be the
position in the threefold tensor $\Psi(b)$ where $\et_0$ acts (see
section \ref{sec:tensor crystal}). Let $\Delta
H_2=H_2(b)-H_2(\et_0 b)$, $\Delta H_1\R_2=H_1(\R_2
b)-H_1(\R_2\et_0 b)$, $\Delta=\Delta
 H_2+\Delta H_1\R_2$. Write $\R=\R_{B_3,(B_2\otimes B_1)}$ and
let $\pos'(id)$ (resp. $\pos'(\R)$) have value $L$ (left) and $R$
(right) on the two-fold tensor product $B_3 \otimes (B_2\otimes
B_1)$ (resp. $(B_2\otimes B_1)\otimes B_3$). Then the rows of the
following table give all the possibilities:
\begin{equation*}
\begin{array}{|c|c|c||c|c|c|c|c|} \hline
\pos(id) & \pos(\R_2) & \pos(\R_1\R_2) & \pos'(id)&\pos'(\R) &
\Delta H_2 & \Delta H_1\R_2 & \Delta \\ \hline
  3 & 3 & 3 & L & L & 1 & 0 & 1 \\ \hline
  3 & 2 & 2 & L & L & 0 & 1 & 1 \\ \hline
  3 & 2 & 1 & L & R & 0 & 0 & 0 \\ \hline
  2 & 3 & 3 & R & L & 0 & 0 & 0 \\ \hline
  2 & 2 & 2 & R & L & 0 & 0 & 0 \\ \hline
  2 & 2 & 1 & R & R &-1 & 0 &-1 \\ \hline
  1 & 1 & 2 & R & L & 0 & 0 & 0 \\ \hline
  1 & 1 & 1 & R & R & 0 &-1 &-1 \\ \hline
\end{array}
\end{equation*}
The columns $\pos'(id)$, $\pos'(\R)$, and $\Delta$ agree with the
defining conditions for $H$ in Theorem \ref{thm:two tensor}.
\end{proof}

The next case is $m=1$.

\begin{prop}\label{pp:RH ell 1}
\begin{align}
\label{eq:R ell 1}
\R_{B_{\ell+1},(B_\ell\otimes\dotsm\otimes B_1)}&=\R_1\R_2\dotsm\R_\ell \\
\label{eq:H ell 1} H_{B_{\ell+1},(B_\ell\otimes\dotsm\otimes B_1)}
&= \sum_{j=1}^\ell H_j \R_{j+1}\R_{j+2}\dotsm\R_\ell.
\end{align}
\end{prop}
\begin{proof} Formula \eqref{eq:R ell 1} holds by Remark \ref{rem:simple}.
Formula \eqref{eq:H ell 1} is proved by induction on $\ell$, by
applying \eqref{eq:R 2 1} and \eqref{eq:H 2 1} for $B'_3 \otimes
(B'_2\otimes B'_1)$ where $B'_3=B_{\ell+1}$, $B'_2=B_\ell$, and
$B'_1=B_{\ell-1}\otimes\dotsm \otimes B_1$.
\end{proof}

\begin{proof}[Proof of Proposition \ref{pp:weave}] Again
\eqref{eq:weave R} follows by Remark \ref{rem:simple}. The proof
of \eqref{eq:weave H} follows by induction on $m$. The base of the
induction ($m=1$) is given by Proposition \ref{pp:RH ell 1}. For
the induction step, \eqref{eq:R 1 2} and \eqref{eq:H 1 2} are
applied for $(B'_3\otimes B'_2)\otimes B'_1$ where
$B'_3=B_{\ell+m}$, $B'_2=B_{\ell+m-1}\otimes\dots\otimes
B_{\ell+1}$ and $B'_1=B_\ell\otimes \dotsm \otimes B_1$.
\end{proof}

\subsection{Graded simple crystals} \label{subsec:graded simple}
Let $B$ be a simple $U'_q(\gggg)$-crystal, equipped with a
function $D=D_B:B\rightarrow\Z$ called its intrinsic energy, which
is required to be constant on $J$-components and defined up to a
global additive constant. Call the pair $(B,D)$ a graded simple
$U'_q(\gggg)$-crystal. We normalize the intrinsic energy function
by the requirement that
\begin{equation}\label{eq:extremal D}
D_B(u(B))=0.
\end{equation}

Let $(B_j,D_j)$ be graded simple $U'_q(\gggg)$-crystals for
$j=1,2$. Write $\pi$ for the projection onto the rightmost tensor
factor. Define $D_{B_2\otimes B_1}:B_2\otimes B_1\rightarrow\Z$ by
\begin{equation} \label{eq:intrinsic 2}
D_{B_2\otimes B_1} = H_{B_2,B_1} + D_1 \pi + D_2 \pi \R_{B_2,B_1}.
\end{equation}
{}From now on we suppress the operator $\pi$, observing that the
operator $D_j$ always acts on the rightmost tensor factor, which
will be in $B_j$. The subscripts of $H$ and $\R$ indicate
positions.

Observe that \eqref{eq:extremal D} holds for $B=B_2\otimes B_1$ if
it is assumed that it holds for $B_1$ and $B_2$ and that
\eqref{eq:extremal H} holds.
This can be seen using \eqref{eq:extremal tensor} and \eqref{eq:extremal R}.

\begin{prop} \label{pp:graded tensor} Graded simple crystals form
a tensor category.
\end{prop}
\begin{proof} It must be shown that the tensor product
construction for graded simple crystals is associative. This is
true if the grading is ignored. It suffices to show that
\begin{equation*}
D_{B_3\otimes(B_2\otimes B_1)} = D_{(B_3\otimes B_2)\otimes B_1}.
\end{equation*}
 By Proposition \ref{pp:RH 3} and
\eqref{eq:intrinsic 2} we have
\begin{equation}\label{eq:D 1 2}
\begin{split}
D_{B_3 \otimes (B_2\otimes B_1)} &= H_{B_3,B_2\otimes B_1} +
D_{B_2\otimes B_1} + D_3 \R_{B_3,B_2\otimes B_1} \\
&= H_2+H_1\R_2+H_1+D_1+D_2\R_1+D_3 \R_1\R_2
\end{split}
\end{equation}
and
\begin{equation}\label{eq:D 2 1}
\begin{split}
D_{(B_3 \otimes B_2)\otimes B_1} &= H_{(B_3\otimes B_2),B_1} +
D_1 + D_{B_3\otimes B_2} \R_{B_3\otimes B_2,B_1} \\
&= H_1+H_2\R_1+D_1+(H_1+D_2+D_3 \R_1)\R_2\R_1 \\
&= H_1+H_2\R_1+H_1\R_2\R_1+D_1+D_2\R_2\R_1+D_3\R_2\R_1\R_2\\
&= H_1+H_2\R_1+H_1\R_2\R_1+D_1+D_2\R_1+D_3\R_1\R_2.
\end{split}
\end{equation}
In the second computation, \eqref{eq:iso YB} is used. Also used is
$D_j \R_2=D_j$, which holds because $D_j$ acts on the rightmost
tensor factor, which is not changed by $\R_2$. One sees that
\eqref{eq:D 1 2} and \eqref{eq:D 2 1} are equal by
\eqref{eq:energy YB2}.
\end{proof}

Let $(B_j,D_j)$ be graded simple $U'_q(\gggg)$-crystals for $1\le
j\le L$ and let $u_j=u(B_j)$. Let $B=B_L\otimes\dotsm\otimes B_1$.
Following \cite{NY} define the energy function
$E_B:B\rightarrow\Z$ by
\begin{equation}\label{eq:NY}
E_B = \sum_{1\le i<j\le L} H_i\R_{i+1}\R_{i+2}\dotsm\R_{j-1}.
\end{equation}
By the normalization assumption \eqref{eq:extremal H} it follows
that
\begin{equation}\label{eq:extremal NY}
E_B(u(B))=0.
\end{equation}

The following formula was motivated by the definition of the $D$
energy function in \cite[Section 3.3]{HKOTT}.

\begin{prop} \label{pp:intrinsic gen}
The intrinsic energy $D_B$ for the $L$-fold tensor product
$B=B_L\otimes\dotsm\otimes B_1$ is given by
\begin{equation}\label{eq:energy}
  D_B = E_B + \sum_{j=1}^L D_j \R_1\R_2\dotsm\R_{j-1}.
\end{equation}
\end{prop}
\begin{proof} The proof proceeds by induction on $L$. For $L\le 2$
\eqref{eq:energy} holds by definition. Suppose it holds for $L-1$
tensor factors. Let $B'_2=B_L$ and
$B'_1=B_{L-1}\otimes\dotsm\otimes B_1$. Applying the definition
\eqref{eq:intrinsic 2} for $B'_2\otimes B'_1$, induction, and
Proposition \ref{pp:weave}, we have
\begin{equation*}
\begin{split}
D_B &= H_{B'_2,B'_1} + D_{B'_1} + D_{B'_2} \R_{B'_2,B'_1} \\
&= \sum_{i=1}^{L-1} H_i \R_{i+1}\R_{i+2}\dotsm\R_{L-1} + E_{B'_1} \\
&+ \sum_{j=1}^{L-1} D_{B_j} \R_1\R_2\dotsm\R_{j-1} + D_{B_L}
\R_1\dotsm\R_{L-1}
\end{split}
\end{equation*}
which is evidently the right hand side of \eqref{eq:energy}.
\end{proof}

The energy functions are insensitive to reordering of tensor
factors.

\begin{prop} \label{pp:reorder DE}
Let $B'$ be a reordering of $B=B_L\otimes\dots\otimes B_1$
and $\R:B\rightarrow B'$ any composition of local isomorphisms.
Then $D_{B'} \R = D_B$ and $E_{B'}\R=E_B$.
\end{prop}
\begin{proof} The proof immediately reduces to the case $\R=\R_j$.
Let us prove the first assertion, as the second follows from it.
Write $B'_3=B_L\otimes\dotsm\otimes B_{j+2}$, $B'_1=B_{j-1}\otimes
\dotsm\otimes B_1$, $B'_2=B_{j+1}\otimes B_j$ and
$B''_2=B_j\otimes B_{j+1}$. By \eqref{eq:HR=H} it follows that
$D_{B''_2} \circ \R'=D_{B'_2}$ where $\R':B_{j+1}\otimes
B_j\rightarrow B_j\otimes B_{j+1}$. Therefore $B'_2\cong B''_2$ as
graded simple crystals. Tensoring both on the left by $B'_3$ and
on the right by $B'_1$, the result is a pair of isomorphic graded
simple crystals $B\cong B'$ via the map $\R=1_{B'_3}\otimes
\R'\otimes 1_{B'_1}$. In particular $D_B=D_{B'}\circ\R$. This
argument implicitly uses Proposition \ref{pp:graded tensor} and
Remark \ref{rem:simple}.
\end{proof}

\subsection{The crystals $B^{r,s}$}
\label{sec:finite crystals}

We recall conjectures regarding the crystals $B^{r,s}$. Recall
Chevalley's partial order on $\Pfin$, defined by
$\mu\trianglerighteq\la$ if and only if
$\mu-\la\in\bigoplus_{i\in J} \N \alpha_i$.

\begin{conjecture} \label{conj:Brs} \cite{HKOTT} For each
$r\in J$ and $s\ge1$, there exists an irreducible
finite dimensional integrable $U'_q(\gggg)$-module $W^{(r)}_s$
with simple crystal basis $B^{r,s}$ having a unique extremal vector
$u(B^{r,s})$ of weight $s\Lab_r$, and a prescribed
$U_q(\gfin)$-crystal decomposition of the form $B^{r,s} \cong
B(s\Lab_r)\oplus B$, where $B$ is a direct sum of
$U_q(\gfin)$-crystals of the form $B(\la)$ where $\la\in \Pfin^+$
and $s \Lab_r \vartriangleright \la$. Moreover there is a
prescribed intrinsic energy function
$D=D_{B^{r,s}}:B^{r,s}\rightarrow \Z$, that is constant on
$J$-components, such that $0=D(u(B^{r,s}))>D(b)$ and $b$ is any
element not in the $J$-component of $u(B^{r,s})$.
\end{conjecture}

\begin{conjecture} \cite{Ka2} $W^{(r)}_s$ is the universal
$U'_q(\gggg)$-module generated by an extremal vector of weight $s
\Lab_r$.
\end{conjecture}

This is known to hold for $s=1$ \cite{Ka2}.

\begin{conjecture} \label{conj:grade} \cite{HKOTT} There is a unique
$b^\natural \in B^{r,s}$ such that
$\varphi(b^\natural)=\lev(B^{r,s}) \La_0$. Moreover
\begin{equation} \label{eq:grade}
  D_{B^{r,s}}(b) = H(b\otimes b^\natural)-H(u(B^{r,s})\otimes b^\natural)
\end{equation}
where $H=H_{B^{r,s},B^{r,s}}$ is the local energy function.
\end{conjecture}

\begin{remark} \label{rem:grade} Observe that the function
$B^{r,s}\rightarrow \Z$ given by $b\mapsto H(b\otimes b^\natural)$
is constant on $J$-components. By definition
$\varphi_i(b^\natural)=0$ for $i\in J$. By \eqref{eq:f on two
factors}, $\ft_i(b\otimes b^\natural)=\ft_i(b)\otimes b^\natural$
for all $b\in B^{r,s}$. The assertion follows from the fact that
$H$ is constant on $J$-components of $B^{r,s}\otimes B^{r,s}$, by
Theorem \ref{thm:two tensor}.
\end{remark}

\subsection{One dimensional sums}
Let $\la\in\Pfin$ and $B=B_L\otimes \dotsm\otimes B_1$ where
$B_j=B^{r_j,s_j}$ for $1 \le j\le L$. The set $\Path(B,\la)$ of
(classically restricted) paths in $B$ of weight $\la$ is the set
of $J$-highest weight vectors in $B$ of weight $\la$. Define the
classically restricted one dimensional sum
\begin{equation}\label{eq:onedimsum}
\X(B,\la;q)=\sum_{b\in\Path(B,\la)} q^{D_B(b)}.
\end{equation}

\section{Crystals of type $A_n^{(1)}$} \label{sec:type A}
In this section we discuss the combinatorial properties
of type $A_n^{(1)}$ in detail.

\subsection{The root system $A_{n-1}$}
\label{sec:root A} In this section let $\gfin$ be the complex
simple Lie algebra of type $A_{n-1}$. Let $J=\{1,2,\dots,n-1\}$ be
the index set for the vertices of its Dynkin diagram. Let us
realize the weight lattice $\Pfin$ explicitly inside the real span
of the hyperplane in $\RR^n$ orthogonal to the vector
$e=(1,1,\dotsc,1)\in \RR^n$. Let $\{\varepsilon_i\mid 1\le i\le n\}$ be the
standard basis of $\RR^n$. Then one may take
$\alpha_i=\varepsilon_i-\varepsilon_{i+1}$ for $i\in J$. Identifying the
weight lattice of $\mathfrak{gl}(n)$ with $\Z^n$, there is a projection
$\slwt:\Z^n\rightarrow \Pfin$ defined by $v \mapsto v -
\frac{(e,v)}{n} e$ where $(\cdot,\cdot)$ is the standard bilinear form on
$\RR^n$. In coordinates it is given by $\slwt(a_1,a_2,\dots,a_n) =
\sum_{i=1}^{n-1} (a_i-a_{i+1})\Lab_i$. In particular
$\Lab_i=\slwt(1^i,0^{n-i})$ for $i\in J$. The highest root of
$\gfin$ is $\theta=(1,0^{n-2},-1)$. If $\la$ is a partition
$(\la_1\ge\la_2\ge\dotsm\ge\la_n)\in\N^n$ then $\slwt(\la)\in
\Pfin^+$. By abuse of notation an element of $\Z^n$ is sometimes
identified with its image under $\slwt$.

\subsection{Crystal graphs of type $A_{n-1}$}
\label{sec:crystal A}
Recall that for $\la\in \Pfin^+$, $B(\la)$ is the crystal of the
irreducible integrable $U_q(\gfin)$-module of highest weight $\la$.
We now review the explicit realization of $B(\la)$ given in \cite{KN}.
The vertices of $B(\la)$ are given by the set of (semistandard) tableaux of
shape $\la$ in the alphabet $[n]=\{1,2,\dots,n\}$. The
$U_q(\gfin)$-crystal structure shall be defined by embedding
tableaux into the set of words.

The crystal graph of $B(\Lab_1)$ (where recall that $\Lab_1$ is
the image under $\slwt$ of the partition $(1,0^{n-1})$) is given
by the set $[n]$ with $\ft_i(i)=i+1$ for all $i\in J$ and all
other elements sent to $\emptyset$. One defines
$\wt(i)=\Lab_i-\Lab_{i-1}$ for all $i\in[n]$.

The set of words in the alphabet $[n]$ is the crystal graph of the
tensor algebra of $B(\Lab_1)$; its explicit structure is given by
the signature rule (see section \ref{sec:tensor crystal}). The
weight function $\wt$ on this crystal graph is given as follows.
Let the content of a word $u$ be the element
$\cont(u)=(m_1(u),m_2(u)\dots,m_n(u))\in\N^n$ where $m_j(u)$ is
the number of occurrences of the letter $j$ in $u$. Then
$\wt(u)=\slwt(\cont(u))\in \Pfin$.

We identify the tableau $t\in B(\la)$ with its column-reading
word, which is by definition $t=c_1c_2\dotsm$ where $c_j$ is the
strictly decreasing word given by the $j$-th column of $t$. This
defines an embedding of $B(\la)$ into the set of words in the
alphabet $[n]$. It is easy to see that $\et_i$ and $\ft_i$
stabilize the image of this embedding. The $U_q(\gfin)$-crystal
structure on $B(\la)$ is defined by declaring this embedding to be
a morphism of $U_q(\gfin)$-crystals. The highest weight vector
$u_\la$ of $B(\la)$ is given explicitly by the Yamanouchi tableau,
the tableau of shape $\la$ whose $i$-th row is filled with the
letter $i$ for all $i$.

\subsection{Schensted's insertion}
\label{sec:schen} Schensted \cite{Schen} gave an insertion
algorithm that associates to a word $u$ a tableau $\PSchen(u)$ of
partition shape. Given an alphabet (totally ordered set) $A$, let
$\Knuth$ be the equivalence relation on words in the alphabet $A$
defined by $u\Knuth v$ if and only if $\PSchen(u)=\PSchen(v)$. Let
$\Plac(A)$ denote the set of $\Knuth$-classes. It follows from
\cite{Schen} that $\PSchen(\PSchen(u)\PSchen(v))=\PSchen(uv)$, so
that the multiplication on words given by juxtaposition, descends
to a multiplication in $\Plac(A)$. $\Plac(A)$ is called the
plactic monoid.  Knuth \cite{K} showed that the set of relations
in $\Plac(A)$ is generated by relations of the following form,
where $x,y,z\in A$ and $u,v$ are words:
\begin{equation} \label{Knuth relations}
\begin{aligned}
  u x z y v & \Knuth u z x y v & \qquad & \text{for $x \le y < z$} \\
  u y z x v & \Knuth u y x z v & \qquad & \text{for $x<y\le z$}
\end{aligned}
\end{equation}
A column word is one that is strictly decreasing.

Let $B\subseteq A$ be a subinterval.  For a word $u$ in the
alphabet $A$, denote by $u|_B$ the subword of $u$ obtained
by erasing all of the letters not in $B$.  It follows
immediately from the relations \eqref{Knuth relations}
that if $u\Knuth v$ then $u|_B \Knuth v|_B$.

It is a well-known and remarkable fact that Schensted's algorithm
gives a morphism of type $A_{n-1}$ crystal graphs.

\begin{theorem} \label{thm:crystal P} For all words $u$
in the alphabet $[n]$ and all $i\in J$, $\PSchen(\et_i(u))=\et_i
\PSchen(u)$.
\end{theorem}
\begin{proof} It suffices to show that
if $u\Knuth v$ is one of the relations \eqref{Knuth relations}
then $e_i(u)\Knuth e_i(v)$ is also. The result is entirely
straightforward unless $\{x,y,z\}=\{i,i+1\}$ and $\et_i$ changes
one of $x,y,z$. In this case (ignoring identical left and right
factors) one has $(i+1) (i+1) i \Knuth (i+1) i (i+1)$. Applying
$\et_i$ to both sides one obtains $i (i+1) i \Knuth (i+1) i i$.
\end{proof}

\subsection{The crystals $B^{r,s}$ of type $A^{(1)}_{n-1}$}
\label{sec:A perf} Let $\gggg$ be the affine algebra of type
$A^{(1)}_{n-1}$. Write $I=J\cup\{0\}=\{0,1,\dots,n-1\}$.

\begin{theorem} \label{thm:A perf} \cite{KKMMNN2} For every $r\in J$
and $s\ge 1$ there is a finite dimensional integrable
$U'_q(A^{(1)}_{n-1})$-module $W^{(r)}_s$ with crystal basis
$B^{r,s}$ such that
\begin{enumerate}
\item $B^{r,s}\cong B(s\Lab_r)$ as $U_q(A_{n-1})$-crystals.
\item $D_{B^{r,s}}=0$.
\item $B^{r,s}$ is perfect of level $s$.
\end{enumerate}
Moreover, let $\autoA$ be the rotation of the Dynkin diagram of
type $A^{(1)}_{n-1}$ that sends the simple root $\alpha_i$ to
$\alpha_{i+1}$ where subscripts are taken modulo $n$. Then there
is a unique bijection $\autoA:B^{r,s}\rightarrow B^{r,s}$ such
that
\begin{align}
\label{eq:autoA order} \autoA^n&=1 \\
\label{eq:autoA wt} \wt\circ \autoA &= \autoA \circ \wt \\
\label{eq:autoA f} \ft_i &= \autoA^{-1} \circ \ft_{i+1} \circ \autoA \\
\label{eq:autoA e} \et_i &= \autoA^{-1} \circ \et_{i+1} \circ \autoA
\end{align}
where $\autoA:\Pfin\rightarrow\Pfin$ is induced by the map
$\Z^n\rightarrow \Z^n$ given by $(a_1,a_2,\dots,a_n)\mapsto
(a_2,a_3,\dots,a_n,a_1)$.
\end{theorem}
The operators $\et_0$ and $\ft_0$ are computed explicitly in
\cite{S3}. It suffices to define $\autoA^{-1}$. For $t\in B^{r,s}$
define $\autoA^{-1}(t)\in B^{r,s}$ by the property
\begin{equation} \label{eq:autoA rotate}
  \autoA^{-1}(t)|_{[n-1]}=\PSchen(t|_{[2,n]})-1
\end{equation}
where $t+r$ is the tableau obtained by adding the integer $r$ to
each letter. This definition determines the positions in
$\autoA^{-1}(t)$ of the letters in the subinterval $[n-1]$.  But
it also determines the positions of the letters $n$ as they must
fill up the rest of the rectangular shape.

\begin{example} Let $r=3$, $s=4$, $n=5$.  Then a tableau
$t\in B^{r,s}$ and $\autoA^{-1}(t)$ are given below.
\begin{equation*}
\begin{split}
  t =
  \begin{array}{cccc} \hline
    \mcbn{1} & \mcrn{1} & \mcrn{2} &\mcrn{3} \\ \cline{1-4}
    \mcbn{2} & \mcrn{2} & \mcrn{3} &\mcrn{4} \\ \cline{1-4}
    \mcbn{3} & \mcrn{4} & \mcrn{5} &\mcrn{5} \\ \cline{1-4}
  \end{array}
\qquad
  t|_{[2,n]} =
  \begin{array}{cccc} \hline
    \mcb & \mcr & \mcrn{2} &\mcrn{3} \\ \cline{1-4}
    \mcbn{2} & \mcrn{2} & \mcrn{3} &\mcrn{4} \\ \cline{1-4}
    \mcbn{3} & \mcrn{4} & \mcrn{5} &\mcrn{5} \\ \cline{1-4}
  \end{array}
\qquad
  \PSchen(t|_{[2,n]}) =
  \begin{array}{cccc} \hline
    \mcbn{2} & \mcrn{2} & \mcrn{2} &\mcrn{3} \\ \cline{1-4}
    \mcbn{3} & \mcrn{3} & \mcrn{4} &\mcrn{5} \\ \cline{1-4}
    \mcbn{4} & \mcrn{5} & \mcr &\mcr \\ \cline{1-4}
  \end{array}
\\ \\
 \PSchen(t|_{[2,n]})-1 =
  \begin{array}{cccc} \hline
    \mcbn{1} & \mcrn{1} & \mcrn{1} &\mcrn{2} \\ \cline{1-4}
    \mcbn{2} & \mcrn{2} & \mcrn{3} &\mcrn{4} \\ \cline{1-4}
    \mcbn{3} & \mcrn{4} & \mcr &\mcr \\ \cline{1-4}
  \end{array}
\qquad \autoA^{-1}(t)=
  \begin{array}{cccc} \hline
    \mcbn{1} & \mcrn{1} & \mcrn{1} &\mcrn{2} \\ \cline{1-4}
    \mcbn{2} & \mcrn{2} & \mcrn{3} &\mcrn{4} \\ \cline{1-4}
    \mcbn{3} & \mcrn{4} & \mcrn{5} &\mcrn{5} \\ \cline{1-4}
  \end{array}
\end{split}
\end{equation*}
\end{example}

Since $B^{r,s}$ is perfect of level $s$, there is a unique
$b^\natural\in B^{r,s}$ such that $\varphi(b^\natural)=s\La_0$; it
is the tableau in $B^{r,s}$ whose $i$-th row from the bottom
consists of $s$ copies of the value $n+1-i$. Since $B^{r,s}$
consists of a single $J$-component, Conjecture \ref{conj:grade}
holds by Remark \ref{rem:grade}.

\subsection{Inhomogeneous paths}
\label{sec:inhom paths} Let $R=(R_1,R_2,\dots,R_L)$ be a sequence
of rectangular partitions. Say that $R_j$ has $r_j$ rows and $s_j$
columns where $s_j\ge1$ and $1\le r_j\le n-1$ for all $1\le j\le
L$. Define the $U'_q(\gggg)$-crystal
\begin{equation}\label{eq:B_R}
  B_R = B^{r_L,s_L} \otimes \dotsm\otimes B^{r_2,s_2} \otimes
  B^{r_1,s_1}.
\end{equation}
The elements of $B_R$ are called inhomogeneous paths. Of course
one could compute $\et_0$ and $\ft_0$ on $B_R$ using the signature
rule and the above rule for $\et_0$ and $\ft_0$ on single tensor
factors, but it is simpler to compute them as follows. Let
$\autoA:B_R\rightarrow B_R$ be the bijection given by
$\autoA(b_L\otimes \dotsm\otimes
b_1)=\autoA(b_L)\otimes\dots\otimes\autoA(b_1)$. Then $\et_0$ and
$\ft_0$ as operators on $B_R$ may be defined by the equations
\eqref{eq:autoA e} and \eqref{eq:autoA f}.

\subsection{Simplicity of tensor products of $B^{r,s}$}
\label{sec:A simple} Consider the finite dimensional integrable
$U'_q(\gggg)$-module $W^{(r)}_s$ with crystal basis $B^{r,s}$ (see
Theorem \ref{thm:A perf}). We observe that the crystal $B^{r,s}$
is simple (see section \ref{sec:simple crystal}). Let
$\la=s\Lab_r$, with corresponding partition
$(s^r,0^{n-r})\in\Z^n$. It is not hard to check that the extremal
vectors of $B^{r,s}$ are in bijection with the orbit $\Wb\la$. In
the notation of section \ref{sec:extremal vector}, for every $w\in \Wb$,
let $b_w$ be the unique tableau in $B^{r,s}$ of content
$w(s^r,0^{n-r})\in\Z^n$. This shows that $B^{r,s}$ is simple. The
crystal $B_R$ of the previous section is a tensor product of
simple crystals. So $B_R$ is simple and hence connected.

\subsection{Dual crystals}
Let $N\ge\la_1$ and $\la^\cd=(N-\la_n,\dotsc,N-\la_2,N-\la_1)$.
For a column word $c$ in the alphabet $[n]$, let $c^\cd$ be the
column word that uses precisely those letters that are not in $c$.
For a sequence of column words $u=c_1c_2\dots c_N$, define
\begin{equation}\label{eq:dual tableau}
u^\cd=c_N^\cd\dotsm c_1^\cd.
\end{equation}

\begin{prop} \label{pp:dual irr} There is a unique $U_q(\gfin)$-crystal
isomorphism $B(\la)^\cd\rightarrow B(\la^\cd)$.
\end{prop}
\begin{proof} It is easy to see that the map $B(\la)\rightarrow
B(\la^\cd)$ given by $t\mapsto t^\cd$ (where $t^\cd$ is defined in
\eqref{eq:dual tableau}) is a well-defined bijection that
satisfies \eqref{eq:dual crystal}. For uniqueness, by definition
the map $\cd$ takes highest weight vectors to lowest weight
vectors. But $B(\la)$ has a unique highest weight vector and
$B(\la^\cd)$ has a unique lowest weight vector.
\end{proof}

\begin{prop} \label{pp:dual P}
Let $u=c_1\dotsm c_N$ where each $c_j$ is a column word. Then
$\PSchen(u^\cd)=\PSchen(u)^\cd$.
\end{prop}
\begin{proof} In the statement, $\PSchen(u)^\cd$ is computed with respect to
a factorization into $N$ column words, the last several of which may be empty.
It suffices to show that if $u,u',v,v'$ are column
words such that $uv \Knuth u'v'$ with $u'$ one letter longer than $u$,
then $v^\cd u^\cd \Knuth {v'}^\cd {u'}^\cd$.  This reduction,
together with the proof of this special case, can both be seen
by considering a jeu de taquin on skew tableaux having at most two columns.
\end{proof}

\begin{example} Let $n=6$.  One has
$uv=(42)(6521)\Knuth(642)(521)\Knuth (6421)(52)=\PSchen(uv)$ and
dually $(uv)^\cd=(43)(6531)\Knuth(643)(531)\Knuth
(6431)(53)=\PSchen(uv)^\cd$.
\end{example}

\begin{prop} \label{pp:dual perfect}
The $U'_q(\gggg)$-crystal $B^{r,s}$ has dual crystal $B^{n-r,s}$.
\end{prop}
\begin{proof} Observe that as $U_q(\gfin)$-crystals one has
$(B^{r,s})^\cd\cong B^{n-r,s}$ using $N=s$ columns. By the
uniqueness in Proposition \ref{pp:dual irr} it follows that the
crystals are also dual as $U'_q(\gggg)$-crystals by the same map.
\end{proof}

\subsection{Reversing the Dynkin diagram}
\label{sec:A rev}

\begin{theorem} \label{thm:flip}
Let $B$ be the crystal basis of a finite dimensional $U_q(\gfin)$-module
(resp. $U'_q(\gggg)$-module). Then there is an
involution $B\rightarrow B$ denoted $b\mapsto b^\ed$, such that
\begin{equation} \label{eq:flip crystal}
\begin{split}
  \ft_i(b^\ed) &= \et_{n-i}(b)^\ed \\
  \et_i(b^\ed) &= \ft_{n-i}(b)^\ed \\
  \wt(b^\ed) &= w_0 \wt(b)
\end{split}
\end{equation}
for all $i\in J$ (resp. $i\in I$ with subscripts taken modulo $n$)
where $w_0$ is the longest element of $\Wb$.
\end{theorem}
\begin{proof}
Consider the orientation-reversing automorphism of the Dynkin
diagram of type $A^{(1)}_{n-1}$ that sends $\alpha_i$ to
$\alpha_{n-i}$ for all $i$ with subscripts taken modulo $n$.  The
existence of this automorphism implies the existence of a map
${}^\ed:B\rightarrow B$ satisfying \eqref{eq:flip crystal}. Since
the Dynkin diagram automorphism is an involution, so is the
induced map ${}^\ed$.
\end{proof}

The following result is easily verified.

\begin{prop} \label{pp:flip tensor}
Suppose there are crystal graphs $B_j$ and
bijections ${}^\ed:B_j\rightarrow B_j$ satisfying \eqref{eq:flip
crystal}. Let $B=B_L\otimes \dotsm\otimes B_1$ and $B^\ed=B_1
\otimes \dots \otimes B_L$. Then the map ${}^\ed:B\rightarrow
B^\ed$ given by
\begin{equation*}
  b_L\otimes\dotsm\otimes b_1 \mapsto b_1^\ed \otimes \dotsm\otimes b_L^\ed
\end{equation*}
satisfies \eqref{eq:flip crystal}.
\end{prop}

Given a word $u$, let $u^\ed$ be the word obtained by replacing
each letter $i$ by $n+1-i$, and reversing the resulting word.
Clearly if $u$ is a column word then so is $u^\ed$. If $u=c_1
c_2\dots c_N$ where $c_j$ is a column word for all $j$, then by
definition $u^\ed=c_N^\ed \dots c_1^\ed$, which is a sequence of
column words. In particular, let $t\in B(\la)$.  It is easy to
check that $t^\ed$ is a skew tableau whose shape is given by the
180 degree rotation of $\la$.  Define $t^\ev=\PSchen(t^\ed)$.  The
tableau $t^\ev$ is called the evacuation \cite{LS} of the tableau
$t$. Observe that if $\la$ is a rectangle then $t^\ev=t^\ed$.

\begin{prop} \label{pp:flip P}\mbox{}
\begin{enumerate}
\item The map $u\mapsto u^\ed$ is an involution satisfying
\eqref{eq:flip crystal}.
\item For any word $u$, $\PSchen(u^\ed)=\PSchen(u)^\ev$.
\item The map $B(\la)\rightarrow B(\la)$ given by
$t\mapsto t^\ev$ is the unique map satisfying \eqref{eq:flip
crystal} (in which $\ed$ is replaced by $\ev$).
\end{enumerate}
\end{prop}
\begin{proof} Part 1 is straightforward.
For part 2 it suffices to show that if $u\Knuth v$ is a relation
of the form \eqref{Knuth relations} then so is $u^\ed \Knuth
v^\ed$, but this is also straightforward. Part 3 follows from part
2, the connectedness of $B(\la)$, and the fact that $\PSchen$ is a
morphism of crystal graphs.
\end{proof}

We remark that the uniqueness in part 3, together with Theorem
\ref{thm:flip}, imply that $t\mapsto t^\ev$ is an involution.

\subsection{Combining duality and Dynkin reversal}

\begin{prop} \label{pp:dual flip}
There is a bijection $B(\la)\rightarrow B(\la^\cd)$ given by $t
\mapsto t^{\ev\cd}=t^{\cd\ev}$ such that
\begin{equation} \label{eq:crystal dual flip}
\begin{split}
\ft_i(t^{\ev\cd})&=\ft_{n-i}(t)^{\ev\cd} \\
\et_i(t^{\ev\cd})&=\et_{n-i}(t)^{\ev\cd} \\
\wt(t^{\ev\cd})&=-w_0(\wt(t))
\end{split}
\end{equation}
for all $t\in B(\la)$ and $i\in J$.  Moreover if $\la$ is a
rectangle then \eqref{eq:crystal dual flip} also holds for $i=0$.
\end{prop}
\begin{proof} ${}^\cd$ and ${}^\ed$ obviously commute on
sequences of column words. Proposition \ref{pp:dual P} implies
that $t^{\ev\cd}=t^{\cd\ev}$ for all $t\in B(\la)$. Equation
\eqref{eq:crystal dual flip} follows from \eqref{eq:dual crystal}
and Theorem \ref{thm:flip}.
\end{proof}

A similar statement holds for tensor products.

\subsection{Local isomorphism}
Let $B_j=B^{r_j,s_j}$ for $j=1,2$.
The existence and uniqueness of the local isomorphism
$\R:B_2\otimes B_1\rightarrow B_1\otimes B_2$ (see section
\ref{sec:comb R def}) is guaranteed by the simplicity of $B^{r,s}$
(see section \ref{sec:A simple}). Observe that as a
$U_q(\gfin)$-crystal $B_2\otimes B_1$ is
multiplicity-free, since $B_j\cong B(s_j \Lab_j)$
is indexed by a rectangular partition
for $j=1,2$ \cite{St}.
Therefore the isomorphism $\R:B_2\otimes
B_1\rightarrow B_1\otimes B_2$ is uniquely specified by the
property that $\R(b_2\otimes b_1)$ is the unique element
$b_1'\otimes b_2'\in B_1\otimes B_2$ such that $\PSchen(b_2
b_1)=\PSchen(b_1'b_2')$. The bijection $\R$ is described quite
explicitly in \cite{S1}.
\begin{example}
Let
\begin{equation*}
b_2=\begin{array}{cc} \hline
    \mcbn{1} & \mcrn{2} \\ \cline{1-2}
    \mcbn{3} & \mcrn{3} \\ \cline{1-2}
  \end{array}
\qquad \text{and} \qquad
b_1=\begin{array}{cc} \hline
    \mcbn{1} & \mcrn{1} \\ \cline{1-2}
  \end{array}\; .
\end{equation*}
Then
\begin{equation*}
\PSchen(b_2 b_1)=\begin{array}{ccc} \hline
    \mcbn{1} & \mcrn{1} & \mcrn{1}\\ \cline{1-3}
    \mcbn{2} & \mcrn{3} \\ \cline{1-2}
    \mcbn{3} \\ \cline{1-1}
  \end{array}.
\end{equation*}
Hence
\begin{equation*}
b'_1=\begin{array}{cc} \hline
    \mcbn{1} & \mcrn{3} \\ \cline{1-2}
  \end{array}
\qquad \text{and} \qquad
b'_2=\begin{array}{cc} \hline
    \mcbn{1} & \mcrn{1} \\ \cline{1-2}
    \mcbn{2} & \mcrn{3} \\ \cline{1-2}
  \end{array}\; .
\end{equation*}
\end{example}

\begin{prop} \label{pp:R dual}
The following diagram commutes:
\begin{equation}
\begin{CD}
  B_2 \otimes B_1 @>{\R}>> B_1 \otimes B_2 \\
  @V{\cd}VV    @VV{\cd}V \\
  B_1^{\cd} \otimes B_2^\cd @>>{\R}> B_2^\cd \otimes B_1^\cd
\end{CD}
\end{equation}
\end{prop}
\begin{proof} Let $b_j\in B_j$ for $j=1,2$ and
$\R(b_2 \otimes b_1)=b_1'\otimes b_2'$.  Then
\begin{equation*}
  \PSchen((b_2 b_1)^\cd) = \PSchen(b_2 b_1)^\cd = \PSchen(b_1' b_2')^\cd
 =\PSchen((b_1' b_2')^\cd)
\end{equation*}
by the definition of $\R$ and Proposition \ref{pp:dual P}. This
suffices by the definition of $\R$.
\end{proof}

\begin{prop} \label{pp:R ev}
The following diagram commutes:
\begin{equation}
\begin{CD}
  B_2 \otimes B_1 @>{\R}>> B_1 \otimes B_2 \\
  @V{\ed}VV    @VV{\ed}V \\
  B_1 \otimes B_2 @>>{\R}> B_2 \otimes B_1
\end{CD}
\end{equation}
\end{prop}
\begin{proof} Let $b_j,b_j'$ be as before for $j=1,2$.
\begin{equation*}
  \PSchen((b_2 b_1)^\ed) = \PSchen(b_2 b_1)^\ev = \PSchen(b_1' b_2')^\ev =
  \PSchen((b_1' b_2')^\ed)
\end{equation*}
by the definition of $\R$ and Proposition \ref{pp:flip P}. This
suffices by the definition of $\R$.
\end{proof}

\subsection{Local isomorphism for single columns}
\label{sec:R sc}

In the case that $B_j=B^{r_j,1}$ for $j=1,2$ there is an explicit
construction for the local isomorphism $\R$ which will be useful
later. For $\gfin=A_{n-1}$, represent $b\in B^{r,1}$ by a column
of height $n$ with a dot at height $i$ if and only if the letter
$i$ occurs in $b$. Now $\R:B^{r_2,1}\otimes B^{r_1,1}\to B^{r_1,1}
\otimes B^{r_2,1}$ with $r_2\ge r_1$ is obtained as follows
\cite{NY}. Let $b_2\otimes b_1\in B_2\otimes B_1$ and $b_1'\otimes
b_2'=\R(b_2\otimes b_1)$.
\begin{enumerate}
\item Pick the highest dot in $b_1$. Call it $\bullet_a$. Connect $\bullet_a$
with the highest dot in $b_2$ not higher than $\bullet_a$ (assuming
periodic boundary conditions if necessary).
\item Repeat this with all unconnected dots in $b_1$.
\item Slide all $(r_2-r_1)$ unpaired dots from $b_2$ to $b_1$.
The result is $b'_1\otimes b'_2$.
\end{enumerate}

\begin{example}
Let $n=7$, $b_2=7532$ and $b_1=651$. Then
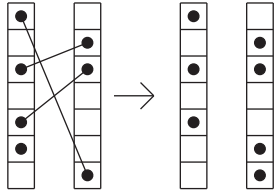
\begin{equation*}
\begin{picture}(100,70)(0,0)
\Line(0,0)(10,0)
\Line(0,10)(10,10)
\Line(0,20)(10,20)
\Line(0,30)(10,30)
\Line(0,40)(10,40)
\Line(0,50)(10,50)
\Line(0,60)(10,60)
\Line(0,70)(10,70)
\Line(0,0)(0,70)
\Line(10,0)(10,70)
\CCirc(5,15){2}{Black}{Black}
\CCirc(5,25){2}{Black}{Black}
\CCirc(5,45){2}{Black}{Black}
\CCirc(5,65){2}{Black}{Black}
\Line(25,0)(35,0)
\Line(25,10)(35,10)
\Line(25,20)(35,20)
\Line(25,30)(35,30)
\Line(25,40)(35,40)
\Line(25,50)(35,50)
\Line(25,60)(35,60)
\Line(25,70)(35,70)
\Line(25,0)(25,70)
\Line(35,0)(35,70)
\CCirc(30,5){2}{Black}{Black}
\CCirc(30,45){2}{Black}{Black}
\CCirc(30,55){2}{Black}{Black}
\Line(30,55)(5,45)
\Line(30,45)(5,25)
\Line(30,5)(5,65)
\Line(40,35)(55,35)
\Line(55,35)(50,40)
\Line(55,35)(50,30)
\Line(65,0)(75,0)
\Line(65,10)(75,10)
\Line(65,20)(75,20)
\Line(65,30)(75,30)
\Line(65,40)(75,40)
\Line(65,50)(75,50)
\Line(65,60)(75,60)
\Line(65,70)(75,70)
\Line(65,0)(65,70)
\Line(75,0)(75,70)
\CCirc(70,25){2}{Black}{Black}
\CCirc(70,45){2}{Black}{Black}
\CCirc(70,65){2}{Black}{Black}
\Line(90,0)(100,0)
\Line(90,10)(100,10)
\Line(90,20)(100,20)
\Line(90,30)(100,30)
\Line(90,40)(100,40)
\Line(90,50)(100,50)
\Line(90,60)(100,60)
\Line(90,70)(100,70)
\Line(90,0)(90,70)
\Line(100,0)(100,70)
\CCirc(95,5){2}{Black}{Black}
\CCirc(95,15){2}{Black}{Black}
\CCirc(95,45){2}{Black}{Black}
\CCirc(95,55){2}{Black}{Black}
\end{picture}
\end{equation*}
so that $b'_1=753$ and $b'_2=6521$.
\end{example}

\subsection{Energy function}
Explicit formulas may be given for the energy functions in type
$A$. Let $B_j=B^{r_j,s_j}$ with $J$-highest weight vector $u_j$
for $j=1,2$. Given the normalization \eqref{eq:extremal H} one has
the explicit formula \cite{S3,SW}
\begin{equation} \label{eq:A local energy}
  H(b_2\otimes b_1) = -\min(r_1,r_2)\min(s_1,s_2) +
  d_{\max\{s_1,s_2\}}(\shape(\PSchen(b_2 b_1)))
\end{equation}
where $d_s(\la)$ is the number of cells of $\la$ that lie in
columns strictly to the right of the $s$-th column. It is not hard
to see that the values of $H$ are nonpositive integers.

\begin{remark} \label{rem:different energy}
In \cite{S3,SW} a different normalization was used, namely, only
the $d$ function was present.
\end{remark}

Let $R=(R_1,\dots,R_L)$ and $b=b_L\otimes \dots \otimes b_1\in
B_R$ where $B_R$ is defined in \eqref{eq:B_R}. Write $D_R$ and
$E_R$ for $D_{B_R}$ and $E_{B_R}$ in \eqref{eq:energy} and
\eqref{eq:NY}. Since $D_{B^{r,s}}=0$,
\begin{equation} \label{eq:A D=E}
D_R = E_R.
\end{equation}

\begin{prop} \label{pp:energy dual} For all $b\in B_R$,
\begin{equation} \label{eq:energy dual}
  E_{R^{\cd\ed}}(b^{\cd\ed})=E_R(b).
\end{equation}
\end{prop}
\begin{proof} Let $B=B_R$, $B'$ be any reordering of the tensor factors of
$B_R$, and $\R:B\rightarrow B'$ any composition of local
isomorphisms. By Propositions \ref{pp:R dual} and \ref{pp:R ev},
\begin{equation} \label{eq:iso dual}
  \R(b^{\cd\ed}) = \R(b)^{\cd\ed}.
\end{equation}
Using this one may reduce to the case $R=(R_1,R_2)$. The crystal
$B_R$ is connected (see section \ref{sec:A simple}). It then
suffices to check \eqref{eq:energy dual} for the single element
$u_2\otimes u_1$ where $u_j$ is the $J$-highest weight vector of
$B_j$, by \eqref{eq:iso dual} and Proposition \ref{pp:dual flip}.
Again by Proposition \ref{pp:dual flip} $(u_2\otimes
u_1)^{\cd\ed}$ is a $J$-highest weight vector of the same weight
(namely $-w_0(\wt(u_2\otimes u_1))$) as $v_2\otimes v_1$ where
$v_j$ is the unique $J$-highest weight vector in $B_j^{\cd}$. But
there is only one such $J$-highest weight vector in
$B_{R^{\cd\ed}}$ of this weight, so the two vectors must agree. It
follows that both sides of \eqref{eq:energy dual} evaluate to
zero.
\end{proof}

\subsection{Kostka polynomials} Let $\la$ be a partition (or an
element of $\Pfin^+$ via $\slwt$ as in section \ref{sec:root A}),
$R=(R_1,\dotsc,R_L)$ a sequence of rectangles and $B_R$ as in
\eqref{eq:B_R}. Define
\begin{align}
  K_{\la,R}(q) &= q^{||R||} \X(B_R,\la;q) \\
  \widetilde{K}_{\la,R}(q) &= \X(B_R,\la;q^{-1}).
\end{align}
with $X(B_R,\la;q)$ as in \eqref{eq:onedimsum} and
\begin{equation*}
||R|| = \sum_{1\le i<j\le L} \min(r_i,r_j) \min(s_i,s_j).
\end{equation*}
It was shown in \cite{NY} (see also \cite{S3,SW}) that
$K_{\la,R}(q)$ is the Kostka polynomial $K_{\la \mu}(q)$ if $R$ is
a sequence of single rowed partitions where $r_i=1$ and
$s_i=\mu_i$. Hence we call $K_{\la,R}(q)$ the generalized
\textit{Kostka polynomial}. Both $K_{\la,R}(q)$ and
$\widetilde{K}_{\la,R}(q)$ have nonnegative integer coefficients.

\section{Littlewood--Richardson tableaux} \label{sec:LR}

\subsection{Definition of Littlewood--Richardson tableaux}
The $J$-components of $B_R$ are parametrized by the combinatorial
objects called Littlewood--Richardson (LR) tableaux. These
objects, being in the multiplicity space, should be considered a
completely different kind of object than the elements of $B_R$. It
is a beautiful coincidence in type $A$ that operators that act on the
representation space such as crystal operators and duality, also
act in the multiplicity space in certain situations.

We recall the notion of an $R$-LR word \cite{S1}, which can be
formulated as follows.
Let $(\eta_1,\eta_2,\dots,\eta_L)\in\N^L$ be a sequence of
positive integers summing to $N$, $A=\{1<2<\dots<N\}$, $A_1$ the
subinterval of $A$ given by the first $\eta_1$ numbers, $A_2$ the
next $\eta_2$ numbers, etc. Say that a word $u$ in the alphabet
$A$ is $\eta$-balanced if $\et_i(u)=\ft_i(u)=0$ for all $i\in A$
that are not maximum in one of the subintervals $A_j$.

Let $R=(R_1,\dots,R_L)$ be a sequence of rectangular
partitions such that $R_j$ has $\eta_j$ rows and
$\mu_j$ columns.  
Let $\gamma(R)=(\mu_1^{\eta_1},\dots,\mu_L^{\eta_L})\in\N^N$ be
obtained by juxtaposing the parts of the $R_j$.  Define the
set $W(R)$ of $R$-LR words by $u\in W(R)$ if and only if
$u$ has content $\gamma(R)$ and $u$ is $\eta$-balanced.
This definition is equivalent to the one in \cite{S1}.

Conversely, for an $\eta$-balanced word $u$,
let $\gamma\in\N^N$ be the content of $u$.
Then $\gamma_i=\gamma_{i'}$ for all $i,i'$ in the same
subalphabet $A_j$.  Let $\mu_j$ be this
common value for the subalphabet $A_j$
and $R_j$ be the rectangular partition having
$\eta_j$ rows and $\mu_j$ columns.  Then
$u$ is an $R$-LR word.

Let $\LRT(R)$ be the subset of $W(R)$ consisting of tableaux of
partition shape (via the identification of a tableau with its
column-reading word) and $\LRT(\la,R)$ the subset of $\LRT(R)$
consisting of tableaux of partition shape $\la$.  The set
$\LRT(\la,R)$ is empty unless $\la$ has at most $N$ parts.

\subsection{Recording tableaux}
\label{sec rt}

The way in which the LR tableaux parametrize the
multiplicity space of $B_R$, is by suitable
recording tableaux for Schensted's insertion.

Consider a given factorization $u=c_N\dotsm c_2 c_1$ of $u$ where
$c_j$ is a column word in the alphabet $[n]$. Let $Q=\QQ(c_N\dotsm
c_2 c_1)$ be the unique filling of the Ferrers diagram of the
shape of $\PSchen(c_N \dotsm c_2c_1)$ such that the shape of
$Q|_{[j]}$ is the shape of $\PSchen(c_j \dotsm c_1)$ for all $j$.
Since each $c_j$ is a column word it follows that the transpose
$Q^t$ of $Q$ is a semistandard tableau. In more traditional
language, $Q$ is the recording tableau for the column insertion of
the word $c_N \dotsm c_2 c_1$ such that the insertion of the
letters in $c_j$ are recorded by the letter $j$.

Fix $R=(R_1,\ldots,R_L)$ and let $R^t=(R_1^t,\ldots,R_L^t)$ where
${}^t$ denotes transpose.
Let $n$ be a positive integer such that $n\ge \mu_j$ for
all $j$. One may regard elements of $B_{R^t}$ as tableaux of a
fixed skew shape in the alphabet $[n]$. The following result, up
to labeling of the recording tableaux, is a special case of the
main theorem in \cite{Wh}.

\begin{prop} \label{path LR}
There is a bijection $B_{R^t}\rightarrow \bigcup_\la B(\la^t)
\times \LRT(\la,R)$ given by $b\mapsto (\PSchen(b),\QQ(b)^t)$
where $\la$ runs over partitions with $\la_1\le n$.
\end{prop}

Recall that the operation ${}^\cd$ was defined for
tableaux $t$ in the alphabet $[n]$ with at most $N$ columns
and for partitions with at most $n$ rows and $N$ columns.
We define an operation ${}^\pd$ which is
a similar kind of operation but with $n$ and $N$ exchanged.

Let $\la^\pd$ be the partition obtained by skewing
the rectangle $(n^N)$ by $\la$ and rotating 180 degrees,
or equivalently $\la^\pd_i = n-\la_{N+1-i}$.
Let $R^\pd=(R^\pd_1,\ldots,R^\pd_L)$ where
$R^\pd_j = ((n-\mu_j)^{\eta_j})$ for each $j$.
This complements the widths of the rectangles in $R$.

For a sequence $u$ of $n$ column words in the alphabet $[N]$, we
write $u^\pd$ for the dual sequence of column words instead of
$u^\cd$.

\begin{prop} \label{pp:record dual}
\begin{equation*}
  \QQ((c_N \dotsm c_2c_1)^{\cd\ed})^t =
  \QQ(c_N\dotsm c_2c_1)^{t \pd}
\end{equation*}
\end{prop}
\begin{proof} Let $Q_1$ and $Q_2$ be the tableaux on the left and
right hand sides.  Let $v$ be obtained from $u$ by replacing
$c_N$ by $\emptyset$ and let $Q_1'$ and $Q_2'$ be the corresponding
tableaux for $v$.  By induction on the number of
nonempty factors $Q_1'=Q_2'$.  By definition the first
$N-1$ column words of $u^{\cd\ed}$ and $v^{\cd\ed}$ agree.
To conclude $Q_1=Q_2$ it suffices to show:
\begin{enumerate}
\item $Q_i'=Q_i|_{[N-1]}$ for $i=1,2$.
\item $Q_1$ and $Q_2$ have the same shape.
\end{enumerate}
Part 1 holds by the definition of $\QQ$. For part 2 one has
\begin{equation*}
\begin{split}
  \shape(\QQ(u^{\cd\ed})^t)&= \shape(\QQ(u^{\cd\ed}))^t=
\shape(\PSchen(u^{\cd\ed}))^t\\
  &= \shape(\PSchen(u^\cd)^\ev)^t=  \shape(\PSchen(u^\cd))^t=  \shape(\PSchen(u)^\cd)^t \\
  &=  \shape(\PSchen(u))^{\cd t}=  \shape(\PSchen(u))^{t \pd}=
\shape(\QQ(u))^{t\pd} \\
  &=  \shape(\QQ(u)^{t\pd})
\end{split}
\end{equation*}
by the fact that the $\PSchen$ and $\QSchen$ tableaux of a given
word have the same shape, Propositions \ref{pp:flip P} and
\ref{pp:dual P}, the fact that ${}^\ev$ is shape-preserving, and
the rule for how ${}^\cd$ and ${}^\pd$ change shapes.
\end{proof}

\begin{prop} There is a bijection
$\LRT(\la,R)\rightarrow \LRT(\la^\pd,R^\pd)$ given by
$t\mapsto t^\pd$.
\end{prop}
\begin{proof} Let $t\in\LRT(\la,R)$. Then $t^\pd$ is a tableau of
shape $\la^\pd$ by Proposition \ref{pp:dual irr}.
It has the correct content to be $R^\pd$-LR by definition since the
corresponding rectangles in $R$ and $R^\pd$ have the same heights.
Since $t$ is $\eta$-balanced, $t^\pd$ is also, by \eqref{eq:dual crystal}.
Thus the desired map is well-defined. Since ${}^\pd$ is an
involution the map is bijective.
\end{proof}

\begin{example} Let $\mu=(3,2,1)$ and $\eta=(2,3,1)$ and $n=5$.  Then
$N=6$ and $R$ and $R^\pd$ (whose rectangles have their rows are filled with the
letters of their corresponding subintervals) are given by
\begin{align*}
  & R &\qquad & R^\pd \\
  &\begin{array}{ccc} \hline
     \mcbn{1} & \mcrn{1} & \mcrn{1} \\
     \cline{1-3} \mcbn{2} & \mcrn{2} & \mcrn{2} \\
     \cline{1-3}
  \end{array}
  & \qquad
  &\begin{array}{cc} \hline
     \mcbn{1} & \mcrn{1}  \\
     \cline{1-2} \mcbn{2} & \mcrn{2} \\ \cline{1-2}
  \end{array} \\
  &\begin{array}{cc}
     \hline \mcbn{3} & \mcrn{3} \\
     \cline{1-2} \mcbn{4} & \mcrn{4} \\
     \cline{1-2} \mcbn{5} & \mcrn{5} \\ \cline{1-2}
  \end{array}
  & \qquad
  &\begin{array}{ccc} \hline
     \mcbn{3} & \mcrn{3} & \mcrn{3} \\
     \cline{1-3} \mcbn{4} & \mcrn{4} & \mcrn{4} \\
     \cline{1-3} \mcbn{5} & \mcrn{5} & \mcrn{5} \\ \cline{1-3}
  \end{array} \\
  &\begin{array}{c}
     \cline{1-1} \mcbn{6} \\ \cline{1-1}
  \end{array}
  & \qquad
  &\begin{array}{cccc} \hline
     \mcbn{6} & \mcrn{6} & \mcrn{6} & \mcrn{6} \\ \cline{1-4}
  \end{array}
\end{align*}
Let $\la=(5,4,3,1,0,0)$.  Then $\la^\pd=(5,5,4,2,1,0)$.
Below is a tableau $t\in \LRT(\la,R)$ and its dual
$t^\pd\in\LRT(\la^\pd,R^\pd)$.
\begin{align*}
& t & \qquad & t^\pd \\
&
\begin{array}{ccccc} \hline
  \mcbn{1} & \mcrn{1} & \mcrn{1} & \mcrn{3} & \mcrn{3} \\ \cline{1-5}
  \mcbn{2} & \mcrn{2} & \mcrn{2} & \mcrn{4} & \mc \\ \cline{1-4}
  \mcbn{4} & \mcrn{5} & \mcrn{6} & \mc & \mc \\ \cline{1-3}
  \mcbn{5} & \mc & \mc & \mc & \mc \\ \cline{1-1}
\\ \\
\end{array}
& \qquad &
\begin{array}{ccccc} \hline
  \mcbn{1} & \mcrn{1} & \mcrn{3} & \mcrn{3} & \mcrn{3} \\ \cline{1-5}
  \mcbn{2} & \mcrn{2} & \mcrn{4} & \mcrn{4} & \mcrn{6} \\ \cline{1-5}
  \mcbn{4} & \mcrn{5} & \mcrn{5} & \mcrn{6} & \mc \\ \cline{1-4}
  \mcbn{5} & \mcrn{6} & \mc & \mc & \mc \\ \cline{1-2}
  \mcbn{6} & \mc & \mc & \mc & \mc \\ \cline{1-1} \\
\end{array}
\end{align*}
\end{example}

\subsection{Generalized automorphisms of conjugation}

Given $1\le p\le L-1$, let $\tau_p R$ be obtained from $R$ by
exchanging the $p$-th and $(p+1)$-st rectangles. There is a
bijection \cite{S1} given by the generalized automorphisms of
conjugation $\tau_p:\LRT(\la,R)\rightarrow \LRT(\la,\tau_p R)$. It
extends uniquely to a bijection $\tau_p:W(R)\rightarrow W(\tau_p
R)$ by $\PSchen(\tau_p u)=\tau_p \PSchen(u)$ and $\QQ(\tau_p
u)=\QQ(u)$.

The bijections $\tau_p$ and the combinatorial $R$-matrices $\R_p$
are related as follows (which may also be used as the definition
of $\tau_p$).

\begin{prop} \label{R auto} For any $b\in B_{R^t}$ and $1\le p\le L-1$,
$\QQ(\R_p(b))^t = \tau_p \QQ(b)^t$.
\end{prop}

\begin{prop} \label{gen auto dual}
The following diagram commutes:
\begin{equation}
\begin{CD}
  \LRT(\la,R) @>{\tau_p}>> \LRT(\la,\tau_p R) \\
  @V{\pd}VV    @VV{\pd}V \\
  \LRT(\la^\pd,R^\pd) @>>{\tau_p}> \LRT(\la^\pd,\tau_p R^\pd)
\end{CD}
\end{equation}
\end{prop}
\begin{proof} Let $Q\in\LRT(R)$.
By Proposition \ref{path LR} there is a $b\in B_{R^t}$
such that $\QQ(b)^t=Q$.  Then
\begin{equation*}
\begin{split}
  \tau_p(Q^\pd)&=  \tau_p(\QQ(b)^{t\pd})=  \tau_p(\QQ(b^{\cd\ed})^t)
   =\QQ(\R_p(b^{\cd\ed}))^t\\
   &=\QQ(\R_p(b)^{\cd\ed}))^t=\QQ(\R_p(b))^{t\pd}
    =\tau_p(\QQ(b)^t)^\pd = \tau_p(Q)^\pd
\end{split}
\end{equation*}
by Propositions \ref{pp:record dual}, \ref{R auto}, \ref{pp:R
dual} and \ref{pp:R ev}.
\end{proof}

\subsection{Embeddings}

There are embeddings $\theta_R:\LRT(\la,R)\rightarrow
\LRT(\la,\rows(R))$ where $\rows(R)$ is the sequence of
single-rowed shapes of sizes given by $\gamma(R)$; they are
defined as compositions of two kinds of elementary embeddings
\cite{S2,SW}. The second kind of elementary embedding is given by
the $\tau_p$. The first kind is given as follows. Suppose that
$R_1=(k^a)$ and $R_2=(k^b)$ where $a-1\ge b+1$. Write
$R'=((k^{a-1}),(k^{b+1}),R_3,\dots,R_L)$.  Then there is an
embedding $\iota_{k,a,b}:\LRT(\la,R)\hookrightarrow \LRT(\la,R')$.
This extends to an embedding $W(R)\rightarrow W(R')$ by $u\mapsto
u'$ where $\PSchen(u')=\iota_{k,a,b}(\PSchen(u))$ and
$\QQ(u')=\QQ(u)$.

Again fix $n$.
We have $R^\pd=(((n-k)^a),((n-k)^b),R_3^\pd,\dots,R_L^\pd)$
and ${R'}^\pd=
(((n-k)^{a-1}),((n-k)^{b+1}),R_3^\pd,\dots,R_L^\pd)$.

\begin{prop} \label{elem embed dual} The following diagram commutes:
\begin{equation*}
\begin{CD}
  \LRT(\la,R) @>{\iota_{k,a,b}}>> \LRT(\la,R') \\
  @V{\pd}VV    @VV{\pd}V \\
  \LRT(\la^\pd,R^\pd) @>>{\iota_{n-k,a,b}}> \LRT(\la^\pd,{R'}^\pd)
\end{CD}
\end{equation*}
\end{prop}
\begin{proof} One reduces to the two rectangle case by
restriction to the subinterval $A_1\cup A_2$.
Recalling that the tensor product of two rectangles is multiplicity-free
and using the duality symmetry of tensor product multiplicities,
one has
\begin{equation*}
|\LRT(\la,R)|=|\LRT(\la^\pd,R^\pd)| \le
|\LRT(\la,R')|=|\LRT(\la^\pd,{R'}^\pd)| \le 1
\end{equation*}
from which the result follows.
\end{proof}

\begin{cor} \label{embedding dual}
The following diagram commutes:
\begin{equation*}
\begin{CD}
  \LRT(\la,R) @>{\theta_R}>> \LRT(\la,\rows(R)) \\
@V{\pd}VV @VV{\pd}V    \\
  \LRT(\la^\pd,R^\pd) @>>{\theta_{R^\pd}}> \LRT(\la^\pd,\rows(R)^\pd)
\end{CD}
\end{equation*}
\end{cor}
\begin{proof} This holds by
Propositions \ref{gen auto dual} and \ref{elem embed dual}.
\end{proof}

\subsection{Generalized charge and duality}

Let $R=(R_1,\dots,R_L)$ such that $R_j$ has $\eta_j$ rows and
$\mu_j$ columns. The generalized charge is a function
$\LRT(\la,R)\rightarrow\N$ defined as follows. For $L=1$
$\charge_R$ is the zero function. For $L=2$ and $u\in W(R_1,R_2)$,
let $d_{R_2,R_1}(u)$ be the number of cells in the shape of
$\PSchen(u)$ to the right of the $\max\{\mu_1,\mu_2\}$-th column
(see \eqref{eq:A local energy}). For general $L$ and
$t\in\LRT(\la,R)$, define \cite{S1,SW}
\begin{equation*}
  \charge_R(t) = \dfrac{1}{L!} \sum_{\tau\in \Symm_L}
  \sum_{j=1}^{L-1} (L-j) d_j(\tau t)
\end{equation*}
where $d_j(\tau t)$ refers to the above function $d$ applied to
the shape of the $\PSchen$ tableau of the restriction of $\tau t$
to the $j$-th and $(j+1)$-st subalphabets for $\tau R$.

\begin{prop} \label{charge energy} \cite{S3}
For all $b\in B_{R^t}$, $-E_{R^t}(b) = \charge_R(\QQ(b)^t)$.
\end{prop}

\begin{cor} \label{charge dual} For all $t\in \LRT(R)$,
$\charge_R(t) = \charge_{R^\pd}(t^\pd)$.
\end{cor}
\begin{proof} This follows from
Propositions \ref{path LR}, \ref{charge energy}, and
\ref{pp:energy dual}.
\end{proof}

\section{Rigged configurations} \label{sec:RC}

Rigged configurations are combinatorial objects which
were first introduced in the context of Bethe Ansatz
studies of spin models in statistical mechanics \cite{KKR,KR}.
Here we recall their definition and a bijection from
Littlewood--Richardson tablaux (or equivalently classically
restricted paths of type $A_n^{(1)}$) to rigged configurations.
The relation between duality and this bijection is given
in Theorem \ref{cdual RC}. We will need these results
later in Section \ref{sec:fermi} in the proof of fermionic
formulas of type $\Dt$, $\At$ and $\Cn$. In Section \ref{sec:fermi}
we will also introduce rigged configurations associated to
$\Dt$, $\At$ and $\Cn$. The rigged configurations of this section
correspond to the algebra $A_n^{(1)}$.

\subsection{Definition of rigged configurations}

Let $\la$ be a partition and let $R=(R_1,\dots,R_L)$ be a sequence of
rectangles with $R_j$ having $r_j$ rows and $s_j$ columns.

The set of admissible configurations $\Conf(\la,R)$ is the set of
all sequences of partitions $\nu=(\nu^{(1)},\nu^{(2)},\dotsc)$ subject to
the constraints
\begin{equation*}
\begin{split}
|\nu^{(k)}| & = -\sum_{j=1}^k \la_j + \sum_{a=1}^L s_a \min\{r_a,k\}\\
P_i^{(k)}(\nu) & \ge 0
\end{split}
\end{equation*}
for $k\ge 1$ and $i\ge 0$. Here $P_i^{(k)}(\nu)$ is the vacancy number
of the parts (strings) of length $i$ in $\nu^{(k)}$ defined as
\begin{equation*}
P_i^{(k)}(\nu) = Q_i(\nu^{(k-1)})-2 Q_i(\nu^{(k)})+Q_i(\nu^{(k+1)})
 +Q_i(\xi^{(k)}(R)),
\end{equation*}
where $\nu^{(0)}$ is the empty partition, $Q_i(\rho)$ is the size of
the first $i$ columns of the partition $\rho$, and
$\xi^{(k)}(R)$ is the partition whose parts are the widths of the
rectangles in $R$ of height $k$.

The set $\RC(\la,R)$ of rigged configurations is defined as follows.
An element of $(\nu,J)\in \RC(\la,R)$ consists of a configuration
$\nu\in\Conf(\la,R)$ together with a rigging $J$. The rigging $J$ is
a double sequence of partitions
\begin{equation*}
J=\{ J^{(k,i)}\}_{i,k\ge 1}
\end{equation*}
such that $J^{(k,i)}$ is a partition in a box of width $P_i^{(k)}(\nu)$
and height $m_i(\nu^{(k)})$ where $m_i(\nu^{(k)})$ is the number
of parts of $\nu^{(k)}$ of size $i$.

The cocharge of a rigged configuration $(\nu,J)\in \RC(\la,R)$ is
defined by
\begin{equation*}
\begin{split}
\cc(\nu,J)&=\cc(\nu)+\sum_{i,k\ge 1}|J^{(k,i)}|\\
\text{with}\quad \cc(\nu)&=
 \sum_{k,i\ge 1}\alpha_i^{(k)}(\alpha_i^{(k)}-\alpha_i^{(k+1)}),
\end{split}
\end{equation*}
where $\alpha_i^{(k)}$ is the size of the $i$-th column of $\nu^{(k)}$.

A rigged configuration $(\nu,J)$ can also be viewed as follows.
Each part of $J^{(k,i)}$ labels a part of length $i$ in $\nu^{(k)}$.
The pair $(i,x)$ of a part of $\nu^{(k)}$ together with its label
is called a string. Then $(\nu,J)$ is the multiset of strings.
A string $(i,x)$ is called singular if $x=P_i^{(k)}(\nu)$.

\subsection{Duality and rigged configurations}

Let $\la$ be a partition and $R$ a sequence of rectangles $R_a$ with
$r_a$ rows and $s_a$ columns. Fix $n$ such that $n\ge \la_1$ and
$n\ge s_a$ for all $1\le a\le L$, and let $N=\sum_{a=1}^L r_a$.

\begin{prop} \label{rc dual}\mbox{}
\begin{enumerate}
\item There is a bijection
${}^\pd:\Conf(\la^t,R^t) \rightarrow \Conf(\la^{\pd t},R^{\pd t})$ defined by
 $\nu^{\pd(k)}=\nu^{(n-k)}$ for all $1\le k\le n-1$.
\item $P^{(k)}_i(\nu^\pd)=P^{(n-k)}_i(\nu)$ for all $i\ge 1$
 and $1\le k\le n-1$.
\item There is an induced bijection
 ${}^\pd:\RC(\la^t,R^t) \rightarrow \RC(\la^{\pd t},R^{\pd t})$ defined by
 $\nu^{\pd(k)}=\nu^{(n-k)}$ and $J^{\pd(k,i)}=J^{(n-k,i)}$
 for all $1\le k\le n-1$ and $i\ge 1$.
\item $\cc((\nu,J)^\pd)=\cc(\nu,J)$.
\end{enumerate}
\end{prop}

\begin{proof}
Note that $|\nu^{\pd(k)}|=|\nu^{(n-k)}|=-\sum_{j=1}^{n-k}\la_j^t
+\sum_{a=1}^L r_a\min\{s_a,n-k\}$.
Since $\sum_j \la_j^t=\sum_a r_a s_a$,
$\sum_{j>n-k}\la_j^t = Nk-\sum_{j>n-k}\la^{\pd t}_{n+1-j}$
and $N=\sum_a r_a$ it follows that
$|\nu^{\pd(k)}|=-\sum_{j=1}^k\la_j^{\pd t}+\sum_{a=1}^L r_a
\min\{n-s_a,k\}$ so that indeed
$\nu^{\pd(k)}\in \Conf(\la^{\pd t},R^{\pd t})$. The other assertions
follow easily from the definitions.
\end{proof}

\subsection{Bijection between LR tableaux and rigged configurations}

In ref.~\cite{KSS} the existence of several bijections between
LR tableaux and rigged configurations was established.  Let
$\phib$ and $\phit$ be the columnwise quantum and coquantum number
bijections, respectively.

Let us describe the algorithm for the bijection
$\phib_R:\LRT(\la,R)\to \RC(\la^t,R^t)$ explicitly.

Let $t\in\LRT(\la,R)$ be an LR tableau for a partition $\la$ and
a sequence of rectangles $R=(R_1,\dotsc,R_L)$ where $R_a$ has $r_a$ rows
and $s_a$ columns. Set
$\tilde{A}_j=[|R_1|+\dotsm+|R_{j-1}|+1,|R_1|+\dotsm+|R_j|]$.
Relabel the letters in $t$ such that the $i$-th occurrence (from the left)
of the $a$-th largest letter in $A_j$ is mapped to
$|R_1|+\dotsm+|R_{j-1}|+(i-1)r_j+a$. Call the relabelled tableau $T$.
Set $|\la|=M$. Let $\ZC_j$ be the tableau of shape $R_j$ such that the
columns of $\ZC_j$ are successively labelled by the letters in $\tilde{A}_j$.

The rigged configuration $(\nu,J)=\phib_R(t)$ is obtained recursively
by constructing a rigged configuration $(\nu,J)_{(x)}$ for each letter
$1\le x\le M$ occurring in $T$.
Set $(\nu,J)_{(0)}=\emptyset$. Suppose that $x\in \tilde{A}_j$, and
denote the column index of $x$ in $T$ by $c$ and the column index
of $x$ in $\ZC_j$ by $c'$. Define the numbers $\ell^{(k)}$ for
$c'\le k<c$ as follows. Let $\ell^{(c-1)}$ be the length of the longest
singular string in $(\nu,J)_{(x-1)}^{(c-1)}$. Now select inductively
a singular string in $(\nu,J)_{(x-1)}^{(k)}$ for $k=c-2,c-3,\dotsc,c'$ whose
length $\ell^{(k)}$ is maximal such that $\ell^{(k)}\le \ell^{(k+1)}$; if no
such string exists set $\ell^{(k)}=0$.
Then $(\nu,J)_{(x)}$ is obtained from $(\nu,J)_{(x-1)}$ by adding one
box to the selected strings with labels such that they remain
singular, leaving all other strings unchanged.
Then the image of $t$ under $\phib_R$ is given by $(\nu,J)=(\nu,J)_{(M)}$.

For the above algorithm it is necessary to be able to compute
the vacancy numbers of an intermediate configuration $\nu_{(x)}$.
Suppose $x$ occurs in $\ZC_j$ in column $c'$. In general
$R_{(x)}=(R_1,\dotsc,R_{j-1},\shape(\ZC_j|_{[1,x]}))$ is not a sequence
of rectangles. If $\shape(\ZC_j|_{[1,x]})$ is not a rectangle
one splits it into two rectangles, one of width $c'$ and
one of width $c'-1$. The vacancy numbers are calculated with respect to
this new sequence of rectangles.

By definition set $\ell^{(k)}=\infty$ for $k\ge c$ and
$\ell^{(k)}=0$ if $k<c'$.

\begin{example}\label{ex:bij}
Let $\la=(4,3,2,2,1,1)$, $R=((1,1,1),(2),(2,2,2,2))$. Consider
$t\in \LRT(\la,R)$ given by
\begin{equation*}
t=
\begin{array}{cccc} \hline
  \mcbn{1} & \mcrn{4} & \mcrn{5} & \mcrn{5} \\ \cline{1-4}
  \mcbn{2} & \mcrn{6} & \mcrn{6} \\ \cline{1-3}
  \mcbn{3} & \mcrn{7} \\ \cline{1-2}
  \mcbn{4} & \mcrn{8} \\ \cline{1-2}
  \mcbn{7} \\ \cline{1-1}
  \mcbn{8} \\ \cline{1-1}
\end{array}.
\qquad \text{Then} \qquad
T=
\begin{array}{cccc} \hline
  \mcbn{1} & \mcrn{5} & \mcrn{6} & \mcrn{10} \\ \cline{1-4}
  \mcbn{2} & \mcrn{7} & \mcrn{11} \\ \cline{1-3}
  \mcbn{3} & \mcrn{12} \\ \cline{1-2}
  \mcbn{4} & \mcrn{13} \\ \cline{1-2}
  \mcbn{8} \\ \cline{1-1}
  \mcbn{9} \\ \cline{1-1}
\end{array}.
\end{equation*}
The non-trivial steps of the above algorithm for $\phib_R$
are given in Table~\ref{table:ex}.
A rigged partition is represented by its Ferrers diagram where
to the right of each part the corresponding rigging is indicated.
The vacancy numbers are given to the left of each part.
For example $R_{(10)}=((1,1,1),(2),(2,1,1,1))$ so that the vacancy numbers
of $(\nu,J)_{(10)}$ are calculated with respect to the sequence of
rectangles $((1,1,1),(2),(2),(1,1,1))$.
\begin{table}
\begin{center}
\begin{equation*}
\begin{array}{|r|lll|} \hline &&&\\[-3mm]
 x & (\nu,J)_{(x)}^{(1)} & (\nu,J)_{(x)}^{(2)} & (\nu,J)_{(x)}^{(3)}
 \\[1mm] \hline &&&\\
 6  & \begin{array}{r|c|l} \cline{2-2} 1& &1 \\ \cline{2-2} \end{array}
    & \begin{array}{r|c|l} \cline{2-2} 0& &0 \\ \cline{2-2} \end{array}
    & \\ &&&\\
 7  & \begin{array}{r|c|c|l} \cline{2-3} 1& & &1 \\ \cline{2-3} \end{array}
    & \begin{array}{r|c|l} \cline{2-2} 0& &0 \\ \cline{2-2} \end{array}
    & \\ &&&\\
 9  & \begin{array}{r|c|c|l} \cline{2-3} 1& & &1 \\ \cline{2-3} \end{array}
    & \begin{array}{r|c|l} \cline{2-2} 0& &0 \\ \cline{2-2} \end{array}
    & \\ &&&\\
 10 & \begin{array}{r|c|c|l} \cline{2-3} 2& & &1 \\ \cline{2-3} \end{array}
    & \begin{array}{r|c|l} \cline{2-2} 0& &0 \\ \cline{2-2}
       0& &0 \\ \cline{2-2} \end{array}
    & \begin{array}{r|c|l} \cline{2-2} 0& &0 \\ \cline{2-2} \end{array} \\
    &&&\\
 11 & \begin{array}{r|c|c|l} \cline{2-3} 3& & &1 \\ \cline{2-3} \end{array}
    & \begin{array}{r|c|c|l} \cline{2-3} 0&&&0\\ \cline{2-3}
      0&&\multicolumn{2}{l}{0}\\ \cline{2-2} \end{array}
    & \begin{array}{r|c|l} \cline{2-2} 0& &0 \\ \cline{2-2} \end{array} \\
    &&&\\
 13 & \begin{array}{r|c|c|l} \cline{2-3} 1& & &1 \\ \cline{2-3} \end{array}
    & \begin{array}{r|c|c|l} \cline{2-3} 0&&&0\\ \cline{2-3}
      0&&\multicolumn{2}{l}{0}\\ \cline{2-2} \end{array}
    & \begin{array}{r|c|l} \cline{2-2} 0& &0 \\ \cline{2-2} \end{array} \\
    &&& \\ \hline
\end{array}
\end{equation*}
\end{center}
\caption{\label{table:ex}
 Example for the bijection algorithm (see Example~\ref{ex:bij})}
\end{table}
\end{example}

Let $\flipq$ be the involution on rigged configurations
that  complements the quantum numbers. More precisely
$\flipq(\nu,J)=(\nu,\tilde{J})$ where $\tilde{J}$ is obtained from
$J$ by complementing each partition $J^{(k,i)}$ within the box of
width $P_i^{(k)}(\nu)$ and height $m_i(\nu^{(k)})$.
By definition, the coquantum number bijection is
$\phit_R=\flipq \circ \phib_R$.

It was shown in \cite{KSS} that the bijection $\phit_R$
preserves the statistics.
\begin{theorem}\cite[Theorem 8.3]{KSS}  \label{thm statistics}
For $t\in \LRT(\la,R)$ we have
\begin{equation*}
\charge_R(t)=\cc(\phit_R(t)).
\end{equation*}
\end{theorem}

\subsection{Bijection from paths to rigged configurations}
\label{sec:paths RC}

Proposition \ref{path LR} yields a bijection $B_R\to \bigcup_\la
B(\la) \times \LRT(\la^t,R^t)$. An element $b\in \Path(B_R,\la)$
will be mapped to $(y,t)$ where $y$ is the unique Yamanouchi
tableau of shape and content $\la$ (see the end of section
\ref{sec:crystal A}). Hence the above bijection restricts to a
bijection between paths and LR tableaux $\Path(B_R,\la)\to
\LRT(\la^t,R^t)$. This in turn induces a bijection between paths
and rigged configurations. It will be useful for later to state
this bijection explicitly. By abuse of notation we denote the
bijection induced by $\phib$ also by $\phib:\Path(B_R,\la)\to
\RC(\la,R)$, and similarly for $\phit$.

Let $b=b_L\otimes\cdots\otimes b_1\in \Path(B_R,\la)$.
Recall that $R_i$ in $R=(R_1,\ldots,R_L)$ has $r_i$ rows and
$s_i$ columns.
Factor each step $b_i$ into columns $b_i=u^{s_i}_i\cdots u^1_i$.
Under the bijection $\Path(B_R,\la)\to \LRT(\la^t,R^t)$, which
maps $b\in\Path(B_R,\la)$ to $t\in \LRT(\la^t,R^t)$,
a letter $c\in u_i^j$ is mapped to the letter
$\sum_{k=1}^{i-1}s_k+j$ in column $c$ in $t$.
As before the rigged configuration $(\nu,J)$ corresponding to
$b$ is constructed recursively, this time on $1\le i\le L$,
$1\le j\le s_i$ and increasing $c\in u_i^j$ in this order.
Suppose the letter $c$ is at height $c'$ in $u_i^j$.
Then select singular strings of length $\ell^{(c-1)},\ell^{(c-2)},
\ldots,\ell^{(c')}$ in the $(c-1)$-th, $(c-2)$-th, ..., $c'$-th
rigged partition such that the $\ell^{(k)}$ are maximal subject to
the condition $\ell^{(c-1)}\ge \ell^{(c-2)}\ge \cdots\ge \ell^{(c')}$.
As before add a box to each selected string, keeping it singular and leaving
all other strings unchanged.
If $c\in u_i^j$ and $j<s_i$ the vacancy number is calculated with respect to
the sequence of rectangles $(R_1,\ldots,R_{i-1},(j^{c'}),((s_i-j)^{c'-1}))$
where $c'$ is the height of $c$ in $b_i$.

\begin{example}
Let $\la=(6,4,2,1)$, $R=((3),(1,1),(4,4))$ and
\begin{equation*}
b=
\begin{array}{cccc} \hline
  \mcbn{1} & \mcrn{1} & \mcrn{2} & \mcrn{3} \\ \cline{1-4}
  \mcbn{2} & \mcrn{2} & \mcrn{3} & \mcrn{4} \\ \cline{1-4}
\end{array}
\otimes
\begin{array}{c} \hline
  \mcbn{1} \\ \cline{1-1}
  \mcbn{2} \\ \cline{1-1}
\end{array}
\otimes
\begin{array}{ccc} \hline
  \mcbn{1} & \mcrn{1} & \mcrn{1} \\ \cline{1-3}
\end{array}
\end{equation*}
in $\Path(B_R,\la)$. Then the corresponding LR tableau in
$\LRT(\la^t,R^t)$ is the $t$ of Example \ref{ex:bij}.
For example the letter $4$ in $u_3^1=\begin{array}{c} \hline \mcbn{3} \\
\cline{1-1} \mcbn{4} \\ \cline{1-1} \end{array}$ corresponds to the letter $5$
in the fourth column of $t$ and will select singular string lengths
$\ell^{(3)}\ge \ell^{(2)}$ in the second and third rigged partition.
Similarly, the letter $3$ in $u_3^2=\begin{array}{c} \hline \mcbn{2} \\
\cline{1-1} \mcbn{3} \\ \cline{1-1} \end{array}$ corresponds to the letter
$6$ in the third column of $t$ and will select a singular string of length
$\ell^{(2)}$ in the second rigged partition.
For example, the rigged configuration in Table \ref{table:ex}
for $x=10$ corresponds to the intermediary paths
\begin{equation*}
\begin{array}{cccc} \hline
  \mcbn{1} & \mcrn{1} & \mcrn{2} & \mcrn{3} \\ \cline{1-4}
  \mcbn{} & \mcrn{} & \mcrn{} & \mcrn{4} \\ \cline{1-4}
\end{array}
\otimes
\begin{array}{c} \hline
  \mcbn{1} \\ \cline{1-1}
  \mcbn{2} \\ \cline{1-1}
\end{array}
\otimes
\begin{array}{ccc} \hline
  \mcbn{1} & \mcrn{1} & \mcrn{1} \\ \cline{1-3}
\end{array}\; .
\end{equation*}
\end{example}

\subsection{Duality under the bijection}

\begin{lemma}\label{lem ineq}
Let $R=(R_1,\dotsc,R_L)$ where $R_L$ is a single row with at least two boxes.
Fix $t\in\LRT(\la,R)$ and let $t'$ (resp. $t''$) be obtained from $t$
by removing the rightmost (resp. two rightmost) letters $L$ from $t$.
Denote by $\ell^{(k)}$ (resp. $\lb^{(k)}$) the lengths of the
selected strings under the bijection algorithm for $\phib$
going from $t'$ to $t$ (resp. $t''$ to $t'$).
Then
\begin{equation*}
\ell^{(k)}\le \lb^{(k-1)}.
\end{equation*}
\end{lemma}

\begin{proof}
Let $(\nu,J)=\phib(t)$, $(\nu',J')=\phib(t')$, and $(\nu'',J'')=
\phib(t'')$.
The vacancy numbers change by
\begin{equation}\label{v change}
P_i^{(k)}(\nu')-P_i^{(k)}(\nu'')=
\chi(\lb^{(k-1)}<i\le \lb^{(k)})-\chi(\lb^{(k)}<i\le \lb^{(k+1)}).
\end{equation}
Let $c$ (resp. $\overline{c}$) be the column index of the box
$t/t'$ (resp. $t'/t''$). Then $c>\overline{c}$.
Since $\lb^{(\overline{c})}=\infty$ it follows from \eqref{v change}
that there are no singular strings of length $i$
with $\lb^{(\overline{c}-1)}<i$ in $(\nu',J')^{(\overline{c})}$.
Hence $\ell^{(\overline{c})}\le \lb^{(\overline{c}-1)}$.
By induction on $k$ it follows that
$\ell^{(k)}\le \ell^{(k+1)} \le \lb^{(k)}$ for $k<\overline{c}$.
Since by \eqref{v change} the strings in the range
$\lb^{(k-1)}<i\le \lb^{(k)}$ are nonsingular we must have
$\ell^{(k)} \le \lb^{(k-1)}$.
\end{proof}

\begin{theorem}\cite[Theorem 8.3]{KSS} \label{RC embed}
The following diagram commutes:
\begin{equation*}
\begin{CD}
  \LRT(\la,R) @>{\phit_R}>> \RC(\la^t,R^t) \\
  @V{\theta_R}VV  @VV{\mathrm{inclusion}}V \\
  \LRT(\la,\rows(R)) @>>{\phit_{\rows(R)}}> \RC(\la^t,\rows(R)^t)
\end{CD}
\end{equation*}
\end{theorem}

\begin{theorem} \label{cdual RC}
The following diagram commutes:
\begin{equation*}
\begin{CD}
  \LRT(\la,R) @>{\phit_R}>> \RC(\la^t,R^t) \\
  @V{\pd}VV  @VV{\pd}V \\
  \LRT(\la^\pd,R^\pd) @>>{\phit_{R^\pd}}>
  \RC(\la^{\pd t},R^{\pd t})
\end{CD}
\end{equation*}
The diagram also commutes with $\phib$ replacing $\phit$.
\end{theorem}

\begin{proof}
By Corollary \ref{embedding dual} and Theorem \ref{RC embed},
one may reduce to the case that
$R=(R_1,\dotsc,R_L)$ is a sequence of single-rowed shapes.
Recall that $\phib_R=\flipq \circ \phit_R$.
By Proposition \ref{rc dual} it is obvious that
$\flipq$ commutes with the duality map on rigged configurations.
Thus it suffices to prove the theorem where $\phib$ replaces
$\phit$.

We proceed by induction on $L$. The theorem is true for $L=0$.
Let $t\in\LRT(\la,R)$ and $c_1<\dotsc<c_{\ell}$ be the column indices
of the letters $L$ in $t$.
Then in $t^\pd$ the column indices of the letters $L$ are
in the alphabet $\{1<\dotsm<n\}$ with letters
$n+1-c_{\ell},\dotsc,n+1-c_1$ omitted.
Call them $c_i^\pd$ for $1\le i\le n-\ell$.
Let $t'$ be the tableau obtained from $t$ by removing all letters $L$
and $R'=(R_1,\dotsc,R_{L-1})$.
By induction $\phib_{{R'}^\pd}({t'}^\pd)=\phib_{R'}(t')^\pd$.
Hence it suffices to show that the addition of the letters $L$ to
$t'$ under $\phib_R$ is (up to reversal of the order of
all partitions) equal to the addition of the letters $L$
to ${t'}^\pd$ under $\phib_{R^\pd}$.

Set $(\nu_0,J_0):=\phib_{R'}(t')$ and let $(\nu_{i-1},J_{i-1})$ be
the rigged configuration corresponding to
$t'$ with letters $L$ added in columns $c_1,\dotsc, c_{i-1}$.
Adding the letter $L$ in column $c_i$ has the effect on the rigged
configurations of selecting singular strings in
$(\nu,J)^{(k)}$ for $k=c_i-1,c_i-2,\dotsc,i$ of length
$\ell_i^{(k)}$ maximal such that
$\ell_i^{(c_i-1)}\ge \ell_i^{(c_i-2)} \ge \dotsm \ge \ell_i^{(i)}
\ge \ell_i^{(i-1)}=\dotsm=\ell_i^{(0)}=0$,
adding a box to the selected strings and making them singular again.
By Lemma \ref{lem ineq}, $\ell_i^{(k)}\le \ell_{i-1}^{(k-1)}$.
It follows that $\ell_i^{(k)}$ with $1\le i\le \ell$ and
$i\le k<c_i$ is uniquely defined to be the length of the
maximal singular string in $(\nu_0,J_0)^{(k)}$
such that
\begin{equation}\label{ineq}
\ell_i^{(k)}\le \min\{\ell_i^{(k+1)},\ell_{i-1}^{(k-1)}\}
\end{equation}
where $\ell_i^{(k+1)}=\infty$ if $k\ge c_i-1$ and $\ell_{i-1}^{(k-1)}=\infty$
if $i=1$ or $k>c_{i-1}$.
Define a matrix $M$ with $\ell$ rows and $c_\ell-\ell$ columns
with entries $M_{ij}=\ell_i^{(j+i-1)}$.
The entries are weakly increasing along a row and weakly
decreasing along a column. Some of the entries may be $\infty$.

Now do the analogous construction for $t^\pd$. Call the selected
strings under $\phib$, $s_i^{(k)}$ for $1\le i\le n-\ell$ and
$i\le k< c_i^\pd$. By the same arguments as before they are
uniquely defined as the lengths of the maximal singular string in
$(\nu_0^\pd,J_0^\pd)^{(k)}=\phib_{{R'}^\pd}({t'}^\pd)$ such that
$s_i^{(k)}\le \min\{s_i^{(k+1)},s_{i-1}^{(k-1)}\}$ where
$s_i^{(k+1)}=\infty$ if $k\ge c_i^\pd-1$ and
$s_{i-1}^{(k-1)}=\infty$ if $i=1$ or $k>c_{i-1}^\pd$. This yields
a matrix $M^\pd$ with entries $M_{ij}^\pd=s_i^{(j+i-1)}$.

It is clear from \eqref{ineq} that the entries in $M$
can either be defined inductively row by row, top to bottom,
right to left, or column by column, right to left, top to bottom.
The same is true for $M^\pd$. Since by induction
$(\nu_0,J_0)^{(k)}=(\nu_0^\pd,J_0^\pd)^{(n-k)}$ for $1\le k<n$
it follows that $M_{i,j}=M^\pd_{c_\ell-\ell+1-j,\ell+1-i}$ or
equivalently $\ell_i^{(k)}=s_{c_\ell-\ell-k+i}^{(c_\ell-k)}$.
This implies $(\nu_\ell,J_\ell)^{(k)}=(\nu^\pd_\ell,J^\pd_\ell)^{(n-k)}$
as desired.
\end{proof}

\begin{example}
Let $R=((1),(3),(2),(4),(3),(2),(4))$ and
\begin{equation*}
t=\begin{array}{|c|c|c|c|c|c|} \hline 1&2&2&3&4&5\\ \hline 2&4&4&5&7&7\\
 \hline 3&5&6&7&\mc&\mc\\ \cline{1-4} 4&7&\mc&\mc&\mc&\mc\\ \cline{1-2}
 6&\mc&\mc&\mc&\mc&\mc\\ \cline{1-1} \end{array}
\quad \text{so that} \quad
t'=\begin{array}{|c|c|c|c|c|c|} \hline 1&2&2&3&4&5\\ \hline 2&4&4&5&\mc&\mc\\
 \cline{1-4} 3&5&6&\mc&\mc&\mc\\ \cline{1-3}
 4&\mc&\mc&\mc&\mc&\mc\\ \cline{1-1} 6&\mc&\mc&\mc&\mc&\mc\\ \cline{1-1}
\end{array}.
\end{equation*}
Then $c_1=2$, $c_2=4$, $c_3=5$, $c_4=6$ and
\begin{equation*}
\begin{split}
(\nu_0,J_0)&=
\begin{array}{r|c|l} \cline{2-2} 1&&1\\ \cline{2-2} \end{array}
\qquad
\begin{array}{r|c|c|l}  \cline{2-3} 0&&&0\\ \cline{2-3} 2&&
 \multicolumn{2}{l}{0}\\ \cline{2-2}
\end{array}
\qquad
\begin{array}{r|c|l} \cline{2-2} 0&&0 \\ \cline{2-2} 0&&0\\
\cline{2-2} 0&&0\\ \cline{2-2} \end{array}
\qquad\quad
\begin{array}{r|c|l} \cline{2-2} 1&&1\\ \cline{2-2}1&&0\\ \cline{2-2}
\end{array}
\qquad
\begin{array}{r|c|l} \cline{2-2} 0&&0\\ \cline{2-2} \end{array}\\
(\nu_1,J_1)&=
\begin{array}{r|c|c|l} \cline{2-3} 1&&&1\\ \cline{2-3} \end{array}
\quad
\begin{array}{r|c|c|l}  \cline{2-3} 1&&&0\\ \cline{2-3} 2&&
 \multicolumn{2}{l}{0}\\ \cline{2-2}
\end{array}
\qquad
\begin{array}{r|c|l} \cline{2-2} 0&&0 \\ \cline{2-2} 0&&0\\
\cline{2-2} 0&&0\\ \cline{2-2} \end{array}
\qquad\quad
\begin{array}{r|c|l} \cline{2-2} 1&&1\\ \cline{2-2}1&&0\\ \cline{2-2}
\end{array}
\qquad
\begin{array}{r|c|l} \cline{2-2} 0&&0\\ \cline{2-2} \end{array}\\
(\nu_2,J_2)&=
\begin{array}{r|c|c|l} \cline{2-3} 1&&&1\\ \cline{2-3} \end{array}
\quad
\begin{array}{r|c|c|l}  \cline{2-3} 1&&&0\\ \cline{2-3} 1&&
 \multicolumn{2}{l}{1}\\ \cline{2-2} 1&& \multicolumn{2}{l}{0}\\
 \cline{2-2}
\end{array}
\qquad
\begin{array}{r|c|c|l}  \cline{2-3} 0&&&0\\ \cline{2-3} 1&&
 \multicolumn{2}{l}{0}\\ \cline{2-2} 1&& \multicolumn{2}{l}{0}\\
 \cline{2-2}
\end{array}
\qquad
\begin{array}{r|c|l} \cline{2-2} 1&&1\\ \cline{2-2}1&&0\\ \cline{2-2}
\end{array}
\qquad
\begin{array}{r|c|l} \cline{2-2} 0&&0\\ \cline{2-2} \end{array}\\
(\nu_3,J_3)&=
\begin{array}{r|c|c|l} \cline{2-3} 1&&&1\\ \cline{2-3} \end{array}
\quad
\begin{array}{r|c|c|l}  \cline{2-3} 1&&&0\\ \cline{2-3} 1&&
 \multicolumn{2}{l}{1}\\ \cline{2-2} 1&& \multicolumn{2}{l}{0}\\
 \cline{2-2}
\end{array}
\qquad
\begin{array}{r|c|c|l}  \cline{2-3} 0&&&0\\ \cline{2-3} 0&&
 \multicolumn{2}{l}{0}\\ \cline{2-2} 0&& \multicolumn{2}{l}{0}\\
 \cline{2-2} 0&& \multicolumn{2}{l}{0}\\ \cline{2-2}
\end{array}
\qquad
\begin{array}{r|c|c|l} \cline{2-3} 1&&&1\\ \cline{2-3} 2&&
 \multicolumn{2}{l}{0} \\ \cline{2-2}
\end{array}
\quad
\begin{array}{r|c|l} \cline{2-2} 0&&0\\ \cline{2-2} \end{array}\\
(\nu_4,J_4)&=
\begin{array}{r|c|c|l} \cline{2-3} 1&&&1\\ \cline{2-3} \end{array}
\quad
\begin{array}{r|c|c|l}  \cline{2-3} 1&&&0\\ \cline{2-3} 1&&
 \multicolumn{2}{l}{1}\\ \cline{2-2} 1&& \multicolumn{2}{l}{0}\\
 \cline{2-2}
\end{array}
\qquad
\begin{array}{r|c|c|l}  \cline{2-3} 0&&&0\\ \cline{2-3} 0&&
 \multicolumn{2}{l}{0}\\ \cline{2-2} 0&& \multicolumn{2}{l}{0}\\
 \cline{2-2} 0&& \multicolumn{2}{l}{0}\\ \cline{2-2}
\end{array}
\qquad
\begin{array}{r|c|c|l} \cline{2-3} 1&&&1\\ \cline{2-3} 1&&
 \multicolumn{2}{l}{1} \\ \cline{2-2} 1&& \multicolumn{2}{l}{0}\\
 \cline{2-2}
\end{array}
\quad
\begin{array}{r|c|c|l} \cline{2-3} 0&&&0\\ \cline{2-3} \end{array}
\end{split}
\end{equation*}
so that
\begin{equation*}
M=\begin{pmatrix}
\ell_1^{(1)} & \ell_1^{(2)}\\
\ell_2^{(2)} & \ell_2^{(3)}\\
\ell_3^{(3)} & \ell_3^{(4)}\\
\ell_4^{(4)} & \ell_4^{(5)}
\end{pmatrix}
=
\begin{pmatrix}
1 & \infty\\ 0&1\\ 0&1\\ 0&1
\end{pmatrix}.
\end{equation*}

On the other hand $c_1^\pd=4$, $c_2^\pd=6$ so that
\begin{equation*}
\begin{split}
(\nu_0^\pd,J_0^\pd)&=
\begin{array}{r|c|l} \cline{2-2} 0&&0\\ \cline{2-2} \end{array}
\qquad
\begin{array}{r|c|l} \cline{2-2} 1&&1\\ \cline{2-2}1&&0\\ \cline{2-2}
\end{array}
\qquad
\begin{array}{r|c|l} \cline{2-2} 0&&0 \\ \cline{2-2} 0&&0\\
\cline{2-2} 0&&0\\ \cline{2-2} \end{array}
\qquad
\begin{array}{r|c|c|l}  \cline{2-3} 0&&&0\\ \cline{2-3} 2&&
 \multicolumn{2}{l}{0}\\ \cline{2-2}
\end{array}
\quad
\begin{array}{r|c|l} \cline{2-2} 1&&1\\ \cline{2-2} \end{array}\\
(\nu_1^\pd,J_1^\pd)&=
\begin{array}{r|c|c|l} \cline{2-3} 0&&&0\\ \cline{2-3} \end{array}
\quad
\begin{array}{r|c|c|l} \cline{2-3} 1&&&1\\ \cline{2-3}
 1&& \multicolumn{2}{l}{0}\\ \cline{2-2}
\end{array}
\quad
\begin{array}{r|c|c|l} \cline{2-3} 0&&&0 \\ \cline{2-3}
 0&& \multicolumn{2}{l}{0}\\ \cline{2-2} 0&& \multicolumn{2}{l}{0}\\
 \cline{2-2} \end{array}
\quad
\begin{array}{r|c|c|l}  \cline{2-3} 1&&&0\\ \cline{2-3} 2&&
 \multicolumn{2}{l}{0}\\ \cline{2-2}
\end{array}
\quad
\begin{array}{r|c|l} \cline{2-2} 1&&1\\ \cline{2-2} \end{array}\\
(\nu_2^\pd,J_2^\pd)&=
\begin{array}{r|c|c|l} \cline{2-3} 0&&&0\\ \cline{2-3} \end{array}
\quad
\begin{array}{r|c|c|l} \cline{2-3} 1&&&1\\ \cline{2-3} 1&&
 \multicolumn{2}{l}{1} \\ \cline{2-2} 1&& \multicolumn{2}{l}{0}\\
 \cline{2-2}
\end{array}
\quad
\begin{array}{r|c|c|l}  \cline{2-3} 0&&&0\\ \cline{2-3} 0&&
 \multicolumn{2}{l}{0}\\ \cline{2-2} 0&& \multicolumn{2}{l}{0}\\
 \cline{2-2} 0&& \multicolumn{2}{l}{0}\\ \cline{2-2}
\end{array}
\quad
\begin{array}{r|c|c|l}  \cline{2-3} 1&&&0\\ \cline{2-3} 1&&
 \multicolumn{2}{l}{1}\\ \cline{2-2} 1&& \multicolumn{2}{l}{0}\\
 \cline{2-2}
\end{array}
\quad
\begin{array}{r|c|c|l} \cline{2-3} 1&&&1\\ \cline{2-3} \end{array}
\end{split}
\end{equation*}
and
\begin{equation*}
M^\pd=\begin{pmatrix}
s_1^{(1)} & s_1^{(2)} & s_1^{(3)} & s_1^{(4)}\\
s_2^{(2)} & s_2^{(3)} & s_2^{(4)} & s_2^{(5)}
\end{pmatrix}
=\begin{pmatrix}
1&1&1&\infty\\ 0&0&0&1
\end{pmatrix}.
\end{equation*}
So indeed $\ell_i^{(k)}=s_{2-k+i}^{(6-k)}$.
\end{example}

\subsection{An embedding} \label{subsec:embedding}
For $r<n$ there is a unique embedding of multiplicity-free
$U_q(A_{2n-1})$-crystals
\begin{equation*}
i_{r,s}: B^{2n-r,s}\otimes B^{r,s} \to B^{2n-r-1,s} \otimes
B^{r+1,s}.
\end{equation*}
Explicitly it is given by $i_{r,s}(u\otimes v)=v'\otimes u'$ where
$\PSchen(uv)=\PSchen(v'u')$. By \eqref{eq:A local energy},
\begin{equation*}
H_{B^{2n-r,s},B^{r,s}}(u\otimes v) =
 H_{B^{2n-r-1,s},B^{r+1,s}}(v'\otimes u')+s.
\end{equation*}
This map also preserves $\et_0$ and $\ft_0$ when such operators act
on the left tensor factor.
An analogous map can be defined on rigged configurations. Let
$R=(R_1,\ldots,R_L)$ be a sequence of rectangles such that
$R_{L-1}=(s^r)$ and $R_L=(s^{2n-r})$ and $R^+$ the same sequence
of rectangles with $R_{L-1}$ replaced by $(s^{r+1})$ and $R_L$
replaced by $(s^{2n-r-1})$. Define $j_{r,s}$ by the following
commutative diagram:
\begin{equation*}
\begin{CD}
  \RC(\la^t,R^t) @>{\mathrm{inclusion}}>> \RC(\la^t,R^{+t}) \\
  @V{\RCtr}VV   @VV{\RCtr}V \\
  \RC(\la,R) @>>{j_{r,s}}> \RC(\la,R^+).
\end{CD}
\end{equation*}
Here $\RCtr:\RC(\la^t,R^t)\to \RC(\la,R)$ is the RC-transpose map
as defined in \cite[Section 7]{KSS}. It can be checked by
direct computation that for $(\nu,J)\in\RC(\la,R)$ the rigged
configuration $j_{r,s}(\nu,J)$ is obtained by adding a singular
string of length $s$ to the $k$-th rigged partition for all
$r+1\le k\le 2n-r-1$.

\begin{theorem} \label{thm:emb}
Let $B$ be any tensor product of type $A_{2n-1}^{(1)}$ crystals of
the form $B^{r',s'}$ and let $R$ (resp. $R^+$) be the sequence of
rectangles associated with $B_R=B^{2n-r,s}\otimes B^{r,s}\otimes
B$ (resp. $B_{R^+}=B^{2n-r-1,s}\otimes B^{r+1,s}\otimes B$). The
following diagram commutes:
\begin{equation*}
\begin{CD}
  \Path(B_R,\la) @>{i_{r,s}\otimes \mathrm{id}_B}>> \Path(B_{R^+},\la)\\
  @V{\phib}VV   @VV{\phib}V \\
  \RC(\la,R) @>>{j_{r,s}}> \RC(\la,R^+).
\end{CD}
\end{equation*}
\end{theorem}
\begin{proof}
The proof follows from the Evacuation Theorem \cite[Theorem 5.6]{KSS},
the Transpose Theorem \cite[Theorem 7.1]{KSS} and
$\PSchen(u^*v^*)=\PSchen({v'}^*{u'}^*)$ for $i_{r,s}(u\otimes v)
=v'\otimes u'$ which is true by Proposition \ref{pp:flip P}.
\end{proof}

\section{Crystals of type $\Dt,\At$ and $\Cn$}
\label{sec:three types}
Let $\gggg$ be an affine Lie algebra of type
$\Dt$, $\At$, $\Atd$, or $\Cn$, defined by the
Dynkin diagrams with distinguished vertex $0$ given in Table
\ref{tab:Dynkin}.
\begin{table}
\caption{}\label{tab:Dynkin}
\begin{tabular}[t]{rl}
\begin{minipage}[b]{4em}
\begin{flushright}
$D_{n+1}^{(2)}$:\\$(n \ge 2)$
\end{flushright}
\end{minipage}&
\begin{picture}(126,20)(-5,-5)
\multiput( 0,0)(20,0){3}{\circle{6}}
\multiput(100,0)(20,0){2}{\circle{6}}
\multiput(23,0)(20,0){2}{\line(1,0){14}}
\put(83,0){\line(1,0){14}}
\multiput( 2.85,-1)(0,2){2}{\line(1,0){14.3}} 
\multiput(102.85,-1)(0,2){2}{\line(1,0){14.3}} 
\multiput(59,0)(4,0){6}{\line(1,0){2}} 
\put(10,0){\makebox(0,0){$<$}}
\put(110,0){\makebox(0,0){$>$}}
\put(0,-5){\makebox(0,0)[t]{$0$}}
\put(20,-5){\makebox(0,0)[t]{$1$}}
\put(40,-5){\makebox(0,0)[t]{$2$}}
\put(100,-5){\makebox(0,0)[t]{$n\!\! -\!\! 1$}}
\put(120,-5){\makebox(0,0)[t]{$n$}}
\end{picture}
\\
&
\\
$A^{(2)}_2$:&
\begin{picture}(26,20)(-5,-5)
\multiput( 0,0)(20,0){2}{\circle{6}}
\multiput(2.958,-0.5)(0,1){2}{\line(1,0){14.084}}
\multiput(2.598,-1.5)(0,3){2}{\line(1,0){14.804}}
\put(0,-5){\makebox(0,0)[t]{$0$}}
\put(20,-5){\makebox(0,0)[t]{$1$}}
\put(10,0){\makebox(0,0){$<$}}
\end{picture}
\\
&
\\
\begin{minipage}[b]{4em}
\begin{flushright}
$A_{2n}^{(2)}$:\\$(n \ge 2)$
\end{flushright}
\end{minipage}&
\begin{picture}(126,20)(-5,-5)
\multiput( 0,0)(20,0){3}{\circle{6}}
\multiput(100,0)(20,0){2}{\circle{6}}
\multiput(23,0)(20,0){2}{\line(1,0){14}}
\put(83,0){\line(1,0){14}}
\multiput( 2.85,-1)(0,2){2}{\line(1,0){14.3}} 
\multiput(102.85,-1)(0,2){2}{\line(1,0){14.3}} 
\multiput(59,0)(4,0){6}{\line(1,0){2}} 
\put(10,0){\makebox(0,0){$<$}}
\put(110,0){\makebox(0,0){$<$}}
\put(0,-5){\makebox(0,0)[t]{$0$}}
\put(20,-5){\makebox(0,0)[t]{$1$}}
\put(40,-5){\makebox(0,0)[t]{$2$}}
\put(100,-5){\makebox(0,0)[t]{$n\!\! -\!\! 1$}}
\put(120,-5){\makebox(0,0)[t]{$n$}}
\end{picture}
\\
&
\\
$A^{(2)\dagger}_2$:&
\begin{picture}(26,20)(-5,-5)
\multiput( 0,0)(20,0){2}{\circle{6}}
\multiput(2.958,-0.5)(0,1){2}{\line(1,0){14.084}}
\multiput(2.598,-1.5)(0,3){2}{\line(1,0){14.804}}
\put(0,-5){\makebox(0,0)[t]{$0$}}
\put(20,-5){\makebox(0,0)[t]{$1$}}
\put(10,0){\makebox(0,0){$>$}}
\end{picture}
\\
&
\\
\begin{minipage}[b]{4em}
\begin{flushright}
$A_{2n}^{(2)\dagger}$:\\$(n \ge 2)$
\end{flushright}
\end{minipage}&
\begin{picture}(126,20)(-5,-5)
\multiput( 0,0)(20,0){3}{\circle{6}}
\multiput(100,0)(20,0){2}{\circle{6}}
\multiput(23,0)(20,0){2}{\line(1,0){14}}
\put(83,0){\line(1,0){14}}
\multiput( 2.85,-1)(0,2){2}{\line(1,0){14.3}} 
\multiput(102.85,-1)(0,2){2}{\line(1,0){14.3}} 
\multiput(59,0)(4,0){6}{\line(1,0){2}} 
\put(10,0){\makebox(0,0){$>$}}
\put(110,0){\makebox(0,0){$>$}}
\put(0,-5){\makebox(0,0)[t]{$0$}}
\put(20,-5){\makebox(0,0)[t]{$1$}}
\put(40,-5){\makebox(0,0)[t]{$2$}}
\put(100,-5){\makebox(0,0)[t]{$n\!\! -\!\! 1$}}
\put(120,-5){\makebox(0,0)[t]{$n$}}
\end{picture}
\\
&
\\
\begin{minipage}[b]{4em}
\begin{flushright}
$C_n^{(1)}$:\\$(n \ge 2)$
\end{flushright}
\end{minipage}&
\begin{picture}(126,20)(-5,-5)
\multiput( 0,0)(20,0){3}{\circle{6}}
\multiput(100,0)(20,0){2}{\circle{6}}
\multiput(23,0)(20,0){2}{\line(1,0){14}}
\put(83,0){\line(1,0){14}}
\multiput( 2.85,-1)(0,2){2}{\line(1,0){14.3}} 
\multiput(102.85,-1)(0,2){2}{\line(1,0){14.3}} 
\multiput(59,0)(4,0){6}{\line(1,0){2}} 
\put(10,0){\makebox(0,0){$>$}}
\put(110,0){\makebox(0,0){$<$}}
\put(0,-5){\makebox(0,0)[t]{$0$}}
\put(20,-5){\makebox(0,0)[t]{$1$}}
\put(40,-5){\makebox(0,0)[t]{$2$}}
\put(100,-5){\makebox(0,0)[t]{$n\!\! -\!\! 1$}}
\put(120,-5){\makebox(0,0)[t]{$n$}}

\end{picture}
\end{tabular}
\end{table}

There is an embedding $P\rightarrow P^A$ of the weight lattice of
$\gggg$ into that of $\A$ which preserves distance up to a fixed
scalar factor. This suggests that $U'_q(\gggg)$-crystals can be
embedded into $U'_q(\A)$-crystals. We introduce the notion of a
virtual $U'_q(\gggg)$-crystal, which is by definition a suitable
subset of a $U'_q(\A)$-crystal. We define the virtual
$U'_q(\gggg)$-crystal $\V^{r,s}$ and give some evidence that it
agrees with the crystal $B^{r,s}$ of Conjecture \ref{conj:Brs}
(whose existence was conjectured in \cite[Conj. 2.1]{HKOTT}).
These conjectures are proved for crystals $B^{r,1}$ of types
$\gggg=\Cn,\At,\Atd$.

\subsection{Affine Lie algebras $\Dt,\At$ and $\Cn$}
We use the notation of section \ref{subsec:affine case}. The null
root $\delta$ is given explicitly by
\begin{equation*}
\delta =
\begin{cases}
\alpha_0+\alpha_1+\cdots+\alpha_{n-1}+\alpha_n   &
  \text{for $\gggg=\Dt$} \\
2\alpha_0+2\alpha_1+\cdots+2\alpha_{n-1}+\alpha_n &
  \text{for $\gggg=\At$} \\
\alpha_0+2\alpha_1+\cdots+2\alpha_{n-1}+2\alpha_n &
  \text{for $\gggg=\Atd$} \\
\alpha_0+2\alpha_1+\cdots+2\alpha_{n-1}+\alpha_n &
  \text{for $\gggg=\Cn$.}
\end{cases}
\end{equation*}
To each $\gggg$ we associate two simple Lie algebras $\gfin$ and $\gt$
\begin{center}
\begin{tabular}{c|cccc}
$\gggg$ & $\Dt$ & $\At$ & $\Atd$ & $\Cn$ \\
\hline
$\gfin$ & $B_n$ & $C_n$ & $B_n$ & $C_n$ \\
$\gt$ &   $B_n$ & $B_n$ & $B_n$ & $C_n$.
\end{tabular}
\end{center}
Note that $\gfin=\gt$ except for $\At$. The algebra $\gt$ will be
used in section \ref{sec:fermi}.

Let $(\cdot|\cdot)$ be the standard bilinear form on $P$
normalized by
\begin{equation} \label{eq:bilinear form}
(\delta|\la)=\inner{c}{\la}\qquad\text{for any $\la\in P$.}
\end{equation}

\subsection{Embedding of weight lattices} Let
\begin{equation} \label{eq:gamma def}
\begin{split}
  \gamma &= \begin{cases}
    1 & \text{for $\gggg=\Dt,\Atd$} \\
    2 & \text{for $\gggg=\At,\Cn$}
  \end{cases} \\
  \gamma' &= \begin{cases}
  1 & \text{for $\gggg=\Dt,\At$} \\
  2 & \text{for $\gggg=\Atd,\Cn$.}
  \end{cases}
\end{split}
\end{equation}
We use the superscript and subscript $A$ to denote the root system
$\A$ (or $A_{2n-1}$). There is an embedding of weight lattices
$\wtembed: P\rightarrow P^A$ defined by
\begin{equation}\label{eq:weight embed}
\begin{split}
\wtembed(\La_0)&=\gamma'\La_0^A\\
\wtembed(\La_i)&=\La_i^A+\La_{2n-i}^A \qquad\text{for $0<i<n$} \\
\wtembed(\La_n)&=\gamma\La_n^A\\
\wtembed(\delta)&=\gamma\delta^A.
\end{split}
\end{equation}
This induces the embedding $\wtembed:\Pfin \rightarrow \Pfin^A$
given by
\begin{equation} \label{eq:finite weight embed}
\begin{split}
\wtembed(\Lab_i)&=\Lab_i^A+\Lab_{2n-i}^A\qquad\text{for $1\le i<n$}\\
\wtembed(\Lab_n)&=\gamma\Lab_n^A.
\end{split}
\end{equation}
One may verify that
\begin{equation} \label{eq:embed roots}
\begin{split}
\wtembed(\alpha_0)&=\gamma'\alpha_0^A\\
\wtembed(\alpha_i)&=\alpha_i^A+\alpha_{2n-i}^A\qquad\text{for $0<i<n$} \\
\wtembed(\alpha_n)&=\gamma\alpha_n^A.
\end{split}
\end{equation}
Observe that
\begin{equation}
(\wtembed(u)|\wtembed(v))_A=\gamma\gamma'(u|v)\qquad \text{for
$u,v\in P$}
\end{equation}
under the normalization \eqref{eq:bilinear form}.
Note that if $\la\in \overline{P}^A$, then $\la\in\Image(\wtembed)$
if and only if the following conditions hold:
\begin{equation}\label{eq:wt image}
\begin{split}
  &\inner{h_i}{\la} = \inner{h_{2n-i}}{\la} \qquad
\text{for $1\le i\le n-1$} \\
  &\text{$\inner{h_n}{\la}$ is divisible by $\gamma$} \\
  &\text{$\inner{h_0}{\la}$ is divisible by $\gamma'$.}
\end{split}
\end{equation}
The third condition is a consequence of the other two. Let $w_0$
be the longest element in the Weyl group of type $A_{2n-1}$. Then
for $\la\in \Image(\wtembed)\cap \Pfin^A$,
\begin{equation}\label{eq:self dual wt}
  -w_0(\la) = \la
\end{equation}
since this holds for $\Psi(\Lab_i)$ for $1\le i\le n$.

\subsection{Virtual $U'_q(\gggg)$-crystals} We wish to define a category of
$U'_q(\gggg)$-crystals which are realized as subsets of
$U'_q(\A)$-crystals. Let $\Vh$ be a finite $U'_q(\A)$-crystal. In
view of \eqref{eq:embed roots} we define the virtual $U'_q(\gggg)$
raising and lowering operators on $\Vh$ by
\begin{equation} \label{eq:virtual e}
\begin{aligned}
\et_0&=(\et_0^A)^{\gamma'}&\qquad\ft_0&=(\ft_0^A)^{\gamma'} \\
\et_i&=\et_{2n-i}^A\et_i^A&\qquad\ft_i&=
\ft_{2n-i}^A\ft_i^A\qquad\text{for $0<i<n$}\\
\et_n&=(\et_n^A)^{\gamma}&\qquad\ft_n&=(\ft_n^A)^{\gamma}.
\end{aligned}
\end{equation}
Let $\V\subset\Vh$ be a nonempty subset such that
\begin{equation}\label{eq:wt in image}
  \Image(\wt_A|_\V) \subset \Image(\wtembed).
\end{equation}
Then $\V$ has a virtual weight function
$\wt:\V\rightarrow P$ defined by
\begin{equation} \label{eq:virtual wt}
  \wtembed(\wt(b)) = \wt_A(b).
\end{equation}
This is well-defined since $\wtembed$ is injective.

A virtual $U'_q(\gggg)$-crystal is a pair $(\V,\Vh)$ such that $\V$
satisfies \eqref{eq:wt in image} and is closed under the virtual
raising and lowering operators $\et_i$ and $\ft_i$ for $0\le i\le
n$. Sometimes the larger crystal $\Vh$ (called the ambient
crystal) is omitted in the notation.

\begin{prop} \label{pp:virtual crystal} A virtual $U'_q(\gggg)$-crystal
$(\V,\Vh)$ is a $U'_q(\gggg)$-crystal in the sense of section
\ref{subsec:crystals} (see also \ref{subsec:affine case}).
\end{prop}
\begin{proof} All the properties are immediate except for
\eqref{eq:string length}. Since $\et_i$ and $\et_{2n-i}$ (resp.
$\ft_i$ and $\ft_{2n-i}$) commute for $1\le i\le n-1$, by the
definitions we have for $b\in\V$
\begin{align} \label{eq:virtual eps}
\epsilon_0(b)&=\floor{\epsilon_0^A(b)/\gamma'} & \qquad
&\varphi_0(b)=\floor{\varphi_0^A(b)/\gamma'}  \\
\epsilon_i(b)&=\min\{\epsilon_i^A(b),\epsilon_{2n-i}^A(b)\}&\qquad&
\text{for $1\le i\le n-1$} \\
\varphi_i(b)&=\min\{\varphi_i^A(b),\varphi_{2n-i}^A(b)\}&\qquad&
\text{for $1\le i\le n-1$} \\
\epsilon_n(b)&=\floor{\epsilon_n^A(b)/\gamma} & \qquad
&\varphi_0(b)=\floor{\varphi_n^A(b)/\gamma}.
\end{align}
Here $\floor{x}$ denotes the largest integer smaller or equal to $x$.
First let us verify \eqref{eq:string length} for $i=0$.
We have that
\begin{equation*}
  \varphi_0^A(b)-\epsilon_0^A(b)=\inner{h_0^A}{\wt_A(b)}
\end{equation*}
is a multiple of $\gamma'$ by \eqref{eq:wt image} applied to
$\wt_A(b)$ and \eqref{eq:string length}
for $b$ viewed as an element of $\Vh$. For $b\in \V$ we have
\begin{equation*}
 \varphi_0(b)-\epsilon_0(b) = \dfrac{\varphi_0^A(b)-\epsilon_0^A(b)}{\gamma'}
 = \dfrac{1}{\gamma'}\inner{h_0^A}{\wt_A(b)} =\inner{h_0}{\wt(b)}
\end{equation*}
by \eqref{eq:virtual eps}, \eqref{eq:virtual wt}, and
\eqref{eq:weight embed}.
A similar calculation establishes \eqref{eq:string length} for $i=n$.
Let $1\le i\le n-1$. Observe that
\begin{equation} \label{eq:strlen image}
  \varphi_i^A(b)-\epsilon_i^A(b)=\inner{h_i^A}{\wt_A(b)}=
\inner{h_{2n-i}^A}{\wt_A(b)} =\varphi_{2n-i}^A(b)-\epsilon_{2n-i}^A(b)
\end{equation}
using \eqref{eq:string length} for $b\in \Vh$ and \eqref{eq:wt image}.
Suppose first that $\varphi_i^A(b)\le \varphi_{2n-i}^A(b)$. Then
$\epsilon_i^A(b)\le\epsilon_{2n-i}^A(b)$ and one has
\begin{equation*}
  \varphi_i(b)-\epsilon_i(b)=\varphi_i^A(b)-\epsilon_i^A(b).
\end{equation*}
If $\varphi_i^A(b)>\varphi_{2n-i}^A(b)$ then
$\epsilon_i^A(b)>\epsilon_{2n-i}^A(b)$ and
\begin{equation*}
  \varphi_i(b)-\epsilon_i(b)=\varphi_{2n-i}^A(b)-\epsilon_{2n-i}^A(b).
\end{equation*}
In either case \eqref{eq:string length} follows by \eqref{eq:strlen image}
and \eqref{eq:weight embed}.
\end{proof}

A morphism of virtual $U'_q(\gggg)$-crystals
$(\V_1,\Vh_1)\rightarrow (\V_2,\Vh_2)$ is a morphism
$\Vh_1\rightarrow \Vh_2$ of the ambient $U'_q(\A)$-crystals which
restricts to a morphism $\V_1\rightarrow \V_2$ of
$U'_q(\gggg)$-crystals.

\subsection{Tensor products of virtual crystals}
Let $(\V_1,\Vh_1)$ and $(\V_2,\Vh_2)$ be virtual
$U'_q(\gggg)$-crystals. We wish to define their tensor product.
Define the virtual crystal operators $\et_i$ and $\ft_i$ on
$\Vh_2\otimes \Vh_1$ by \eqref{eq:virtual e}. Define
$\wt_{\V_2\otimes \V_1}:\V_2\otimes \V_1\rightarrow P$ by
\eqref{eq:tensor wt}. Unfortunately the subset $\V_2\otimes
\V_1\subset \Vh_2\otimes \Vh_1$ need not be closed under $\et_i$
and $\ft_i$.

\begin{example} Let $n=2$, $\gggg=C^{(1)}_2$, and let $\Vh_2=\Vh_1=B^{2,2}_A$
be $U'_q(\A)$-crystals. Let $(\V_2,\Vh_2)$ and $(\V_1,\Vh_1)$ be the
virtual crystals generated by the elements $b_2$ and $b_1$ where
\begin{equation*}
  b_2=\begin{array}{|c|c|} \hline 1&1 \\ \hline 2&3 \\ \hline
\end{array} \qquad
  b_1=\begin{array}{|c|c|} \hline 1&1 \\ \hline 2&2 \\ \hline
\end{array}.
\end{equation*}
Then $\ft_2(b_2\otimes b_1)=\ft_2^A(b_2)\otimes
\ft_2^A(b_1)\not\in \V_2\otimes \V_1$. This disagrees with the usual
tensor product structure on $U'_q(\gggg)$-crystals defined by
\eqref{eq:f on two factors}.
\end{example}

Let $b\in \Vh$. Consider the following conditions.
\begin{equation}\label{eq:aligned}
\begin{split}
&\text{$\varphi_i^A(b)=\varphi_{2n-i}^A(b)$ and
$\epsilon_i^A(b)=\epsilon_{2n-i}^A(b)$ for $1\le i\le
n-1$} \\
&\text{$\varphi_0^A(b)$ and $\epsilon_0^A(b)$ are divisible by
$\gamma'$} \\
&\text{$\varphi_n^A(b)$ and $\epsilon_n^A(b)$ are divisible by
$\gamma$.}
\end{split}
\end{equation}
Say that $b$ is $i$-aligned for $1\le i\le n-1$ if the first
condition holds, $0$-aligned if the second holds, and $n$-aligned
if the third holds. Say that $(\V,\Vh)$ is an aligned virtual
$U'_q(\gggg)$-crystal if, for every $b\in \V$, $b$ is $i$-aligned
for all $0\le i\le n$.

Let $1\le i\le n-1$. If $b$ is $i$-aligned then
\begin{equation} \label{eq:i align}
  \epsilon_i(b)=\epsilon_i^A(b)=\epsilon_{2n-i}^A(b) \qquad\text{
  and }\qquad
  \varphi_i(b)=\varphi_i^A(b)=\varphi_{2n-i}^A(b).
\end{equation}
If $b$ is $n$-aligned then
\begin{equation} \label{eq:n align}
  \epsilon_n(b)=\frac{1}{\gamma}\epsilon_n^A(b) \qquad \text{ and
  } \qquad \varphi_n(b) =\frac{1}{\gamma}\varphi_n^A(b).
\end{equation}
If $b$ is $0$-aligned then
\begin{equation} \label{eq:0 align}
  \epsilon_0(b)=\frac{1}{\gamma'}\epsilon_0^A(b) \qquad\text{ and
  } \qquad
  \varphi_0(b) =\frac{1}{\gamma'}\varphi_0^A(b).
\end{equation}

\begin{remark} \label{rem:aligned} Let $(\V,\Vh)$ be a virtual
$U_q(\gggg)$-crystal and $b\in \Vh$. If $\gamma=1$ then $b$ is
$n$-aligned and if $\gamma'=1$ then $b$ is $0$-aligned.
\end{remark}

\begin{prop} \label{pp:virtual tensor} Aligned virtual
$U'_q(\gggg)$-crystals form a tensor category.
\end{prop}
\begin{proof} Let $(\V_j,\Vh_j)$ be aligned virtual $U'_q(\gggg)$-crystals
for $j=1,2$. Since they are aligned, it follows that for all
$b_2\otimes b_1\in \V_2\otimes \V_1$ and all $i\in I$, $\et_i$ and
$\ft_i$ either act entirely on the left factor or on the right
factor of $b_2\otimes b_1$. Therefore the operators $\et_i$ and
$\ft_i$ coincide with those given by the $U'_q(\gggg)$-crystal
structure on the tensor product $\V_2\otimes \V_1$ defined by
\eqref{eq:e on two factors} and \eqref{eq:f on two factors}, which define
a tensor category. Finally, it is easy to verify that
$(\V_2\otimes \V_1,\Vh_2\otimes \Vh_1)$ is aligned.
\end{proof}

If $(\V_2,\Vh_2)$ and $(\V_1,\Vh_1)$ are aligned virtual
$U'_q(\gggg)$-crystals then define their tensor product to be the
aligned virtual $U'_q(\gggg)$-crystal $(\V_2\otimes \V_1,
\Vh_2\otimes \Vh_1)$.

\subsection{Virtual combinatorial $R$-matrix and energy
function}

Say that the virtual $U'_q(\gggg)$-crystal $(\V,\Vh)$ is simple if
\begin{enumerate}
\item $\Vh$ is a simple $U'_q(\A)$-crystal.
\item $\V$, with its virtual $U'_q(\gggg)$-crystal structure, is
isomorphic to a simple $U'_q(\gggg)$-crystal.
\item $u(\V)=u(\Vh)$.
\end{enumerate}

The following result shows that local isomorphisms and energy
functions exist for aligned simple virtual crystals.

\begin{prop} \label{pp:virtual R}
Let $(\V_1,\Vh_1)$ and $(\V_2,\Vh_2)$ be aligned simple virtual
$U'_q(\gggg)$-crystals. Write
\begin{equation*}
\begin{split}
\Rh &=\R_{\Vh_2,\Vh_1}:\Vh_2\otimes \Vh_1\rightarrow
  \Vh_1\otimes \Vh_2 \\
\Hh&=H_{\Vh_2,\Vh_1}:\Vh_2\otimes\Vh_1\rightarrow\Z \\
\R&=\R_{\V_2,\V_1}:\V_2\otimes\V_1\rightarrow\V_1\otimes\V_2 \\
H&=H_{\V_2,\V_1}:\V_2\otimes\V_1\rightarrow\Z
\end{split}
\end{equation*}
for the local isomorphisms and energy functions for the pair of
simple $U'_q(\A)$-crystals $\Vh_1$ and $\Vh_2$, and the pair of
simple $U'_q(\gggg)$-crystals $\V_1$ and $\V_2$. Then
\begin{align}
\label{eq:virtual R}
  \R &= \Rh|_{\V_2\otimes \V_1} \\
\label{eq:virtual H}
  H &= \dfrac{1}{\gamma'} \Hh|_{\V_2\otimes\V_1}.
\end{align}
\end{prop}
\begin{proof}
$(\V_2\otimes \V_1,\Vh_2\otimes \Vh_1)$ is an aligned simple
virtual $U'_q(\gggg)$-crystal, by Proposition \ref{pp:virtual tensor}
and Theorem \ref{thm:simple}, generated from the extremal vector
\begin{equation*}
u:=u(\V_2\otimes\V_1) = u(\V_2)\otimes u(\V_1) = u(\Vh_2)\otimes
u(\Vh_1) = u(\Vh_2\otimes \Vh_1),
\end{equation*}
which holds by \eqref{eq:extremal tensor} and property 3 of the
definition of simple virtual crystal $(\V_j,\Vh_j)$. Similarly
$(\V_1\otimes\V_2,\Vh_1\otimes\Vh_2)$ is an aligned simple virtual
$U'_q(\gggg)$-crystal with generator
$u':=u(\V_1\otimes\V_2)=u(\Vh_1\otimes\Vh_2)$. Now $\Rh$ is a
$U'_q(\A)$-crystal isomorphism such that $\Rh(u)=u'$. As such it
intertwines with the virtual operators $\et_i$ and $\ft_i$ for
$0\le i\le n$. Therefore $\Rh|_{\V_2\otimes\V_1}$ is a
$U'_q(\gggg)$-crystal isomorphism from the $U'_q(\gggg)$-component
of $u$ to that of $u'$, that is, from $\V_2\otimes \V_1$ to
$\V_1\otimes \V_2$. But there is a unique such map, namely, $\R$.
This proves \eqref{eq:virtual R}.

By abuse of notation, let $H:\V_2\otimes
\V_1\rightarrow\frac{1}{\gamma'}\Z$ be defined by
\eqref{eq:virtual H}. It suffices to show that $H$ satisfies the
defining properties of the energy function given in Theorem
\ref{thm:two tensor}. This is entirely straightforward except for
\eqref{eq:local energy}. Let $b_j\in \V_j$ for $j=1,2$ and write
$\R(b_2\otimes b_1)=b_1'\otimes b_2'$.

Suppose first that $\epsilon_0(b_2)>\varphi_0(b_1)$ and
$\epsilon_0(b_1')>\varphi_0(b_2')$. Since $(\V_j,\Vh_j)$ is
aligned for $j=1,2$, this means
\begin{equation*}
  \epsilon_0^A(b_2)=\gamma'\epsilon_0(b_2)>\gamma'\varphi_0(b_1)=
  \varphi_0^A(b_1)
\end{equation*}
and
\begin{equation*}
  \epsilon_0^A(b_1')=\gamma'\epsilon_0(b_1')>\gamma'\varphi_0(b_2')=
  \varphi_0^A(b_2')
\end{equation*}
by \eqref{eq:0 align}. Applying $\et_0^A$ to $b_2\otimes b_1$
$\gamma'$ times and using \eqref{eq:local energy}, one has
\begin{equation*}
  \Hh(\et_0(b_2\otimes b_1))=
  \Hh((\et_0^A)^{\gamma'}(b_2 \otimes b_1))=
  \Hh(b_2\otimes b_1) - \gamma'.
\end{equation*}
Dividing by $\gamma'$ one obtains \eqref{eq:local energy} for $H$
in the first case. The other cases are similar.
\end{proof}

\subsection{Virtual intrinsic energy}
Let $(\V,\Vh)$ be an aligned simple virtual $U'_q(\gggg)$-crystal.
Suppose the ambient simple crystal $\Vh$ is graded by the
intrinsic energy function $D_{\Vh}$. Then $\V$ has an inherited
intrinsic energy function $D_\V:\V\rightarrow\frac{1}{\gamma'}\Z$,
defined by
\begin{equation} \label{eq:virtual energy}
   D_\V(b) = \dfrac{1}{\gamma'} D_{\Vh}(b).
\end{equation}
Call $(\V,\Vh;D_{\V},D_{\Vh})$ an aligned graded simple virtual
crystal.

Let $(\V_j,\Vh_j;D_{\V_j},D_{\Vh_j})$ be an aligned graded simple
virtual $U'_q(\gggg)$-crystal for $j=1,2,\dotsc,L$. Consider the
simple $U'_q(\gggg)$-crystal $\V=\V_L\otimes \dotsm\otimes \V_1$
and the simple $U'_q(\A)$-crystal $\Vh=\Vh_L\otimes\dotsm\otimes
\Vh_1$. There are two ways to form an intrinsic energy function
for $\V$. One way is to use the fact that it is a tensor product of
graded simple $U'_q(\gggg)$-crystals $(\V_j,D_{\V_j})$ using
\eqref{eq:energy}. The other way is to use the inherited intrinsic
energy function \eqref{eq:virtual energy} coming from the fact
that $(\V,\Vh)$ is an aligned virtual crystal and $(\Vh,D_{\Vh})$
is a graded simple crystal (with $D_{\Vh}$ defined by
\eqref{eq:energy} in terms of $D_{\Vh_j}$). It follows from
\eqref{eq:virtual H} that the two definitions agree.

\subsection{Virtual crystal $\V^{r,s}$} \label{subsec:vc}
For $\gggg$ of type $\Dt$, $\At$, $\Atd$, or $\Cn$, we wish to define a
virtual $U'_q(\gggg)$-crystal $\V^{r,s}$ ($1\le r\le n,s\ge1$).
Let $\Vh^{r,s}$ be the $U'_q(\A)$-crystal
\begin{equation}\label{eq:ac}
\Vh^{r,s} = \begin{cases}
  B^{2n-r,s}_A \otimes B^{r,s}_A & \text{if $r<n$} \\
  B^{n,s}_A & \text{if $r=n$ and $\gggg=\Dt$} \\
  (B^{n,s}_A)^{\otimes 2} & \text{if $r=n$ and $\gggg=\At,\Atd$} \\
  B^{n,2s}_A & \text{if $r=n$ and $\gggg=\Cn$.}
\end{cases}
\end{equation}

The virtual crystal $(\V^{r,s},\Vh^{r,s})$ is defined as the
subgraph of $\Vh^{r,s}$ generated from the extremal vector
$u(\Vh^{r,s})$ by applying the virtual operators $\et_i$ and
$\ft_i$ ($0\le i\le n$).

\begin{conjecture} \label{conj:aligned} $\V^{r,s}$ is aligned.
\end{conjecture}

For $1\le i\le n$, define
\begin{equation} \label{eq:mu def}
  \mu_i = \begin{cases}
  2 & \text{if $\gggg=\Atd$ and $i=n$} \\
  1 & \text{otherwise.}
  \end{cases}
\end{equation}
The crystals $B^{r,s}$ of Conjecture \ref{conj:Brs} come equipped
with the following prescribed decomposition as
$U_q(\gfin)$-crystals:
\begin{equation}\label{eq:decomp}
  B^{r,s} \cong
\begin{cases}
  B(s\Lab_n) & \text{if $r=n$ and $\gggg=\Dt,\Cn$} \\
  \bigoplus_\la B(\la) & \text{otherwise}
\end{cases}
\end{equation}
where $\la\in\Pfin^+$ has the form $\la=\sum_{i=1}^r m_i \mu_i
\Lab_i$ such that $\sum_{i=1}^r m_i \le s$, $s-m_r\in \gamma'\Z$,
and $m_i\in\gamma'\Z$ for $1\le i\le r-1$. Moreover, $B^{r,s}$
comes equipped with the following intrinsic energy function:
\begin{equation} \label{eq:Brs D}
  D_{B^{r,s}}(b) =
  \begin{cases}
  0 & \text{if $r=n$ and $\gggg=\Dt,\Cn$} \\
  (\Lab_n|\la-s\mu_r\Lab_r) & \text{if $b\in B(\la)\subseteq B^{r,s}$.}
  \end{cases}
\end{equation}

\begin{conjecture} \label{conj:virtual} The virtual crystal
$\V^{r,s}$ is the simple crystal $B^{r,s}$ of the
$U'_q(\gggg)$-module $W^{(r)}_s$, with extremal vector
$u(\V^{r,s})=u(\Vh^{r,s})$ of weight $s\mu_r \Lab_r$. In
particular, $\V^{r,s}$ has the $U_q(\gfin)$-crystal decomposition
given by \eqref{eq:decomp} and $D_{\V^{r,s}}=D_{B^{r,s}}$.
\end{conjecture}

It will be shown in Proposition \ref{pp:big enough} that
$\V^{r,s}$ has at least the $U_q(\gfin)$-components specified by
\eqref{eq:decomp} and that the virtual intrinsic energy on these
components agrees with \eqref{eq:Brs D}. If either $\gggg=\Dt$ or
$s=1$, it will be shown in Theorems \ref{thm:D decomp}, \ref{thm:C
decomp}, \ref{thm:A2 decomp}, and \ref{thm:A2d decomp} that
$\V^{r,s}$ is aligned and has exactly the decomposition
\eqref{eq:decomp}.

\subsection{Self-duality} We now consider the problem of
giving explicit conditions for membership in the virtual crystal
$\V^{r,s}$ as a subset of $\Vh^{r,s}$. We shall see that a
necessary condition for membership in $\V^{r,s}$ is self-duality up
to local isomorphism.

Let $\Vho^{r,s}$ be the $U'_q(\A)$-crystal obtained by reversing
the order of the tensor factors in $\Vh^{r,s}$. Let $\Vo^{r,s}$ be
defined as $\V^{r,s}$ is, except with $\Vh^{r,s}$ replaced by
$\Vho^{r,s}$. There is a $U'_q(\A)$-crystal isomorphism
$\R:\Vh^{r,s} \cong \Vho^{r,s}$ where $\R$ is the identity or the
local isomorphism according as $\Vh^{r,s}$ has one or two tensor
factors. The map $\R$ induces a $U'_q(\gggg)$-crystal
isomorphism $\V^{r,s}\rightarrow \Vo^{r,s}$.

\begin{prop} \label{pp:Brs self dual} For every $b\in \V^{r,s}$,
\begin{equation} \label{eq:self dual R}
  b^{\cd\ed}=\R(b).
\end{equation}
\end{prop}
\begin{proof}
First let $b=u(\Vh^{r,s})$. Since $\Vho^{r,s}$ is
multiplicity-free as a $U_q(A_{2n-1})$-crystal, it suffices to
show that both $b^{\cd\ed}$ and $\R(b)$ are $A_{2n-1}$-highest
weigh vectors of the same weight. Since $\R$ is a
$U'_q(\A)$-crystal isomorphism, $\R(b)$ is a
$U_q(A_{2n-1})$-highest weight vector of weight $\wt_A(b)$.
$b^{\cd\ed}\in\Vho^{r,s}$ is an $A_{2n-1}$-highest weight vector
of weight $\wt_A(b^{\cd\ed})=-w_0(\wt_A(b))$ by \eqref{eq:crystal
dual flip}. Since $\wt_A(b) = \wtembed(s\mu_r\Lab_r)$ it follows
that $-w_0(\wt_A(b))=\wt_A(b)$ by \eqref{eq:self dual wt}. So
\eqref{eq:self dual R} holds.

Now suppose $b\in \Vh^{r,s}$ satisfies \eqref{eq:self dual R} and
$\ft_i(b)\not=\emptyset$ for some $0\le i\le n$. It will be shown that
$\ft_i(b)$ also satisfies \eqref{eq:self dual R}. This, together
with a similar proof with $\et_i$ replacing $\ft_i$, suffices. Let
$0\le j\le 2n$. We have
\begin{equation} \label{eq:dual raise}
  \ft_j^A(b)^{\cd\ed}=\et_j^A(b^\cd)^\ed=\ft_{2n-j}^A(b^{\cd\ed})=
  \ft_{2n-j}^A(\R(b))=\R(\ft_{2n-j}^A(b))
\end{equation}
by \eqref{eq:dual crystal}, Theorem \ref{thm:flip}, \eqref{eq:self
dual R}, and the fact that $\R$ is a $U'_q(\A)$-crystal
isomorphism. \eqref{eq:self dual R} holds for $\ft_0(b)$ and
$\ft_n(b)$ by applying \eqref{eq:dual raise} with $j=0$ and $j=n$
respectively. \eqref{eq:self dual R} holds for $\ft_i(b)$ for
$1\le i\le n-1$ by observing that $\ft_j^A$ and $\ft_{2n-j}^A$
commute, and applying \eqref{eq:dual raise} for $\ft_j^A(b)$ and
then for $\ft^A_{2n-j}\ft^A_j(b)$.
\end{proof}

\begin{prop} \label{pp:virtual eps} For all $b\in \Vh^{r,s}$
such that \eqref{eq:self dual R} holds and for all $1\le i\le
n-1$, $b$ is $i$-aligned.
\end{prop}
\begin{proof} Let $1\le i\le n-1$. Since $\R$ is an isomorphism
of $U'_q(\A)$-crystals we have
\begin{equation*}
  \epsilon_i^A(b)=\epsilon_i^A(\R(b))=\epsilon_i^A(b^{\cd\ed})=
  \varphi_{2n-i}^A(b^\cd)=
  \epsilon_{2n-i}^A(b)
\end{equation*}
by Theorem \ref{thm:flip} and \eqref{eq:dual crystal}. The proof
for $\varphi$ is similar.
\end{proof}

\subsection{Virtual $U_q(\gfin)$-crystals} Recall that for the
affine algebras $\gggg$ under consideration $\gfin$ is
either of type $B_n$ or $C_n$. A pair $(\V,\Vh)$ is a virtual
$U_q(\gfin)$-crystal provided that $\Vh$ is a finite
$U_q(A_{2n-1})$-crystal, $\V\subset \Vh$ is a subset that
satisfies \eqref{eq:wt in image} and is closed under the virtual
operators $\et_i$ and $\ft_i$ for $1\le i\le n$ defined by
\eqref{eq:virtual e}. An aligned virtual $U_q(\gfin)$-crystal is
one that consists of elements that are $i$-aligned for $1\le i\le
n$.

Let $\la\in\Pfin^+$ and $\la^A=\wtembed(\la)\in (\Pfin^A)^+$.
Since the virtual operators $\et_i$ and $\ft_i$ (for $1\le i\le
n$) are comprised of type $U_q(A_{2n-1})$ operators and the latter
preserve the property of being a tableau of a given shape, it
follows that $B_A(\la^A)$ is closed under these virtual operators.
Therefore we may define $\V(\la)$ to be the virtual
$U_q(\gfin)$-crystal given by the subset of $B_A(\la^A)$ generated
from $u_{\la^A}$ using the virtual operators $\et_i$ and $\ft_i$
for $1\le i\le n$.

The following crucial theorem is due to T. Baker \cite{B}. Our
approach, which emphasizes the self-duality property, could be
used to give an elegant proof of this result.

\begin{theorem} \label{thm:virtual fin}
$\V(\la)$ is aligned and is isomorphic to $B(\la)$ as a
$U_q(\gfin)$-crystal.
\end{theorem}

\begin{prop} \label{pp:align hwv} Let $(\V,\Vh)$ be an aligned
$U_q(\gfin)$-crystal. Then every $b\in \V$ is in the connected
component of a virtual $U_q(\gfin)$-highest weight vector $u$, and
$u$ is a $U_q(A_{2n-1})$-highest weight vector when viewed as an
element of $\Vh$.
\end{prop}
\begin{proof} Let $b\in \V$. Suppose $\epsilon_j^A(b)>0$ for some
$1\le j\le 2n-1$. Since $(\V,\Vh)$ is an aligned virtual
$U_q(\gfin)$-crystal, it follows that $\et_i(b)\in \V$ where $i=j$
for $1\le j\le n$ and $i=2n-j$ for $n<j\le 2n-1$. Replacing $b$
with $\et_i(b)$ and continuing in this manner, eventually one has
$\epsilon_j^A(b)=0$ for $1\le j\le 2n-1$.
\end{proof}

\begin{prop} \label{pp:virtual B} Let $\la\in \Pfin^+$
and $b\in B_A(\la^A)$. If $b\in \V(\la)$ then
\begin{equation} \label{eq:PSD}
  \PSchen(b^{\cd\ed})=b.
\end{equation}
Moreover if $\gfin=B_n$ and $b$ satisfies \eqref{eq:PSD} then
$b\in \V(\la)$.
\end{prop}
\begin{proof} The proof that elements of $\V(\la)$ satisfy
\eqref{eq:PSD} is similar to the proof of Proposition \ref{pp:Brs
self dual} and is omitted. Suppose $\gfin=B_n$. Let $b\in
B_A(\la^A)$ be such that \eqref{eq:PSD} holds. Let $\V'\subset
B_A(\la^A)$ be the virtual $U_q(B_n)$-crystal generated by $b$. It
suffices to show that $u_{\la^A}\in\V'$. By the proof of
Proposition \ref{pp:Brs self dual}, \eqref{eq:PSD} holds for any
element in $\V'$. But \eqref{eq:PSD} implies that $\V'$ is
$i$-aligned for $1\le i\le n-1$, in a manner similar to the proof
of Proposition \ref{pp:virtual eps}. By Remark \ref{rem:aligned}
$\V'$ is $n$-aligned. Therefore $(\V',B_A(\la^A))$ is an aligned
$U_q(B_n)$-crystal. By Proposition \ref{pp:align hwv}
$u_{\la^A}\in\V'$.
\end{proof}

\subsection{$U_q(C_n)$-crystals} \label{sec:C crystals} For this
section let $\gfin=C_n$. The crystal graphs of the irreducible
integrable $U_q(C_n)$-modules $B(\la)$ for $\la\in \Pfin^+$, were
constructed explicitly by Kashiwara and Nakashima \cite{KN} in
terms of certain tableaux.

The crystal graph $B(\Lab_1)$ is given by
\newcommand{\bx}[1]{\begin{array}{|c|}\hline #1 \\ \hline \end{array}}
\begin{equation} \label{eq:vector rep}
\bx{1} \overset{1}\longrightarrow \bx{2} \overset{2}\longrightarrow \dotsm
\overset{n-1}\longrightarrow \bx{n} \overset{n}\longrightarrow
\bx{\overline{n}} \overset{n-1}\longrightarrow \dotsm
\overset{2}\longrightarrow \bx{\overline{2}} \overset{1}\longrightarrow
\bx{\overline{1}}
\end{equation}
with $\wt(i)=\Lab_i-\Lab_{i-1}=-\wt(\overline{i})$ for $1\le i\le n$.

Let $\AZ$ be the set of vertices of $B(\Lab_1)$, totally ordered
so that the elements appear from smallest to largest going from
left to right in \eqref{eq:vector rep}. Let $\len{u}$ denote the
length of the word $u$. A subset $A \subset [n]$ is often
identified with the column word having precisely the letters in
$A$. Write $A^c= [n]-A$, $\ba{A}=\{\ba{i}\mid i\in A\}$, and
$\ba{A}^c=\ba{[n]}-\ba{A}=\ba{A^c}$. In particular if $u$ and $v$
are column words in the alphabet $[n]$ then $\ba{v} u$ is a column
word in the alphabet $\AZ$ and $(\ba{v}u)^{\cd\ed}=\ba{u}^c v^c$.
If $w$ is a column word in the alphabet $\AZ$ then it can be
written uniquely in the form $w=\ba{v}u$ where $u$ and $v$ are column
words in the alphabet $[n]$. In this case write $u=w_+$ and $v=w_-$.

The crystal $B(\Lab_r)$ may be defined as the connected component in the
tensor product $B(\Lab_1)^{\otimes r}$ of the column word
$r(r-1)\dotsm 2 1$ where the tensor symbols are omitted.
Explicitly, $B(\Lab_r)$ is the set of column words $P$ of length
$r$ in the alphabet $[n]\cup\overline{[n]}$ that satisfy the
one-column condition \cite{KN}
\begin{equation} \label{1CC}
\text{If $i$ and $\overline{i}$ are both in $P$ then
$\len{{P_+}|_{[i]}} + \len{{P_-}|_{[i]}} \le i$.}
\end{equation}

Let $\la\in\Pfin^+$. Write $\la=\sum_{j=1}^p \Lab_{m_j}$ where
$m_1\ge m_2\ge\dotsm\ge m_p$. Then $B(\la)$ may be defined as the
connected component of $u_{\Lab_{m_1}}\otimes\dotsm \otimes
u_{\Lab_{m_p}}$ in $B(\Lab_{m_1})\otimes\dotsm\otimes
B(\Lab_{m_p})$.

\begin{lemma} \label{lem:1CC tableau}
$P$ satisfies \eqref{1CC} if and only if $P^c_- P_+$ is a tableau.
\end{lemma}
\begin{proof} The following are equivalent:
\begin{enumerate}
\item $P$ does not satisfy \eqref{1CC}.
\item There is an index $1\le i\le n$ such that
$i\in P_+$, $i\not\in P^c_-$, and $\len{{P_+}|_{[i]}} >
\len{{P^c_-}|_{[i]}}$.
\item There is an index $1\le i\le n$ such that
$\len{{P_+}|_{[i]}} > \len{{P^c_-}|_{[i]}}$.
\item $P^c_- P_+$ is not a tableau.
\end{enumerate}
Each of the above conditions is obviously equivalent to the next
except for 3 implies 2. For that case just take $i$ to be minimal.
If $i \not\in P_+$ or $i\in P^c_-$ then the index $i-1$ satisfies
3, contradicting minimality.
\end{proof}

We now recall the explicit isomorphism $B(\la)\rightarrow \V(\la)$
given in \cite{B}. The first case is $\la=\Lab_r$. For $P\in
B(\Lab_r)$, let
\begin{equation} \label{eq:pre split column}
\begin{split}
  K&=P_+\cap P_- \\
  J&=\max\{A\subset (P_+\cup P_-)^c\mid |A|=|K| \text{ and }A<K \}
\end{split}
\end{equation}
where the maximum is computed with respect to the partial order on
column words that says $u \le v$ if and only if $uv$ is a tableau.
This given, let
\begin{equation} \label{eq:split column}
  Q_\pm = (P_\pm - K) \cup J.
\end{equation}
Define the map $\embedfin:B(\Lab_r)\rightarrow
B_A(\Lab_{2n-r}+\Lab_r)$ by
\begin{equation}\label{eq:single column embed}
  P\mapsto \ba{Q}^c_+ P^c_- \otimes \ba{Q}_- P_+
\end{equation}
where $B_A(\Lab_{2n-r}+\Lab_r)$ is regarded as the
$U_q(A_{2n-1})$-subcrystal of $B_A(\Lab_{2n-r})\otimes
B_A(\Lab_r)$ by taking the column-reading word. Lemma \ref{lem:1CC
tableau} ensures that the column-reading word indeed corresponds
to a tableau.

Let $\la=\Lab_{m_1}+\Lab_{m_2}+\dotsm+\Lab_{m_p}$ for $1\le m_j\le
n$, $\la^A=\wtembed(\la)$, and $b=P^{(1)} \otimes\dotsm\otimes
P^{(p)}\in B(\la)$ where $P^{(j)}\in B(\Lab_{m_j})$. Define
$\embedfin:B(\la)\rightarrow B_A(\la^A)$ by
\begin{equation} \label{eq:embed}
  \embedfin(b) = \PSchen(\embedfin(P^{(1)}) \embedfin(P^{(2)}) \dotsm
  \embedfin(P^{(p)})).
\end{equation}

We observe that $Q_+$ and $Q_-$ may be computed from $P$ by local
isomorphisms of type $A$ for tensor products of columns.

\begin{lemma} \label{lem:one col iso}
\begin{equation} \label{eq:P to Q}
\R(P^c_- \otimes P_+)= Q_+ \otimes Q_-^c.
\end{equation}
\end{lemma}
\begin{proof} The construction in the definition of $Q_+$ and $Q_-$
is equivalent to the algorithm of section \ref{sec:R sc} used to
compute $\R(P^c_- \otimes P_+)$.
\end{proof}

\begin{example} Let $n=9$ and $P=\ba{3}\ba{8}8763$.
Then $P_-=83$, $P_+=8763$, $K=83$, $J=52$, $Q_-=52$, and
$Q_+=7652$. On the other hand $\R(P^c_- \otimes P_+)$ may be
computed using the algorithm described in section \ref{sec:R sc}.
Note that $P^c_-=9765421$ so that
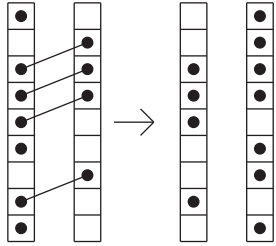
\begin{equation*}
\begin{picture}(100,90)(0,0)
\Line(0,0)(10,0) \Line(0,10)(10,10) \Line(0,20)(10,20)
\Line(0,30)(10,30) \Line(0,40)(10,40) \Line(0,50)(10,50)
\Line(0,60)(10,60) \Line(0,70)(10,70) \Line(0,80)(10,80)
\Line(0,90)(10,90) \Line(0,0)(0,90) \Line(10,0)(10,90)
\CCirc(5,5){2}{Black}{Black} \CCirc(5,15){2}{Black}{Black}
\CCirc(5,35){2}{Black}{Black} \CCirc(5,45){2}{Black}{Black}
\CCirc(5,55){2}{Black}{Black} \CCirc(5,65){2}{Black}{Black}
\CCirc(5,85){2}{Black}{Black}
\Line(25,0)(35,0) \Line(25,10)(35,10) \Line(25,20)(35,20)
\Line(25,30)(35,30) \Line(25,40)(35,40) \Line(25,50)(35,50)
\Line(25,60)(35,60) \Line(25,70)(35,70) \Line(25,80)(35,80)
\Line(25,90)(35,90) \Line(25,0)(25,90) \Line(35,0)(35,90)
\CCirc(30,25){2}{Black}{Black} \CCirc(30,55){2}{Black}{Black}
\CCirc(30,65){2}{Black}{Black} \CCirc(30,75){2}{Black}{Black}
\Line(30,75)(5,65) \Line(30,65)(5,55) \Line(30,55)(5,45)
\Line(30,25)(5,15)
\Line(40,45)(55,45) \Line(55,45)(50,50) \Line(55,45)(50,40)
\Line(65,0)(75,0) \Line(65,10)(75,10) \Line(65,20)(75,20)
\Line(65,30)(75,30) \Line(65,40)(75,40) \Line(65,50)(75,50)
\Line(65,60)(75,60) \Line(65,70)(75,70) \Line(65,80)(75,80)
\Line(65,90)(75,90) \Line(65,0)(65,90) \Line(75,0)(75,90)
\CCirc(70,15){2}{Black}{Black} \CCirc(70,45){2}{Black}{Black}
\CCirc(70,55){2}{Black}{Black} \CCirc(70,65){2}{Black}{Black}
\Line(90,0)(100,0) \Line(90,10)(100,10) \Line(90,20)(100,20)
\Line(90,30)(100,30) \Line(90,40)(100,40) \Line(90,50)(100,50)
\Line(90,60)(100,60) \Line(90,70)(100,70) \Line(90,80)(100,80)
\Line(90,90)(100,90) \Line(90,0)(90,90) \Line(100,0)(100,90)
\CCirc(95,5){2}{Black}{Black} \CCirc(95,25){2}{Black}{Black}
\CCirc(95,35){2}{Black}{Black} \CCirc(95,55){2}{Black}{Black}
\CCirc(95,65){2}{Black}{Black} \CCirc(95,75){2}{Black}{Black}
\CCirc(95,85){2}{Black}{Black}
\end{picture}
\end{equation*}
$K$ is the set of all heights which has a dot on the right, but no
dot on the left. $(P_+\cup P_-)^c$ is the set of heights with a
dot on the left, but no dot on the right. Hence $J$ is obtained by
selecting from top to bottom for each element $i\in K$ the maximal
$h<i$ with a dot on the left and no dot on the right at height
$h$. The elements in $J$ are all selected by the algorithm for the
computation of $\R$ as described in section \ref{sec:R sc}, and
hence both computations of $J$ are equivalent.
\end{example}

For the following Lemma, we identify $\AZ$ with $[2n]$ via the
bijection $i\mapsto i$ and $\ba{i}\mapsto 2n+1-i$ for $i\in[n]$.
\begin{lemma} \label{lem:self dual}
Let $u$ and $v$ be column words in the alphabet $[n]$.  Then there
is a unique pair of column words $u'$ and $v'$ in the alphabet
$[n]$ such that $\len{u}=\len{u'}$, $\len{v}=\len{v'}$, and
(defining $s=\ba{u'}^c v$ and $t=\ba{v'}^c u$)
\begin{equation} \label{self dual}
\R(s\otimes t)=s^{\cd\ed} \otimes t^{\cd\ed}
\end{equation}
where $\R:B_A^{2n-k,1}\otimes B_A^{k,1}\rightarrow
B_A^{k,1}\otimes B_A^{2n-k,1}$ is the local isomorphism of type
$\A$ with $k=n+\len{u}-\len{v}$.
\end{lemma}
\begin{proof} For existence, define $u'$ and $v'$ by
\begin{equation} \label{eq:new pair}
  \R(v\otimes u)=u'\otimes v'.
\end{equation}
By definition $u'$ and $v'$ are column words of the correct
length. Applying the map ${}^{\cd\ed}$ to \eqref{eq:new pair}, by
Propositions \ref{pp:R dual} and \ref{pp:R ev} and the fact that
$\R$ is an involution, one has
\begin{equation} \label{eq:dual new pair}
  \R(\ba{u'}^c \otimes \ba{v'}^c) =
  \ba{v}^c \otimes \ba{u}^c.
\end{equation}
By definition $s=\ba{u'}^c v$ and $t=\ba{v'}^c u$, so that $s$ and
$t^{\cd\ed}$ are column words of the same length and $t$ and
$s^{\cd\ed}$ are column words of the same length. Thus to prove
\eqref{self dual} it is enough to show that
$\PSchen(st)=\PSchen((st)^{\cd\ed})$.  Since the shapes of both
$\PSchen(st)$ and $\PSchen((st)^{\cd\ed})$ have two columns, it
suffices to show that
\begin{enumerate}
\item $\PSchen(st|_\BZ)=\PSchen((st)^{\cd\ed}|_\BZ)$ for $\BZ=[n]$
and for $\BZ=\ba{[n]}$.
\item $\PSchen(st)$ and $\PSchen((st)^{\cd\ed})$ have the same shape.
\end{enumerate}
Now $s^{\cd\ed}=\ba{v}^c u'$ and $t^{\cd\ed}=\ba{u}^c v'$. We have
\begin{equation*}
  \PSchen((st)^{\cd\ed}|_{[n]})=\PSchen(u'v')=\PSchen(vu)=\PSchen(st|_{[n]})
\end{equation*}
by \eqref{eq:new pair} and
\begin{equation*}
  \PSchen((st)^{\cd\ed}|_{\ba{[n]}})
=\PSchen(\ba{v}^c \ba{u}^c) =\PSchen(\ba{u'}^c \ba{v'}^c) =
\PSchen(st|_{\ba{[n]}})
\end{equation*}
by \eqref{eq:dual new pair}. This proves 1. Condition 2 is
equivalent to: $Q_1 = \QQ(st)^t$ and $Q_2 = \QQ((st)^{\cd\ed})^t$
have the same shape. But $Q_2=\QQ(st)^{t\pd}=Q_1^\pd$ by
Proposition \ref{pp:record dual} with $N=2$. So it is enough to
show that $\shape(Q_1)^\pd = \shape(Q_1)$. But this holds since
$\shape(Q_1)$ has at most two rows and exactly $2n$ cells (since
${}^\pd$ is taken within the $2\times 2n$ rectangle).

For uniqueness, observe that $\QQ(v u)^t$ is a tableau with $|u|$
ones and $|v|$ twos, and $\QQ(u' v')^t$ is a tableau of the same
shape with $|v|$ ones and $|u|$ twos. The second tableau is
uniquely specified by this property.
\end{proof}

\begin{prop} \label{pp:virtual C} Let $b=u\otimes v\in
B_A(\Lab_{2n-r}+\Lab_r)$. Write $u_1=u|_{[n]}$ and
$u_2=u|_{[n+1,2n]}$ and similarly for $v$. The following are
equivalent:
\begin{enumerate}
\item $b\in \V(\Lab_r)$.
\item \eqref{eq:PSD} holds, $u_1v_1$ is a tableau and
$\len{u_1} - \len{v_1} = n - r$.
\item \eqref{eq:PSD} holds, $u_2v_2$ is a tableau and
$\len{u_2} - \len{v_2} = n - r$.
\end{enumerate}
\end{prop}
\begin{proof} 1 implies 2: Let $P\mapsto b$ under the $U_q(C_n)$-crystal
isomorphism $B(\Lab_r)\rightarrow \V(\Lab_r)$ given by
\eqref{eq:single column embed}. Then $u_2=\Qb^c_+$, $u_1=P^c_-$,
$v_2=\Qb_-$, and $v_1=P_+$ with $Q_\pm$ defined by Lemma
\ref{lem:one col iso}. By Lemma \ref{lem:1CC tableau} $u_1 v_1$ is
a tableau, with $\len{u_1}-\len{v_1}=n-\len{P_-}-\len{P_+}=
n-\len{P}=n-r$. 2 implies 1: Let $P$ be such that $P^c_-=u|_{[n]}$
and $P_+=v|_{[n]}$. By Lemma \ref{lem:1CC tableau} $P\in
B(\Lab_r)$. By Lemma \ref{lem:self dual} it follows that $P\mapsto
b$ under the isomorphism $B(\Lab_r)\rightarrow \V(\Lab_r)$.

For the equivalence of 2 and 3, suppose that $P$ and $Q$ are such
that \eqref{eq:P to Q} holds. Then the following are equivalent:
\begin{enumerate}
\item $P^c_-P_+$ is a tableau.
\item $Q_+Q^c_-$ is the column-reading word of a
skew two-column tableau with unique southeast corner.
\item $\Qb^c_- \Qb_+$ is a tableau.
\item $\Qb^c_+ \Qb_-$ is a tableau.
\end{enumerate}
Moreover it is easily seen that $\len{P^c_-}-\len{P_+}=
\len{\Qb^c_+}-\len{\Qb_-}$. Now suppose 2 holds. Write $b$ in the
form $\Qb^c_+ P^c_-\otimes \Qb_-P_+$. By the proof of Lemma
\ref{lem:self dual} it follows that \eqref{eq:P to Q} holds. Thus
3 follows. The proof that 3 implies 2 is entirely similar.
\end{proof}

\subsection{$\V^{r,s}$ contains the prescribed
$U_q(\gfin)$-components} \label{subsec:big enough}

\begin{prop} \label{pp:big enough} Let $\la\in\Pfin^+$ be such
that $B(\la)$ is a summand specified by \eqref{eq:decomp}. Then
there is a unique $U_q(\gfin)$-crystal embedding
$\iota_\la:\V(\la)\rightarrow \V^{r,s}$ defined by $b'\mapsto b$
where $b'=\PSchen(b)$, and $D_{\V^{r,s}}=D_{B^{r,s}}$ on
$\Image(\iota_\la)$.
\end{prop}
\begin{proof} Let $\la=\sum_{j=1}^r m_j\Lab_j\in \Pfin^+$ be
any weight such that $\wtembed(\la)$ occurs in the
$U_q(A_{2n-1})$-crystal decomposition of $\Vh^{r,s}$.
$B_A(\wtembed(\la))$ occurs in $\Vh^{r,s}$ with multiplicity one,
and on this component $D_{\Vh^{r,s}}$ has value $0$ if $\Vh^{r,s}$
has one tensor factor and value $-rs+\sum_{j=1}^r j m_j$ otherwise
\cite{S1,SW}. Let $v_\la\in \Vh^{r,s}$ be the unique
$U_q(A_{2n-1})$-highest weight vector in $\Vh^{r,s}$ with
$\wt_A(v_\la)=\wtembed(\la)$. There is an embedding of virtual
$U_q(\gfin)$-crystals $\V(\la)\rightarrow \V^{r,s}$ given by
$u\mapsto v$ where $u=\PSchen(v)$. The map is a morphism of
virtual $U_q(\gfin)$-crystals since $\PSchen$ is a morphism of
$U_q(A_{2n-1})$-crystals.

Now let $\la\in\Pfin^+$ be a weight appearing in
\eqref{eq:decomp}. Then $B_A(\wtembed(\la))$ occurs in
$\Vh^{r,s}$. By the argument in the previous paragraph it suffices
to show that $v_\la\in \V^{r,s}$. It can be shown that there is an
explicit sequence of virtual operators going from $u(\Vh^{r,s})$
to $v_\la$. Moreover one has $D_{\V^{r,s}}(v_\la)=0$ if
$\Vh^{r,s}$ has a single tensor factor, and otherwise
\begin{equation*}
  D_{\V^{r,s}}(v_\la) = \dfrac{1}{\gamma'} D_{\Vh^{r,s}}(v_\la) =
  \dfrac{1}{\gamma'} (-rs+\sum_{j=1}^r j m_j) =
  (\la-s\mu_r\Lab_r|\Lab_n)=D_{B^{r,s}}(v_\la).
\end{equation*}
\end{proof}

Instead of writing down this sequence of operators in full
generality, we give the following illustrative example.

\begin{example} Let $\gggg=\Cn$, $n=4$, $r=3$, $s=5$, and
$\la=2\Lab_2+\Lab_3$. We work in $\Vo^{r,s}$ and $\Vho^{r,s}$ for
convenience. Using the explicit rules for computing the operators
of type $\A$,
\begin{equation*}
\begin{CD}
  u(\Vho^{3,5}) @>{\et_2^4\et_1^4\et_0^2}>>
v_{4\Lab_2+\Lab_3} @>{\et_1^2\et_0}>> v_{2\Lab_1+2\Lab_2+\Lab_3}
@>{\et_0}>> v_{2\Lab_2+\Lab_3}
\end{CD}
\end{equation*}
Each of the above vectors has the form $b \otimes u(B_A^{5,5})$,
where $b$ is given respectively by
\begin{equation*}
  \begin{array}{|c|c|c|c|c|} \hline
    1&1&1&1&1 \\ \hline 2&2&2&2&2 \\ \hline 3&3&3&3&3 \\ \hline
  \end{array} \quad
  \begin{array}{|c|c|c|c|c|} \hline
    1&1&1&1&1 \\ \hline 2&2&2&2&2 \\ \hline 3&6&6&6&6 \\ \hline
  \end{array} \quad
  \begin{array}{|c|c|c|c|c|} \hline
    1&1&1&1&1 \\ \hline 2&2&2&6&6 \\ \hline 3&6&6&7&7 \\ \hline
  \end{array} \quad
  \begin{array}{|c|c|c|c|c|} \hline
    1&1&1&6&6 \\ \hline 2&2&2&7&7 \\ \hline 3&6&6&8&8 \\ \hline
  \end{array}
\end{equation*}
This shows that $v_\la \in \V^{r,s}$.
\end{example}

\subsection{Characterization of $\V^{r,s}$} \label{subsec:char}
We characterize the elements of $\V^{r,s}$ inside $\Vh^{r,s}$
using the self-duality condition \eqref{eq:self dual R}.

\begin{prop} \label{pp:virtual D} Let $\gggg=\Dt$ and $b\in \Vh^{r,s}$.
Then $b\in\V^{r,s}$ if and only if \eqref{eq:self dual R} holds.
\end{prop}
\begin{proof} The forward direction holds by Proposition
\ref{pp:Brs self dual}. For the reverse direction, suppose $b\in
\Vh^{r,s}$ satisfies \eqref{eq:self dual R}. Let $\V\subseteq
\Vh^{r,s}$ be the virtual $U'_q(\gggg)$-crystal generated by $b$.
It follows from the proof of Proposition \ref{pp:Brs self dual}
that every element of $\V$ satisfies \eqref{eq:self dual R}. But
$\V$ is automatically aligned by Remark \ref{rem:aligned} and
Proposition \ref{pp:virtual eps}. By Proposition \ref{pp:align
hwv} there is a $U_q(\gfin)$-highest weight vector $u \in \V$ that
is also a $U_q(A_{2n-1})$-highest weight vector. Since $u$
satisfies \eqref{eq:self dual R}, $\wt_A(u)\in\Image(\wtembed)\cap
(\Pfin^A)^+$. Since $u\in \Vh^{r,s}$ it follows that
$\wt_A(u)=\wtembed(\la)$ for some $\la$ appearing in
\eqref{eq:decomp}. By Proposition \ref{pp:big enough} $u\in
\V^{r,s}$. It follows that $\V=\V^{r,s}$.
\end{proof}

\begin{theorem} \label{thm:D decomp} Let $\gggg=\Dt$.
Then the $U'_q(\gggg)$-crystal $\V^{r,s}$ is aligned and has the
$U_q(\gfin)$-crystal decomposition given in \eqref{eq:decomp}.
\end{theorem}
\begin{proof} This follows from the proof of Proposition
\ref{pp:virtual D}.
\end{proof}

\begin{lemma} \label{lem:dual wt} Let $b$ be an element of the
$U_q(A_{2n-1})$-crystal $B_A(\Lab_{2n-r}) \otimes B_A(\Lab_r)$
such that \eqref{eq:PSD} holds for $\PSchen(b)$. Let $c_j$ be the
number of occurrences of $j$ in $b$ for $1\le j\le 2n$. Then
$c_j+c_{2n+1-j}=2$ for all $1\le j\le n$.
\end{lemma}
\begin{proof} The lemma follows directly from
$\wt_A(b)=-w_0 \wt_A(b)$, which holds by \eqref{eq:PSD} for
$\PSchen(b)$, \eqref{eq:dual crystal} and \eqref{eq:flip crystal}.
\end{proof}

\begin{theorem} \label{thm:C decomp} Let $\gggg=\Cn$. Then
$\V^{r,1}$ is an aligned $U'_q(\Cn)$-crystal and
$\V^{r,1}\cong\V(\Lab_r)$ as a $U_q(C_n)$-crystal.
\end{theorem}
Recall that the characterization of the $U_q(C_n)$-crystal
$V(\Lab_r)$ was given in Proposition \ref{pp:virtual C}.

\begin{proof} By Proposition \ref{pp:big enough}
$\V(\Lab_r)\subseteq \V^{r,1}$ and the embedding is inclusion. To
show equality it suffices to show that $\V(\Lab_r)$ is closed
under $\ft_0$ and $\et_0$. To check this we use the explicit
computation of the $0$-string in $U'_q(\A)$-crystals given in
section \ref{sec:A perf}. Let $b=u\otimes v\in\V(\Lab_r)\subset
B_A(\Lab_{2n-r}+\Lab_r)$ such that $\ft_0(b)\not=\emptyset$. Write
$\ft_0(b)=w\otimes x\in \Vh^{r,1}$. Then $u$ and $v$ both contain
$2n$ and do not contain $1$. Computing $\ft_0=(\ft_0^A)^2$ on $b$,
in all cases $w$ and $x$ are obtained from $u$ and $v$
respectively by removing $2n$ from the bottom and adding $1$ at
the top. Thus $wx$ is a tableau. $\ft_0(b)\in\V^{r,1}$ satisfies
\eqref{eq:self dual R} by Proposition \ref{pp:Brs self dual}. Let
$w_j$ and $x_j$ be defined as $u_j$ for $j=1,2$. Then $w_2 x_2$ is
a tableau, being obtained from $u_2 v_2$ by removing $2n$ from the
bottom of each column. From this one also sees that
$\len{w_2}-\len{x_2}=n-r$. By Proposition \ref{pp:virtual C},
$\ft_0(b)\in \V(\Lab_r)$. In an entirely similar manner, one may
show that $\V(\Lab_r)$ is closed under $\et_0$.

By Theorem \ref{thm:virtual fin} it remains to show that the
elements of $\V(\Lab_r)\subseteq \Vh^{r,1}$ are $0$-aligned. By
Lemma \ref{lem:dual wt}, we need only check the case that $1$ and
$2n$ occur in $b$ once each. Since $uv$ is a tableau, $1\in u$. If
$2n\in u$ then $\epsilon_0^A(b)=\varphi_0^A(b)=0$ and $b$ is
$0$-aligned. So assume $2n\in v$. If $r<n$, then $\Vh^{r,1}$ has
two tensor factors, and by \eqref{eq:e on two factors} and
\eqref{eq:f on two factors} $b$ is $0$-aligned. $b$ is
$0$-aligned. So let $r=n$. Write $u=\widehat{u}1$ and
$v=(2n)\widehat{v}$. We will show that $\epsilon_0^A(b)=0$, as the
proof of $\varphi_0^A(b)=0$ is similar.
$\et_0^A=\psi\circ\et_{2n-1}^A\circ \psi^{-1}$ by \eqref{eq:autoA
e}. Thus it is enough to show that
$\epsilon_{2n-1}(\psi^{-1}(b))=0$. Since $\psi^{-1}(b)$ has a
single $2n-1$ and a single $2n$, this holds if $2n-1$ is in the
right hand column of $\psi^{-1}(b)$. But
\begin{equation*}
\psi^{-1}(b)|_{[2n-1]}= \PSchen(b|_{[2,2n]})=\PSchen(u
\widehat{v})=\R(u\otimes \widehat{v})
\end{equation*}
by \eqref{eq:autoA rotate}. So it suffices to show that $2n-1$
remains in the right hand column in passing from $u\otimes
\widehat{v}$ to $\R(u\otimes \widehat{v})$. But there is only one
letter that moves from the right hand column to the left under
$\R$, and it must be in $v_1$ (and therefore is not $2n-1$) since
$u_1 v_1=(uv)|_{[n]}$ is a tableau with columns of equal size.
\end{proof}

\begin{prop}\label{pp:virtual A2} Let $\gggg=\At$, $\V=\bigoplus_{k=0}^r
\Image(\iota_{\Lab_k})\subseteq \Vh^{r,1}$, $b=u\otimes v\in
\Vh^{r,1}$, and $u_j$ and $v_j$ defined as in Proposition
\ref{pp:virtual C} for $j=1,2$. Then $b\in\V$ if and only if
\begin{enumerate}
\item $b$ satisfies \eqref{eq:self dual R}.
\item $u_2v_2$ is a tableau.
\item $\len{u_2}-\len{v_2}\ge n-r$.
\end{enumerate}
\end{prop}
\begin{proof} The following are equivalent:
\begin{itemize}
\item $b$ satisfies \eqref{eq:self dual R}.
\item $\PSchen(b^{\cd\ed})=\PSchen(\R(b))$.
\item $\PSchen(b)$ satisfies \eqref{eq:PSD}.
\end{itemize}
The equivalence of the first two items follows from the fact that
$\Vh^{r,1}$ is multiplicity-free as a $U_q(A_{2n-1})$-crystal and
$\PSchen$ is a $U_q(A_{2n-1})$-crystal morphism. The equivalence
of the second and third items is a consequence of Propositions
\ref{pp:dual P} and \ref{pp:flip P}, and the definition of $\R$.

Write $\PSchen(b)=u'\otimes v'$ where $u'$ and $v'$ are column
words. Since $u$ and $v$ are column words of lengths $2n-r$ and
$r$ respectively, it follows that $u'$ and $v'$ have lengths
$2n-k$ and $k$ for some $0\le k\le r$. Moreover one can pass
between $b$ and $\PSchen(b)$ using a two column jeu de taquin. Let
$u'_j$ and $v'_j$ be defined in a manner similar to $u_j$ for
$j=1,2$.

Suppose first that $b\in \V$. Then $\PSchen(b)\in\V(\Lab_k)$ for
some $0\le k\le r$. By Proposition \ref{pp:virtual C},
$\PSchen(b)$ satisfies \eqref{eq:PSD}. By the above argument $b$
satisfies \eqref{eq:self dual R}. 
By Proposition \ref{pp:virtual C} $u'_2v'_2$ is a
tableau and $\len{u'_2}-\len{v'_2}=n-k\ge n-r$. To finish the
forward direction it is enough to show that $u'_2=u_2$ and
$v'_2=v_2$. In passing from $u'\otimes v'$ to $u\otimes v$, some
letters go from the left column to the right. Since all letters in
$v'_2$ are blocking letters of $u'_2$, the only letters that can
block $u'_1$ from moving to the right are those in $v'_1$. But by
Proposition \ref{pp:virtual C}, $\len{u'_1}-\len{v'_1}=n-k$, so
there are $n-k$ letters in $u'_1$ available to move to the right.
The smallest $r-k$ of these, actually do move. In particular
$u'_2=u_2$ and $v'_2=v_2$.

For the reverse direction, suppose $b\in\Vh^{r,1}$ satisfies the
three properties. Since $\PSchen(b)$ is a tableau, it suffices to
show that $\PSchen(b)$ satisfies the properties in Proposition
\ref{pp:virtual C} for $k$. Since $b$ satisfies \eqref{eq:self
dual R}, $\PSchen(b)$ satisfies \eqref{eq:PSD}, arguing as above. This time,
passing from $u\otimes v$ to $u'\otimes v'$, we assume that
$u_2v_2$ is a tableau. That means that all the letters of $v_2$
cannot move to the left, and so stay in their column. The left
column only gets larger, and can only get larger by letters in the
interval $[n]$. Therefore $u'_2=u_2$ and $v'_2=v_2$ again.
\end{proof}

\begin{theorem} \label{thm:A2 decomp} Let $\gggg=\At$. Then
the $U'_q(\At)$-crystal $\V^{r,1}$ is aligned, and as a
$U_q(C_n)$-crystal,
\begin{equation*}
\V^{r,1} \cong \bigoplus_{k=0}^r \V(\Lab_k).
\end{equation*}
Moreover, the $U_q(C_n)$-crystal embedding
$\iota_{\Lab_k}:\V(\Lab_k)\rightarrow \V^{r,1}$ can be computed by
$\iota_{\Lab_k}=i_{r-1,1}\circ i_{r-2,1}\circ \cdots \circ
i_{k,1}$ with $i_{r,s}$ as defined in Section
\ref{subsec:embedding}.
\end{theorem}
\begin{proof} Observe that alignedness follows from the
above $U_q(C_n)$-crystal decomposition, since $\gamma'=1$ for
$\gggg=\At$ and the virtual $U_q(C_n)$-crystals $\V(\Lab_k)$ are
$i$-aligned for $1\le i\le n$. The map $i_{r-1,1}\circ\dotsm\circ
i_{k,1}$ coincides with $\iota_{\Lab_k}$ by Proposition
\ref{pp:big enough}, since both are $U_q(C_n)$-crystal embeddings
$\V(\Lab_k)\rightarrow \Vh^{r,1}$. Let $\V$ be as in Proposition
\ref{pp:big enough}, which asserts that $\V\subseteq\V^{r,1}$. To
show equality it suffices to show that $\V$ is closed under
$\et_0$ and $\ft_0$.

Recall that for $\gggg=\At$, $\Vh^{r,1}=B_A^{2n-r,1}\otimes
B_A^{r,1}$ and $\ft_0=\ft_0^A$. Let $b=u\otimes v$ be such that
$\ft_0(b)\not=\emptyset$. Write $\ft_0(b)=w\otimes x$. Let $w_j$
and $x_j$ be defined as $u_j$ is for $j=1,2$. It will be shown
that $\ft_0(b)$ satisfies the three properties in Proposition
\ref{pp:virtual A2}. Since $b\in \V\subseteq \V^{r,1}$, $\ft_0(b)$
satisfies \eqref{eq:self dual R} by Proposition \ref{pp:Brs self
dual}.

Suppose first that $\ft_0(u\otimes v)=u\otimes\ft_0(v)$. Then
$u=w$ and $v_2=(2n)x_2$. It is easily seen that properties 2 and 3
of Proposition \ref{pp:virtual A2} hold for $\ft_0(b)$, so that
$\ft_0(b)\in \V$.

Otherwise suppose $\ft_0(u\otimes v)=\ft_0(u)\otimes v$. Then
$x=v$, $2n\in u$, $1\not\in u$, and $w$ is obtained from $u$ by
removing $2n$ from the bottom and putting $1$ at the top.
In particular $u_2=(2n)w_2$ and $v_2=x_2$. The only way that
$w_2x_2$ fails to be a tableau is if $n=r$, when the first column
is shorter than the second. Since $u_2v_2$ is a tableau with
columns of equal length, $2n\in v_2$. By Lemma \ref{lem:dual wt}
there are no ones present. But then $\ft_0(u\otimes
v)=u\otimes\ft_0(v)$, contrary to assumption. Therefore $w_2 x_2$
is a tableau.

To check property 3 of Proposition \ref{pp:virtual A2} it is
enough to show $\len{u_2}-\len{v_2}>n-r$. Suppose not, that is,
equality holds. Following the proof of Proposition \ref{pp:virtual
A2}, write $u'\otimes v'=\PSchen(u\otimes v)$ and define $u'_j$
and $v'_j$. Then $u'_1\otimes v'_1=\PSchen(u_1\otimes v_1)$. As in
the aforementioned proof it can be shown that $u'\otimes
v'\in\V(\Lab_k)$ for some $0\le k\le r$ and that $u'_2=u_2$ and
$v'_2=v_2$. By our assumption it follows that $k=r$ and hence that
$u'=u$ and $v'=v$. In particular $uv$ is a tableau. By Lemma
\ref{lem:dual wt} $1\in v$ or $2n\in v$ but not both. But $1\in v$
is not possible since $uv$ is a tableau and $1\not\in u$. And if
$2n\in v$ then $\ft_0^A(u\otimes v)=u\otimes\ft_0^A(v)$, contrary
to our assumption.

The proof that $\V$ is closed under $\et_0$, is similar.
\end{proof}

\begin{prop}\label{pp:virtual A2d} Let $\gggg=\Atd$ and $1\le r\le
n$. Recall that $\V(\mu_r\Lab_r)\subset
B_A(\Lab_{2n-r}+\Lab_r)\subset\Vh^{r,1}$. Then $b\in
\V(\mu_r\Lab_r)$ if and only if $b\in \Vh^{r,1}$, $b$ is a
tableau, and $b$ satisfies \eqref{eq:self dual R}.
\end{prop}
\begin{proof} This follows from Proposition \ref{pp:virtual B}.
\end{proof}

\begin{theorem} \label{thm:A2d decomp}  Let $\gggg=\Atd$. Then
the virtual $U'_q(\Atd)$-crystal $\V^{r,1}$ is aligned and is
isomorphic to $\V(\mu_r\Lab_r)$ as a $U_q(B_n)$-crystal.
\end{theorem}
\begin{proof} It suffices to show that
$(\V(\mu_r\Lab_r,\Vh^{r,1})$ is $0$-aligned and closed under
$\et_0=(\et_0^A)^2$ and $\ft_0=(\ft_0^A)^2$. Let
$b\in\V(\mu_r\Lab_r)$. By Proposition \ref{pp:virtual A2d} and
Lemma \ref{lem:dual wt}, $b$ is a tableau and $c_1+c_{2n}=2$ where
$c_j$ is the number of occurrences of $j$ in $b$ for $1\le j\le
2n$. Write $b=u\otimes v$ where $u$ and $v$ are column words of
type $A_{2n-1}$. Suppose first that $c_1=2$. Then
$\epsilon_0^A(b)=2$ and $\varphi_0^A(b)=0$, so that $b$ is
$0$-aligned. Writing $\et_0(b)=w\otimes x$ the column words $w$
and $x$ are obtained from $u$ and $v$ by removing $1$ from the top
and putting $2n$ at the bottom. Clearly $\et_0(b)$ is a tableau.
Since $\et_0(b)\in\V^{r,1}$ it satisfies \eqref{eq:self dual R} by
Proposition \ref{pp:Brs self dual} and is therefore in
$\V(\mu_r\Lab_r)$ by Proposition \ref{pp:virtual A2d}. An entirely
similar argument applies for the case $c_1=0$. Finally, suppose
$c_1=1$. Since $b$ is a tableau, $1\in u$. Regardless of where the
single symbol $2n$ appears in $b$, one has
$\epsilon_0^A(b)=\varphi_0^A(b)=0$, so that $b$ is again
$0$-aligned.
\end{proof}

\subsection{$B^{r,1}$ for $\gggg=C^{(1)}_n$}
The goal of this section is to establish Conjecture
\ref{conj:virtual} for $\V^{r,1}$ of type $C^{(1)}_n$. The
$U'_q(C^{(1)}_n)$-module $W^{(r)}_1$ exists and has a simple
crystal basis $B^{r,1}$ \cite{AK,Ka2}. We use the symmetry of the
Dynkin diagram $C^{(1)}_n$ to derive properties which uniquely
define the structure of the affine crystal $B^{r,1}$. Then we show
that the virtual crystal $\V^{r,1}$ satisfies these properties.

We state some generalities which must be satisfied by $B^{r,s}$
assuming it exists. Consider the automorphism of the Dynkin
diagram $C^{(1)}_n$ given by $i\mapsto n-i$ for $0\le i\le n$; it
has order two. It induces an involution on crystal bases of
$U'_q(C^{(1)}_n)$-modules. Let $P\rightarrow P$ be the linear
involution given by $\La_i\mapsto \La_{n-i}$ for $i \in I$. There
is an induced linear involution $\wtflip:\Pfin\rightarrow\Pfin$
given by $\Lab_i\mapsto\Lab_{n-i} -\Lab_n$ for all $i\in J$.
Identifying $\Pfin\cong \Z^n$ one has
\begin{equation} \label{eq:weight flip}
\wtflip(a_1,a_2,\dotsc,a_n)=(-a_n,\dotsc,-a_2,-a_1).
\end{equation}
The existence of the Dynkin diagram automorphism implies the
following result.

\begin{prop} \label{pp:C Brs}
Let $B$ be the crystal basis of an irreducible integrable
$U'_q(\Cn)$-module. Then there is an automorphism $\autoC$ of $B$
with the following properties:
\begin{align}
\label{eq:autoC inv} \autoC^2&=1 \\
\label{eq:autoC wt} \wt\circ \autoC &= \wtflip\circ \wt \\
\label{eq:autoC f} \autoC\circ \ft_i&=\ft_{n-i}\circ \autoC
\qquad\text{for all $i\in I$} \\
\label{eq:autoC e} \autoC\circ \et_i &=\et_{n-i}\circ\autoC
\qquad\text{for all $i\in I$.}
\end{align}
\end{prop}
In particular
\begin{equation} \label{eq:C zero ops}
  \ft_0 = \autoC \circ \ft_n \circ \autoC \qquad \text{and} \qquad
  \et_0 = \autoC \circ \et_n \circ \autoC.
\end{equation}
Since the $U_q(C_n)$-crystal structure of $B^{r,s}$ is prescribed
by \eqref{eq:decomp}, it only remains to determine the
automorphism $\autoC$.

\begin{remark} \label{rem:C auto} Let $J'=\{1,2,\dots,n-1\}$.
By Proposition \ref{pp:C Brs} $\autoC$ acts on the set of
$J'$-highest weight vectors in $B^{r,s}$ and is uniquely
determined by this action.
\end{remark}

Define the map $\autoC:\Vh^{r,s}\rightarrow \Vh^{r,s}$ by
\begin{equation} \label{eq:virtual C auto}
  \autoC = \begin{cases}
  \autoA^n \otimes \autoA^n & \text{if $r<n$} \\
  \autoA^n & \text{otherwise}
  \end{cases}
\end{equation}
where $\autoA$ is the order $2n$ automorphism of type
$A^{(1)}_{2n-1}$.

\begin{prop} \label{pp:virtual auto} The map $\autoC$ defined by
\eqref{eq:virtual C auto} stabilizes the virtual crystal $\V^{r,s}$
and satisfies the conditions of Proposition \ref{pp:C Brs} for
$\V^{r,s}$.
\end{prop}
\begin{proof} $\autoC$ satisfies the conditions of Proposition \ref{pp:C Brs}
for all elements of $\Vh^{r,s}$; this follows immediately from the
known properties of $\autoA$ recorded in Theorem \ref{thm:A perf}.
To see that $\autoC$ stabilizes $\V^{r,s}$, by \eqref{eq:autoC f}
and \eqref{eq:autoC e} it is enough to show that
$\autoC(u(\Vh^{r,s}))\in \V^{r,s}$. But this element is of extremal
weight and is easily seen to be in the $U_q(C_n)$-component of
$u(\Vh^{r,s})$ by direct computation.
\end{proof}

\begin{prop} \label{pp:auto col} $\V(\Lab_r)$ is stable under
$\autoC$. In particular $\V^{r,1}=\V(\Lab_r)$ as sets.
\end{prop}
\begin{proof} Let $b=u\otimes v\in \V^{r,1}$. Write $u=u_2 u_1$ and
$v=v_2 v_1$ where $u_1=u|_{[n]}$ and $u_2=u|_{[n+1,2n]}$ and
similarly for $v$. By Proposition \ref{pp:virtual C} $u_2 v_2$ is
a tableau and $\len{u_2}-\len{v_2}=n-r$.

Suppose first that $r<n$. By direct computation
$\autoC(b)=(u_1+n)(u_2-n)\otimes (v_1+n)(v_2-n)$ where $t\pm n$
means to add or subtract $n$ from each entry in the tableau $t$.
Applying Proposition \ref{pp:virtual C} again, it follows that
$\autoC(b)\in \V(\Lab_r)$ since $(u_2 v_2)-n$ has the same above
properties that $u_2 v_2$ does.

Otherwise let $r=n$. Then $\len{u_j}=\len{v_j}$ for $j=1,2$. Now
\begin{equation*}
\autoA^n(uv)|_{[n]}=\PSchen(uv|_{[n+1,2n]})-n=\PSchen(u_2 v_2)-n=
(u_2-n)(v_2-n)
\end{equation*}
which has two columns of equal length. This means that
$\autoA^n(uv)|_{[n+1,2n]}$ has two columns of equal length, so
that it is equal to its $\PSchen$ tableau. But
$\PSchen(\autoA^n(uv)|_{[n+1,2n]})=\PSchen(uv|_{[n]})+n= (u_1+n)
(v_1+n)$. Therefore the columns of $\autoC(uv)$ are given by
$(u_1+n)(u_2-n)$ and $(v_1+n)(v_2-n)$, and the previous argument
goes through.
\end{proof}

\begin{theorem}\label{thm:C} The virtual crystal
$\V^{r,1}$ is the crystal graph $B^{r,1}$ of the $U'_q(C^{(1)}_n)$-module
$W^{(r)}_1$.
\end{theorem}
\begin{proof} By \cite{AK} it is known that there is an
irreducible integrable $U'_q(C^{(1)}_n)$-module $W^{(r)}_1$ with
crystal basis which is isomorphic to $B(\Lab_r)$ as a
$U_q(C_n)$-crystal. Since the $J'$-highest weight vectors in
$B(\Lab_r)$ have distinct weights, it follows by Remark \ref{rem:C
auto} that the map $\autoC$ of Proposition \ref{pp:C Brs} for the
crystal basis of the module $W^{(r)}_1$ is uniquely determined by
\eqref{eq:autoC inv}, \eqref{eq:autoC wt}, \eqref{eq:autoC f}, and
\eqref{eq:autoC e} where only $i\in J'$ are used. By Proposition
\ref{pp:auto col}, the involution $\autoC$ of $\Vh^{r,1}$
restricts to an involution on $\V^{r,1}$. As such it satisfies the
properties of Proposition \ref{pp:C Brs}, by Proposition
\ref{pp:virtual auto}.
\end{proof}

\subsection{$B^{r,1}$ for $\gggg=\Atd$}

The following result is obtained by combining results in
\cite{JMO} and \cite{B}.

\begin{theorem} \label{thm:A2d} For $\gggg=\Atd$, the virtual crystal
$\V^{r,1}$ is isomorphic to the crystal graph $B^{r,1}$ of the
$U'_q(\Atd)$-module $W^{(r)}_1$.
\end{theorem}
\begin{proof} In \cite{JMO} a $U'_q(\Atd)$-module $W_1^{(r)}$ was
constructed. Its crystal basis $B^{r,1}$ was shown to be
isomorphic to $B(\mu_r\Lab_r)$ as a $U_q(B_n)$-crystal. This
agrees with the decomposition specified by \eqref{eq:decomp}.
Combining Theorem \ref{thm:virtual fin} of \cite{B} and the
explicit action of $\et_0$ and $\ft_0$ computed in \cite{JMO}, it
follows that $B^{r,1}\cong \V^{r,1}$ as $U'_q(\Atd)$-crystals.
\end{proof}

\subsection{$B^{r,1}$ for $\gggg=\At$}

The structure of $B^{r,1}$ of type $\At$ can be deduced from that
of the opposite Dynkin labeling $\Atd$. To distinguish various
objects defined for $\Atd$ we shall use the symbol $\dagger$.
There is an isomorphism of Dynkin diagrams $\Atd\rightarrow\At$
given by $i\mapsto n-i$ for $0\le i\le n$. This induces an
isomorphism of weight lattices $\autoAtwo:P^\dagger\rightarrow P$
given by $\autoAtwo(\La_i^\dagger)=\La_{n-i}$ for $0\le i\le n$.

\begin{theorem} \label{thm:A2} For $\gggg=\At$, the virtual crystal
$\V^{r,1}$ is isomorphic to the crystal graph $B^{r,1}$ of the
$U'_q(\At)$-module $W^{(r)}_1$.
\end{theorem}
\begin{proof} The $U'_q(\Atd)$-module $W^{(r)}_{1\dagger}$ is also
the $U'_q(\At)$-module $W^{(r)}_1$. Therefore their crystal bases
$B^{r,1}_\dagger$ and $B^{r,1}$ are equal as sets (but not as
weighted crystals). By abuse of notation we write $B^{r,1}$ for
either crystal, using $\dagger$ to distinguish the two structures
on this set. The isomorphism of Dynkin diagrams implies that
$\ft_i = \ft_{n-i}^\dagger$ on $B^{r,1}$. Therefore there is an
automorphism $\autoAtwo:B^{r,1}\rightarrow B^{r,1}$ such that
\begin{equation}\label{eq:autoA2}
\begin{split}
\autoAtwo(\ft_i(b))&=\ft_{n-i}(\autoAtwo(b))\qquad\text{for $0\le
i\le
n$} \\
\wt(\autoAtwo(b))&=\autoAtwo(\wt_\dagger(b))
\end{split}
\end{equation}
for all $b\in B^{r,1}$. We claim that $B^{r,1} \cong \V^{r,1}$ as
$U'_q(\At)$-crystals. In light of the above facts and the
$U'_q(\Atd)$-crystal isomorphism $B^{r,1}_\dagger\rightarrow
\V^{r,1}_\dagger$, it is enough to show that there is a unique
bijection $\autoAtwo:\V^{r,1}_\dagger\rightarrow\V^{r,1}$ such
that \eqref{eq:autoA2} holds for all $b\in\V^{r,1}_\dagger$. But
the map $\autoAtwo$ is unique since $\V^{r,1}_\dagger$ is
connected as a virtual $U_q(B_n)$-crystal and there is no choice
for $\autoAtwo(u(\V^{r,1}_\dagger))$ by weight considerations.
Recall that $\Vh^{r,1}=\Vh^{r,1}_\dagger=B_A^{2n-r,1}\otimes
B_A^{r,1}$. It is easily verified that the automorphism
$\autoA^n\otimes\autoA^n$ of the $U'_q(\A)$-crystal $\Vh^{r,1}$
satisfies \eqref{eq:autoA2}.
\end{proof}

\section{Fermionic formulas}
\label{sec:fermi}

In this section we prove the fermionic formulas associated with
crystals of type $\Dt$, $\At$ and $\Cn$
of the form $\V^{r_L,1}\otimes \cdots \otimes \V^{r_1,1}$
as conjectured in \cite{HKOTT,HKOTY} and conjecture a new fermionic
formula of type $\Atd$. We use the virtual crystals
discussed in section \ref{sec:three types} and results \cite{KSS} on
bijections between paths of type $A_n^{(1)}$ and rigged configurations.

\subsection{Main theorem}

Let $\V^{r,s}$ be a virtual crystal of type $\Dt$, $\At$ or $\Cn$
as defined in section \ref{subsec:vc}.
Denote by $\Psi^{r,s}:\V^{r,s}\to \Vh^{r,s}$ the inclusion of
the virtual crystal into the corresponding ambient crystal.
The image of $\Psi^{r,1}$ is given in Propositions \ref{pp:virtual C},
\ref{pp:virtual D} and \ref{pp:virtual A2} for type $\Cn$,
$\Dt$ and $\At$, respectively.
Under the bijection $\phib$ from paths of type $A_{2n-1}^{(1)}$
to rigged configurations as described in section \ref{sec:paths RC}
this image can also be explicitly characterized
in terms of rigged configurations as stated in the theorem below.
This theorem will be essential in the proof of the fermionic
formulas.

Let $R=(R_1,\ldots,R_L)$ be a sequence of rectangles where
$R_i$ has $r_i$ rows and $s_i$ columns. Define
$\V_R=\V^{r_L,s_L}\otimes \cdots \otimes \V^{r_1,s_1}$,
$\Vh_R=\Vh^{r_L,s_L}\otimes \cdots \otimes \Vh^{r_1,s_1}$,
and let $\Psi_R:\V_R\to\Vh_R$ denote the inclusion
$\Psi_R=\Psi^{r_L,s_L}\otimes \cdots \otimes \Psi^{r_1,s_1}$.
Let $\Rt$ be the sequence of rectangles corresponding
to the underlying tensor product of crystals $B_A^{r,s}$
(see \eqref{eq:ac}).

\begin{theorem}\label{thm rc prop}
Let $R=(R_1,\ldots,R_L)$ be a sequence of single columns where $R_i$
has height $r_i$ and $\V_R$ the corresponding virtual
crystal of type $C_n^{(1)}$, $A_{2n}^{(2)}$ or $D_{n+1}^{(2)}$. The image
$\Image(\phib\circ\Psi_R)$ of $\phib\circ\Psi_R:
\Path(\V_R,\cdot) \to \RC(\cdot,\Rt)$ is characterized by the
set of rigged configurations $(\nu,J)$ satisfying the properties:
\begin{itemize}
\item 1, 2, 3 for type $C_n^{(1)}$
\item 1, 3 for type $A_{2n}^{(2)}$
\item 1 for type $D_{n+1}^{(2)}$
\end{itemize}
where
\begin{enumerate}
\item[1.] $(\nu,J)^{(k)}=(\nu,J)^{(2n-k)}$, i.e., they are contragredient
self-dual.
\item[2.] All parts of $\nu^{(n)}$ are even.
\item[3.] All riggings in $(\nu,J)^{(n)}$ are even.
\end{enumerate}
\end{theorem}
The proof of Theorem \ref{thm rc prop} is given in Appendix \ref{app proof}.

\begin{conjecture}\label{conj:rc prop}
Theorem \ref{thm rc prop} holds for any sequence of rectangles $R$.
\end{conjecture}

A similar characterization of $\Image(\phib\circ\Psi_R)$ seems to
exist for type $\Atd$ also. The image of $\Psi^{r,1}$ is described
explicitly in Proposition \ref{pp:virtual A2d}.

\begin{conjecture}\label{conj:rc Atd}
Let $V_R$ be a virtual crystal of type $\Atd$. The image
$\Image(\phib\circ\Psi_R)$ of $\phib\circ\Psi_R:
\Path(\V_R,\cdot) \to \RC(\cdot,\Rt)$ is characterized by the
set of rigged configurations $(\nu,J)$ satisfying the properties:
\begin{enumerate}
\item[1.] $(\nu,J)^{(k)}=(\nu,J)^{(2n-k)}$, i.e., they are contragredient
self-dual.
\item[2.] All parts in $J^{(n,i)}$ are congruent to $i$ modulo 2.
\end{enumerate}
\end{conjecture}

We believe that a proof of this conjecture for tensor products of
single columns can be given in a similar fashion to the proof of
Theorem \ref{thm rc prop} as given in Appendix \ref{app proof}.

\subsection{Rigged configurations of type $C_n^{(1)}$}
\label{sec:RC C}

Given a sequence of rectangles $R=(R_1,\ldots,R_L)$ and a
partition $\la$,
let $\nu=(\nu^{(1)},\ldots,\nu^{(n)})$ be a sequence of partitions
with the properties
\begin{enumerate}
\item
$ |\nu^{(a)}| = -\sum_{j=1}^a \la_j + \sum_{i=1}^L s_i
\min(r_i,a)$ for $1\le a\le n$
\item
 $\nu^{(n)}$ has only even parts.
\end{enumerate}
Define vacancy numbers as follows
\begin{equation}\label{vac C}
\begin{split}
 P_i^{(a)}(\nu) &= Q_i(\nu^{(a-1)}) - 2Q_i(\nu^{(a)}) + Q_i(\nu^{(a+1)})
  + Q_i(\xi^{(a)}(R))
 \quad \text{for $1\le a<n$,} \\
 P_i^{(n)}(\nu) &= Q_i(\nu^{(n-1)}) - Q_i(\nu^{(n)}) +
  \frac{1}{2} Q_i(2\xi^{(n)}(R)).
\end{split}
\end{equation}
Denote by $\Conf_\gggg(\la,R)$ where $\gggg=C_n^{(1)}$ the set of
$(\nu^{(1)},\ldots,\nu^{(n)})$ satisfying the two properties listed
above and $P_i^{(a)}(\nu)\ge 0$ for all $1\le a\le n$ and $i\ge 0$.

The set of rigged configurations of type $\gggg=C_n^{(1)}$, denoted by
$\RC_\gggg(\la,R)$, is given by $(\nu,J)$ such that $\nu\in\Conf_\gggg(\la,R)$
and $J$ is a double sequence of partitions
$J=\{J^{(a,i)}\}_{\substack{1\le a\le n\\ i\ge 1}}$
such that the partition $J^{(a,i)}$ lies in a box of width
$P_i^{(a)}(\nu)$ and height $m_i(\nu^{(a)})$.

The cocharge for $(\nu,J)\in\RC_\gggg(\la,R)$ is defined by
\begin{equation*}
\begin{split}
\cc_\gggg(\nu,J) &= \cc_\gggg(\nu)+\sum_{a=1}^n \sum_{i\ge 1}|J^{(a,i)}| \\
\text{where}\quad \cc_\gggg(\nu) &= \sum_{i\ge 1} \Bigl(
\sum_{a=1}^{n-1} \alpha_i^{(a)}(\alpha_i^{(a)}-\alpha_i^{(a+1)})
+\frac{1}{2} \alpha_i^{(n) \;2} \Bigr).
\end{split}
\end{equation*}

\subsection{Rigged configurations of type $A_{2n}^{(2)}$}
\label{sec:RC A}
Given a sequence of rectangles $R=(R_1,\ldots,R_L)$
and a partition $\la$,
let $\nu=(\nu^{(1)},\ldots,\nu^{(n)})$ be a sequence of partitions
such that
\begin{equation*}
 |\nu^{(a)}| = -\sum_{j=1}^a \la_j + \sum_{i=1}^L s_i
\min(r_i,a) \qquad \text{for $1\le a\le n$.}
\end{equation*}
Define vacancy numbers as follows
\begin{equation}\label{vac A2}
\begin{split}
 P_i^{(a)}(\nu) &= Q_i(\nu^{(a-1)}) - 2Q_i(\nu^{(a)}) + Q_i(\nu^{(a+1)})
  + Q_i(\xi^{(a)}(R))
 \quad \text{for $1\le a<n$,} \\
 P_i^{(n)}(\nu) &= Q_i(\nu^{(n-1)}) - Q_i(\nu^{(n)}) + Q_i(\xi^{(n)}(R)).
\end{split}
\end{equation}
Denote by $\Conf_\gggg(\la,R)$ where $\gggg=A_{2n}^{(2)}$ the set of
$(\nu^{(1)},\ldots,\nu^{(n)})$ satisfying the constraint above
and $P_i^{(a)}(\nu)\ge 0$ for all $1\le a\le n$ and $i\ge 0$.

The set of rigged configurations of type $\gggg=A_{2n}^{(2)}$, denoted by
$\RC_\gggg(\la,R)$, is given by $(\nu,J)$ such that $\nu\in\Conf_\gggg(\la,R)$
and $J$ is a double sequence of partitions
$J=\{J^{(a,i)}\}_{\substack{1\le a\le n\\ i\ge 1}}$
such that the partition $J^{(a,i)}$ lies in a box of width
$P_i^{(a)}(\nu)$ and height $m_i(\nu^{(a)})$.

The cocharge for $(\nu,J)\in\RC_\gggg(\la,R)$ is defined by
\begin{equation*}
\begin{split}
\cc_\gggg(\nu,J) &= \cc_\gggg(\nu)+2\sum_{a=1}^n \sum_{i\ge 1}|J^{(a,i)}| \\
\text{where}\quad \cc_\gggg(\nu) &= \sum_{i\ge 1} \Bigl(
\sum_{a=1}^{n-1} 2\alpha_i^{(a)}(\alpha_i^{(a)}-\alpha_i^{(a+1)})
 + \alpha_i^{(n) \;2} \Bigr).
\end{split}
\end{equation*}

\subsection{Rigged configurations of type $D_{n+1}^{(2)}$}
\label{sec:RC D}

Let $\la=(\la_1,\ldots,\la_n)$ be a composition satisfying
$\la_1\ge \la_2\ge\cdots\ge \la_n\ge 0$ and either $\la_i\in \Z$ for all $i$
or $\la_i\in\Z+\frac{1}{2}$ for all $i$, and $R=(R_1,\ldots,R_L)$ a sequence
of rectangles. Then let $\nu=(\nu^{(1)},\ldots,\nu^{(n)})$ be a sequence of
partitions with the properties
\begin{equation*}
 |\nu^{(a)}| = -\sum_{j=1}^a \la_j + \sum_{\substack{i=1\\ r_i<n}}^L s_i
 \min(r_i,a)+\frac{1}{2} \sum_{\substack{i=1\\ r_i=n}}^L s_i a
 \qquad \text{for $1\le a\le n$.}
\end{equation*}
Define vacancy numbers as follows
\begin{equation}\label{vac A}
\begin{split}
 P_i^{(a)}(\nu) &= Q_i(\nu^{(a-1)}) - 2Q_i(\nu^{(a)}) + Q_i(\nu^{(a+1)})
  + Q_i(\xi^{(a)}(R))
 \quad \text{for $1\le a<n$,} \\
 P_i^{(n)}(\nu) &= 2Q_i(\nu^{(n-1)}) - 2Q_i(\nu^{(n)}) + Q_i(\xi^{(n)}(R)).
\end{split}
\end{equation}
Denote by $\Conf_\gggg(\la,R)$ where $\gggg=D_{n+1}^{(2)}$ the set of
$(\nu^{(1)},\ldots,\nu^{(n)})$ satisfying the constraint above
and $P_i^{(a)}(\nu)\ge 0$ for all $1\le a\le n$ and $i\ge 0$.

The set of rigged configurations of type $\gggg=D_{n+1}^{(2)}$, denoted by
$\RC_\gggg(\la,R)$, is given by $(\nu,J)$ such that $\nu\in\Conf_\gggg(\la,R)$
and $J$ is a double sequence of partitions
$J=\{J^{(a,i)}\}_{\substack{1\le a\le n\\ i\ge 1}}$
such that the partition $J^{(a,i)}$ lies in a box of width
$P_i^{(a)}(\nu)$ and height $m_i(\nu^{(a)})$.

The cocharge for $(\nu,J)\in\RC_\gggg(\la,R)$ is defined by
\begin{equation*}
\begin{split}
\cc_\gggg(\nu,J) &= \cc_\gggg(\nu)+\sum_{i\ge 1}(\sum_{a=1}^{n-1}
 2|J^{(a,i)}|+|J^{(n,i)}|) \\
\text{where}\quad \cc_\gggg(\nu) &= \sum_{i\ge 1} \Bigl(
\sum_{a=1}^{n-1} 2\alpha_i^{(a)}(\alpha_i^{(a)}-\alpha_i^{(a+1)})
 + \alpha_i^{(n) \;2} \Bigr).
\end{split}
\end{equation*}

\subsection{Fermionic formulas}

It follows from Theorem \ref{thm statistics} and
Proposition \ref{charge energy} that
\begin{equation} \label{coenergy rccocharge}
-E_R(b) = \cc(\phit(b)) \qquad \text{for $b\in \Path(B_R^A,\La)$.}
\end{equation}

By Theorem \ref{thm rc prop} (and Conjecture \ref{conj:rc prop}),
\eqref{coenergy rccocharge}, \eqref{eq:virtual energy}
and the definitions of rigged configurations in Sections \ref{sec:RC C},
\ref{sec:RC A} and \ref{sec:RC D} we can express the one dimensional sum
of type $\gggg=C_n^{(1)}, A_{2n}^{(2)}, D_{n+1}^{(2)}$ in terms of rigged
configurations
\begin{equation}\label{K rc C}
\X(B_R,\La;q^{-1}) = \sum_{(\nu,J)\in \RC_\gggg(\la,R)} q^{\cc_\gggg(\nu,J)}
\end{equation}
where $\la$ is the ``partition'' corresponding to the weight $\La$.
That is, if $\La=\Lab_{k_1}+\cdots+\Lab_{k_\ell}$ then $\la$ is
the ``partition'' with columns of height $k_1,\ldots,k_\ell$.
In the case of $\gggg=D_{n+1}^{(2)}$ $\Lab_n$ yields a column
of height $n$ and width $\frac{1}{2}$ (explaining the quotation
marks around partition).
The contribution of $\{J^{(a,i)}\}$ to $\cc_\gggg(\nu,J)$ can be evaluated
explicitly by noting that the generating function of partitions in a box of
dimension $m\times p$ is the $q$-binomial coefficient
\begin{equation*}
\qbin{m+p}{m}_q = \frac{(q)_{m+p}}{(q)_m(q)_p}
\end{equation*}
for $p,m\in \N$ and zero otherwise, where
$(q)_m=(1-q)(1-q^2)\cdots (1-q^m)$.
We also introduce $t_a,t_a^\vee$ ($1\le a\le n$) by
\begin{equation*}
\begin{split}
(t_a)_{1\le a\le n}&=\begin{cases}
(2,2,\ldots,2,1) & \text{for $\Cn$}\\
(1,1,\ldots,1,1) & \text{otherwise}\end{cases}\\
(t_a^\vee)_{1\le a\le n}&=\begin{cases}
(2,2,\ldots,2,1) & \text{for $\Dt$}\\
(2,2,\ldots,2,2) & \text{for $\At$}\\
(1,1,\ldots,1,1) & \text{for $\Cn$.}\end{cases}
\end{split}
\end{equation*}
Then \eqref{K rc C} may be recast into the form
\begin{equation}\label{fermi1}
\X(B_R,\La;q^{-1}) = \sum_{\nu\in\Conf_\gggg(\la,R)}
q^{\cc_\gggg(\nu)}
 \prod_{a=1}^n \prod_{i\ge 1}
 \qbin{m_i(\nu^{(a)})+P_i^{(a)}(\nu)}{m_i(\nu^{(a)})}_{q_a}
\end{equation}
where $q_a=q^{t_a^\vee}$.
This formula is known as the fermionic formula.

To compare this with the fermionic formula stated in \cite[Eq.(4.5)]{HKOTT}
several definitions are necessary. Set
\begin{align*}
p_i^{(a)}&=\begin{cases}
 P_i^{(a)}(\nu) & \text{for $1\le a<n$}\\
 P_{\gamma'\;i}^{(n)}(\nu) & \text{for $a=n$} \end{cases}\\
m_i^{(a)}&=\begin{cases}
 m_i(\nu^{(a)}) & \text{for $1\le a<n$}\\
 m_{\gamma'\;i}(\nu^{(n)}) & \text{for $a=n$} \end{cases}
\end{align*}
with $\gamma'$ as in \eqref{eq:gamma def}.
Let $L_i^{(a)}$ be the number of rectangles in $R$ with $a$ rows
and $i$ columns.
Let $\alphat_a,\Lat_a$ ($a=1,\ldots,n$) be the simple roots and fundamental
weights of $\gt$ and let $\tilde{P}$ be its weight lattice.
To relate $\tilde{P}$ with $\Pb$ we introduce a
$\Z$-linear map $\iota:\Pb\rightarrow\tilde{P}$ by
\begin{equation} \label{eq:iota}
\iota(\Lab_a)=\varepsilon_a\Lat_a\qquad\text{for $1\le a\le n$}
\end{equation}
where $\varepsilon_a$ is specified as
\[
\varepsilon_a=\left\{
\begin{array}{ll}
2\quad&\text{if $\gggg=\At$ and $a=n$}\\
1&\text{otherwise.}
\end{array}\right.
\]
Note that \eqref{eq:iota} induces $\iota(\alpha_a)=\varepsilon_a\alphat_a$.
With these definitions we obtain
\begin{equation}\label{pm}
p_i^{(a)}=\sum_{j\ge 1} L_j^{(a)} \min(i,j) - \frac1{t_a^\vee}
\sum_{b=1}^n (\alphat_a|\alphat_b) \sum_{k\ge 1} \min(t_b i, t_a k) m_k^{(b)}
\end{equation}
for all $i\ge 1$ and $1\le a \le n$,
\begin{equation}\label{m const}
\sum_{a=1}^n \sum_{i\ge 1} i m_i^{(a)} \alphat_a =
\iota\left(\sum_{a=1}^n \sum_{i\ge 1} i L_i^{(a)} \Lab_a - \La\right)
\end{equation}
and
\begin{equation}\label{eq:cc}
\cc_C(\nu)=\cc(\{m\}):=\frac{1}{2}\sum_{a,b=1}^n(\alphat_a|\alphat_b)
\sum_{j,k\ge 1} \min(t_b j,t_a k) m_j^{(a)} m_k^{(b)}.
\end{equation}

Then \eqref{fermi1} is equivalent to
\begin{equation}\label{fermi2}
\X(B_R,\La;q^{-1}) = \sum_{\{m\}} q^{\cc(\{m\})}
 \prod_{a=1}^n \prod_{i\ge 1}
 \qbin{m_i^{(a)}+p_i^{(a)}}{m_i^{(a)}}_{q_a}
\end{equation}
where the sum is over all $\{m_i^{(a)}\in\N\mid 1\le a\le n, i\ge 1\}$
subject to the constraint \eqref{m const} and $p_i^{(a)}$ is given
by \eqref{pm}. This is \cite[Eq. (4.5)]{HKOTT} with $q\to 1/q$
and $\ell\to\infty$.

There are also level-restricted versions of the one dimensional
configuration sums which are $q$-analogues of
fusion coefficients. In \cite[Eq. (4.8)]{HKOTY}
and \cite[Conj. 3.1 and Eq. (4.5)]{HKOTT} fermionic formulas
for the level-restricted configuration sums were conjectured for $\la=0$.
For type $C_n^{(1)}$ these and generalizations thereof for general $\la$
are proven in \cite{OSS} using the virtual crystal techniques developed in
this paper.

\subsection{New conjectured fermionic formula for $\Atd$}
\label{sec:fermi Atd}
Similar to the derivation of \eqref{fermi2} from Theorem \ref{thm rc prop}
and Conjecture \ref{conj:rc prop}, a fermionic formula of type 
$\Atd$ can be derived from Conjecture \ref{conj:rc Atd} which
to our knowledge is new. Namely, we conjecture that for $B_R$
a crystal of type $\Atd$ and $\La$ a dominant integral weight
\begin{multline}\label{eq:fermi Atd}
X(B_R,\La;q^{-1})=\sum_{\{m\}} q^{\cc(\{m\})}
\prod_{a=1}^{n-1} \prod_{i\ge 1} \qbin{m^{(a)}_i+p^{(a)}_i}{m^{(a)}_i}\\
\prod_{\substack{i\ge 1\\ \text{$i$ even}}} 
 \qbin{m^{(n)}_i+p^{(n)}_i}{m^{(n)}_i}
\prod_{\substack{i\ge 1\\ \text{$i$ odd}}}
 q^{m^{(n)}_i/2} \qbin{m^{(n)}_i+p^{(n)}_i-1}{m^{(n)}_i}.
\end{multline}
Here $\cc(\{m\})$ is as in \eqref{eq:cc} and
$p_i^{(a)}$ as in \eqref{pm} with 
$(t_a)_{1\le a\le n}=(t_a^\vee)_{1\le a\le n}=(1,1,\ldots,1,1)$.
The sum is over all $\{m_i^{(a)}\in\N\mid 1\le a\le n, i\ge 1\}$
subject to constraint \eqref{m const} with $\varepsilon_n=2$ 
and $\varepsilon_a=1$ for $1\le a<n$.
Note that a summand in \eqref{eq:fermi Atd} is only nonzero if 
$p_i^{(a)}\ge 0$ for all $1\le a\le n$, $i\ge 1$
and $p_i^{(n)}\ge 1$ if $m_i^{(n)}>0$ and $i$ is odd.

\appendix

\section{Proof of Theorem \ref{thm rc prop}}
\label{app proof}

First we show that if $(\nu,J)\in\Image(\phib\circ\Psi_R)$ then
$(\nu,J)$ has the properties of the theorem. Later we will
show that the converse holds, that is, if $(\nu,J)$ satisfies the
the properties then $(\nu,J)\in\Image(\phib\circ\Psi_R)$.

By Propositions \ref{pp:virtual C}, \ref{pp:virtual D} and
\ref{pp:virtual A2} the image
of $\Psi_R$ yields paths that are contragredient self-dual up to application
of the combinatorial $R$-matrix for all three types $\Dt$, $\At$ and $\Cn$.
It was shown in \cite{KSS} that
under the bijection $\phib$ the combinatorial $R$-matrix on LR
tableaux yields the identity map on rigged configurations. Hence
point 1 follows from Theorem \ref{cdual RC} in all three cases.

To prove point 3 for type $\At$ and points 2 and 3 for type $\Cn$
we proceed by induction on $L$.
Let $p=p_L\otimes \cdots \otimes p_1 \in \Path(\V_R,\cdot)$.
By induction $(\nu,J)=\phib\circ \Psi_{R'}(p_{L-1}\otimes\cdots\otimes p_1)$
satisfies point 3 for $\At$ and points 2 and 3 for $\Cn$ where
$R'=(R_1,\ldots,R_{L-1})$. Let $r'$ be the height of $R_L$.
For type $\Cn$, $p_L=u\otimes v\in\V_C(\Lab_r)$ for some one column tableaux
$u$ and $v$ of height $2n-r$ and $r$ respectively where $r=r'$.
For type $\At$, $p_L=i_{r'-1,1}\circ i_{r'-2,1}\circ\cdots i_{r,1}(u\otimes v)$
for some $u\otimes v\in \V_C(\Lab_r)$ by Theorem \ref{thm:A2 decomp}.
Let $\ell_i^{(k)}$ for $1\le i\le r$ (resp. $s_i^{(k)}$ for
$1\le i\le 2n-r$) be the lengths of the selected strings corresponding
to the $i$-th letter in $v$ (resp. $u$) under $\phib$
as described in Section \ref{sec:paths RC}. By definition and
by Lemma \ref{lem ineq} the following inequalities hold
\begin{align}
\label{prop1}
\ell_i^{(k)}&\le \ell_i^{(k+1)}, \qquad s_i^{(k)}\le s_i^{(k+1)}\\
\label{prop2}
\ell_i^{(k)}&\le \ell_{i-1}^{(k-1)}, \qquad s_i^{(k)}\le s_{i-1}^{(k-1)}.
\end{align}

Since $r\le n$, the number of boxes added to $\nu^{(n)}$ by $v$ is
$d:=r-|v|_{[n]}|$. Similarly, the boxes added to the central
partition by $u$ is $n-|u|_{[n]}|$ which is equal to $d$ by
Proposition \ref{pp:virtual C}. Let $a$ (resp. $b$) be such that
$u_a,u_{a+1},\ldots,u_{a+d-1}$ (resp.
$v_b,v_{b+1},\ldots,v_{b+d-1}$) add boxes to $\nu^{(n)}$. Since by
\eqref{prop1} and \eqref{prop2} $\ell_{i+1}^{(k)}\le \ell_i^{(k)}$
and $s_{i+1}^{(k)}\le s_i^{(k)}$ it suffices to show that
\begin{equation}\label{claim eq}
s_{a+i}^{(n)}=\ell_{b+i}^{(n)}+1 \quad  \text{for $0\le i<d$}.
\end{equation}
Hence in this case point 2 follows directly by induction for type $\Cn$.
Since the vacancy numbers $P_m^{(n)}$ are all even by the symmetry of
point 1 and the fact that under the embedding all rectangles in $\Rt$
of height $n$ have even widths, point 3 is also satisfied.
For type $\At$ equation \eqref{claim eq} also ensures that after
the addition of $u\otimes v$ all riggings are even since the
vacancy numbers are even. The step $p_L$ is obtained from $u\otimes v$
by the application of $i_{r'-1,1}\circ\cdots\circ i_{r,1}$. By
Theorem \ref{thm:emb} this induces $j_{r'-1,1}\circ\cdots\circ j_{r,1}$
on rigged configurations which adds singular strings of length 1.
Again, since the $j_{k,1}$ do not change the vacancy numbers and
all vacancy numbers are even, point 3 is satisfied for type $\At$.
It remains to prove \eqref{claim eq}.

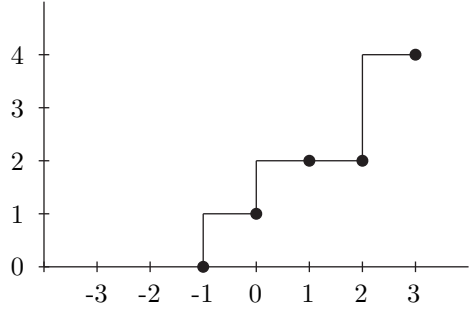
\begin{figure}
\begin{picture}(180,120)(-20,-20)
\Line(0,0)(0,100)
\Line(-2,0)(2,0)
\Line(-2,20)(2,20)
\Line(-2,40)(2,40)
\Line(-2,60)(2,60)
\Line(-2,80)(2,80)
\Text(-10,0)[]{0}
\Text(-10,20)[]{1}
\Text(-10,40)[]{2}
\Text(-10,60)[]{3}
\Text(-10,80)[]{4}
\Line(0,0)(160,0)
\Line(0,-2)(0,2)
\Line(20,-2)(20,2)
\Line(40,-2)(40,2)
\Line(60,-2)(60,2)
\Line(80,-2)(80,2)
\Line(100,-2)(100,2)
\Line(120,-2)(120,2)
\Line(140,-2)(140,2)
\Text(20,-10)[]{-3}
\Text(40,-10)[]{-2}
\Text(60,-10)[]{-1}
\Text(80,-10)[]{0}
\Text(100,-10)[]{1}
\Text(120,-10)[]{2}
\Text(140,-10)[]{3}
\CCirc(60,0){2}{Black}{Black}
\CCirc(80,20){2}{Black}{Black}
\CCirc(100,40){2}{Black}{Black}
\CCirc(120,40){2}{Black}{Black}
\CCirc(140,80){2}{Black}{Black}
\Line(60,0)(60,20)\Line(60,20)(80,20)
\Line(80,20)(80,40)\Line(80,40)(100,40)
\Line(100,40)(120,40)
\Line(120,40)(120,80)\Line(120,80)(140,80)
\end{picture}
\caption{\label{fig paths}Example of an $\ell$-slide for $n=4$ and with
$\ell_3^{(4-1)}=0$, $\ell_3^{(4+0)}=1$, $\ell_3^{(4+1)}=2$, $\ell_3^{(4+2)}=2$,
$\ell_3^{(4+3)}=4$}
\end{figure}

To prove \eqref{claim eq} we need to study the properties of
$\ell_i^{(k)}$ and $s_i^{(k)}$ more closely.
To this end it is useful to think of $\ell_i^{(k)}$ and $s_i^{(k)}$
in geometrical terms. Draw dots, called bolts, at $(j,\ell_i^{(n+j)})$ in the
$xy$-plane whenever $\ell_i^{(n+j)}$ is defined and not infinity.
Connect two bolts $(j,\ell_i^{(n+j)})$ and $(j-1,\ell_i^{(n+j-1)})$
by a horizontal line from $(j-1,\ell_i^{(n+j)})$ to $(j,\ell_i^{(n+j)})$
and a vertical line from $(j-1,\ell_i^{(n+j-1)})$ to $(j-1,\ell_i^{(n+j)})$.
For fixed $i$ call this the $i$-th $\ell$-slide.
In the same fashion define the $i$-th $s$-slide by replacing $\ell_i^{(k)}$
by $s_i^{(k)}$ everywhere. An example is given in figure \ref{fig paths}.

\begin{lemma}
$\ell$-slides (resp. $s$-slides) do not cross.
\end{lemma}
\begin{proof}
This follows directly from the conditions \eqref{prop2}.
\end{proof}

Call the set of all points $(x,y)$ of a slide $S$ with $x>0$ (resp. $x<0$)
the positive (resp. negative) part of $S$.
Let us now define a folding operation on slides by mapping the point
$(x,y)$ of a slide $S$ to $(-|x|,y)$. A bolt that changes under folding
transforms to a ``nut''.

\begin{figure}
\scalebox{1.2}{
\begin{picture}(140,120)(-20,-20)
\Line(120,0)(120,100)
\Line(120,0)(122,0)
\Line(120,20)(122,20)
\Line(120,40)(122,40)
\Line(120,60)(122,60)
\Line(120,80)(122,80)
\Text(130,0)[]{0}
\Text(130,20)[]{1}
\Text(130,40)[]{2}
\Text(130,60)[]{3}
\Text(130,80)[]{4}
\Line(0,0)(120,0)
\Line(20,-2)(20,2)
\Line(40,-2)(40,2)
\Line(60,-2)(60,2)
\Line(80,-2)(80,2)
\Line(100,-2)(100,2)
\Line(120,-2)(120,2)
\Text(20,-10)[]{-5}
\Text(40,-10)[]{-4}
\Text(60,-10)[]{-3}
\Text(80,-10)[]{-2}
\Text(100,-10)[]{-1}
\Text(120,-10)[]{0}
\CCirc(20,43){3}{Black}{Black}
\CCirc(40,60){3}{Black}{Black}
\Line(20,43)(20,60)\Line(20,60)(40,60)
\CCirc(40,43){3}{Black}{Black}
\CCirc(60,43){3}{Black}{Black}
\CCirc(80,80){3}{Black}{Black}
\Line(40,43)(60,43)
\Line(60,43)(60,80)\Line(60,80)(80,80)
\CCirc(60,37){3}{Black}{Black}
\CCirc(80,37){3}{Black}{Black}
\CCirc(100,63){3}{Black}{Black}
\Line(60,37)(80,37)
\Line(80,37)(80,63)\Line(80,63)(100,63)
\CCirc(80,0){3}{Black}{Black}
\CCirc(100,20){3}{Black}{Black}
\Line(80,0)(80,20)\Line(80,20)(100,20)
\Line(100,20)(100,57)\Line(100,57)(120,57)
\CCirc(100,0){3}{Black}{Black}
\Line(100,1)(120,1)
\SetColor{Red}
\Line(80,80)(100,80)\Line(100,80)(100,62)
\Line(100,62)(120,62)
\Line(20,80)(40,80)\Line(40,80)(40,60)
\Line(40,60)(58,60)\Line(58,60)(58,43)
\Line(60,42)(80,42)
\Line(80,40)(97,40)\Line(97,40)(97,20)
\Line(100,20)(119,20)\Line(119,20)(119,0)
\CCirc(80,80){3}{Red}{Black}
\CCirc(100,63){3}{Red}{Black}
\CCirc(120,60){3}{Red}{Red}
\CCirc(20,80){3}{Red}{Red}
\CCirc(40,60){3}{Red}{Black}
\CCirc(60,43){3}{Red}{Black}
\CCirc(80,37){3}{Red}{Black}
\CCirc(100,20){3}{Red}{Black}
\CCirc(120,0){3}{Red}{Red}
\end{picture}
}
\caption{\label{fig nut-bolt}Example of folded $\ell$-slides for
$n=6$ and $v$ the single row with entries $3,5,6,9,12$. The positive
parts are red whereas the negative parts are black. Nut-bolt pairs are denoted
by red circles with black interior.}
\end{figure}
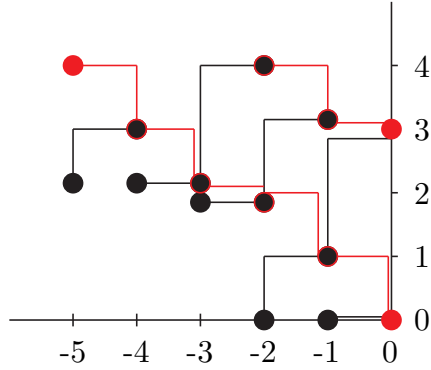

Let us now define nut-bolt pairs as follows. By the symmetry of
the rigged configuration $(\nu,J)$ and the property of $v$ that
there is no $v_k$ with $v_k=2n-v_b+1$ (otherwise the duality
property of Proposition \ref{pp:virtual C} cannot hold) the $b$-th
folded $\ell$-slide must cross each $k$-th folded $\ell$-slide
with $k<b$ and $v_k>2n-v_b$ in at least one bolt. If
$\ell_k^{(n-j)}=\ell_b^{(n+j)}$ with $j>0$ then call
$[\ell_k^{(n-j)},\ell_b^{(n+j)}]$ a nut-bolt pair (note that here
we also allow $k=b$). By induction on $i\ge 0$ the $(b+i)$-th
folded $\ell$-slide must cross each $k$-th folded $\ell$-slide
with $k<b+i$ and $v_k>2n-v_{b+i}$ in at least one unclaimed bolt.
Call $[\ell_k^{(n-j)},\ell_{b+i}^{(n+j)}]$ a nut-bolt pair if
$\ell_k^{(n-j)}=\ell_{b+i}^{(n+j)}$ and
$[\ell_k^{(n-j)},\ell_{b+i-1}^{(n+j)}]$ is not a nut-bolt pair. An
example is given in figure \ref{fig nut-bolt}.

We will use the following two lemmas to prove \eqref{claim eq}.

\begin{lemma}\label{lem nb}
For $i\ge 0$ and $j>0$ the following must be nut-bolt pairs:
\begin{align}
\label{nb 1}
&[\ell_k^{(n-j)},\ell_{b+i}^{(n+j)}] \quad \text{for some $k$}
 \qquad \text{or} \qquad
 [s_{a+i-j}^{(n-j)},\ell_{b+i}^{(n+j)}];\\
\label{nb 2}
&[s_k^{(n-j)},s_{a+i}^{(n+j)}] \quad \text{for some $k$}
 \qquad \text{or} \qquad
 [\ell_{b+i-j}^{(n-j)},s_{a+i}^{(n+j)}].
\end{align}
If $[s_{a+i-j}^{(n-j)},\ell_{b+i}^{(n+j)}]$ is not
a nut-bolt pair then
\begin{equation}\label{s ineq 1}
s_{a+i-j}^{(n-j)}>\ell_{b+i}^{(n+j)},
\end{equation}
and similarly if $[\ell_{b+i-j}^{(n-j)},s_{a+i}^{(n+j)}]$ is not
a nut-bolt pair then
\begin{equation}\label{s ineq 2}
s_{a+i}^{(n+j)}>\ell_{b+i-j}^{(n-j)}.
\end{equation}
\end{lemma}

\begin{lemma}\label{lem ess}
Let $j>1$.
If $[\ell_k^{(n-j)},\ell_{b+i}^{(n+j)}]$ is a nut-bolt pair
and $[\ell_k^{(n-j+1)},\ell_{b+i-1}^{(n+j-1)}]$ is not a nut-bolt pair,
then $\ell_{b+i}^{(n+j)}<s_{a+i-j}^{(n-j)}\le
\max\{\ell_k^{(n-j+1)},\ell_k^{(n-j)}+1\}$.
\end{lemma}

The proof of these lemmas is given by an intertwined
induction on $i$ and $j$. We will give the proofs
separately, but the inductive step of each will
rely on the validity of both other lemmas for
either the same $i$ and greater $j$, or smaller $i$.

\begin{proof}[Proof of Lemma \ref{lem nb}]
We prove \eqref{nb 1} by induction on increasing $i$ and decreasing $j$.
Suppose that for given $i$ and $j$ $[\ell_k^{(n-j)},\ell_{b+i}^{(n+j)}]$
is not a nut-bolt pair for any $k$. We want to show that then
$[s_{a+i-j}^{(n-j)},\ell_{b+i}^{(n+j)}]$ is a nut-bolt pair.

First assume that $j$ is maximal, that is $n+j=v_{b+i}-1$. Note
that from the duality and the tableau condition of $u\otimes v$ as
stated in Proposition \ref{pp:virtual C} it follows that
$u_{a+i-j}>n-j$ and $u_{a+i-j-1}\le n-j$ for $n+j=v_{b+i}-1$.
Since we know already that the final rigged configuration has to
be symmetric, i.e. satisfies condition 1 of Theorem \ref{thm rc
prop}, and since $s_k^{(n-j)} \le s_{a+i-j}^{(n-j)}$ for all
$k>a+i-j$ by \eqref{prop2} we therefore must have
$s_{a+i-j}^{(n-j)}=\ell_{b+i}^{(n+j)}$ so that
$[s_{a+i-j}^{(n-j)},\ell_{b+i}^{(n+j)}]$ forms a nut-bolt pair.

Now assume that $n+j<v_{b+i}-1$. Then by induction \eqref{nb 1} must
hold for $j+1$.

Assume that $[s_{a+i-j-1}^{(n-j-1)},\ell_{b+i}^{(n+j+1)}]$ is a
nut-bolt pair. If $[s_{a+i-j}^{(n-j+1)},\ell_{b+i-1}^{(n+j-1)}]$ is also
a nut-bolt pair then by symmetry $s_{a+i-j}^{(n-j)}=\ell_{b+i}^{(n+j)}$ so
that $[s_{a+i-j}^{(n-j)},\ell_{b+i}^{(n+j)}]$ indeed forms a nut-bolt pair.
If $[\ell_k^{(n-j+1)},\ell_{b+i-1}^{(n+j-1)}]$ is a nut-bolt pair for
some $k$ or $\ell_{b+i-1}^{(n+j-1)}=\infty$ then
$\ell_k^{(n-j+1)}>\ell_{b+i}^{(n+j+1)}$ since otherwise
$[\ell_k^{(n-j)},\ell_{b+i}^{(n+j)}]$ would be a nut-bolt pair which
contradicts our assumptions. Then again by symmetry $s_{a+i-j}^{(n-j)}
=\ell_{b+i}^{(n+j)}$ and $[s_{a+i-j}^{(n-j)},\ell_{b+i}^{(n+j)}]$ forms
a nut-bolt pair.

Otherwise assume that $[\ell_k^{(n-j-1)},\ell_{b+i}^{(n+j+1)}]$ is a nut-bolt
pair for some $k$. If $[\ell_k^{(n-j)},\ell_{b+i-1}^{(n+j)}]$ is also
a nut-bolt pair then $[\ell_{k'}^{(n-j+1)},\ell_{b+i-1}^{(n+j-1)}]$
cannot be a nut-bolt pair for any $k'$ since otherwise
$[\ell_{k'}^{(n-j)},\ell_{b+i}^{(n+j)}]$ would be a nut-bolt pair which
contradicts our assumptions. Hence
$[s_{a+i-j}^{(n-j+1)},\ell_{b+i-1}^{(n+j-1)}]$ must be a nut-bolt pair.
By the change of vacancy numbers there are no singular strings of length
$\ell_k^{(n-j-1)}<h\le \ell_k^{(n-j)}$ in the $(n-j)$-th rigged partition
after the addition of all letters up to $u_{a+i-1}$. Since
$s_{a+i-j}^{(n-j)}\le s_{a+i-j}^{(n-j+1)}=\ell_{b+i-1}^{(n+j-1)}\le
\ell_k^{(n-j)}$ this shows by symmetry that $s_{a+i-j}^{(n-j)}
=\ell_{b+i}^{(n-j)}$ and hence $[s_{a+i-j}^{(n-j)},\ell_{b+i}^{(n+j)}]$ forms
a nut-bolt pair.
If $[\ell_k^{(n-j)},\ell_{b+i-1}^{(n+j)}]$ is not a nut-bolt pair then
$\ell_k^{(n-j-1)}<s_{a+i-j-1}^{(n-j-1)}\le \ell_k^{(n-j)}$
by Lemma \ref{lem ess} (note that $\ell_k^{(n-j)}>\ell_k^{(n-j-1)}$ since
otherwise $[\ell_k^{(n-j)},\ell_{b+i}^{(n+j)}]$ would be a nut-bolt
pair which contradicts our assumptions). By the same arguments as before
there are no singular strings of length $\ell_k^{(n-j-1)}<h\le \ell_k^{(n-j)}$
in the $(n-j)$-th rigged partition after the addition of all letters up to
$u_{a+i-1}$. Using $s_{a+i-j}^{(n-j)}\le s_{a+i-j-1}^{(n-j-1)}$ \eqref{nb 1}
follows again.

The statements \eqref{nb 2}, \eqref{s ineq 1} and \eqref{s ineq 2} are
proven in a similar fashion.
\end{proof}

\begin{proof}[Proof of Lemma \ref{lem ess}]
The lower bound follows directly from \eqref{s ineq 1}.

The upper bound is proven by increasing induction on $i$
and decreasing induction on $j$.

Fix $i$. Let $j$ be maximal, that is $n+j=v_{b+i}-1$.
Assume that $[\ell_{k'}^{(n-j+1)},\ell_{b+i-1}^{(n+j-1)}]$
is a nut-bolt pair for some $k'\le k$. The conditions of the Lemma
require that $k'<k$. By definition
$\ell_{b+i}^{(n+j)}$ is the length of the largest singular string
in $(\nu,J)^{(n+j)}$ subject to the condition
$\ell_{b+i}^{(n+j)}\le \ell_{b+i-1}^{(n+j-1)}$.
There is a singular string of length
$\ell_{k'}^{(n-j)}\ge \ell_k^{(n-j)}$ in $(\nu,J)^{(n-j)}$.
Hence by symmetry $[\ell_{k'}^{(n-j)},\ell_{b+i}^{(n+j)}]$
must form a nut-bolt pair so that $k'=k$. This contradicts $k'<k$.
Therefore, by \eqref{nb 1}, either $\ell_{b+i-1}^{(n+j-1)}=\infty$
or $[s_{a+i-j}^{(n-j+1)},\ell_{b+i-1}^{(n+j-1)}]$ is a nut-bolt pair.
By the definition of $\ell_{b+i}^{(n+j)}$ there is no singular
string of length $\ell_{b+i}^{(n+j)}<h\le \ell_{b+i-1}^{(n+j-1)}$
in $(\nu,J)^{(n+j)}$ and hence by symmetry in $(\nu,J)^{(n-j)}$.
Using that $s_{a+i-j}^{(n-j)}\le s_{a+i-j}^{(n-j+1)}$ and taking into
account the change in vacancy numbers this implies that
$s_{a+i-j}^{(n-j)}\le \max\{\ell_k^{(n-j+1)},\ell_k^{(n-j)}+1\}$.

If $n+j<v_{b+i}-1$ then by \eqref{nb 1} and \eqref{nb 2} and
the crossing conditions on $\ell$-slides
exactly one of the following three must be a nut-bolt pair for
$\ell_{b+i}^{(n+j+1)}$:
(1) $[s_{a+i-j-1}^{(n-j-1)},\ell_{b+i}^{(n+j+1)}]$;
(2) $[\ell_k^{(n-j-1)},\ell_{b+i}^{(n+j+1)}]$;
(3) $[\ell_{k-1}^{(n-j-1)},\ell_{b+i}^{(n+j+1)}]$.
We are going to treat each case separately.

Assume that case (1) holds.
If $\ell_{b+i-1}^{(n+j-1)}\ge \ell_{b+i}^{(n+j+1)}$,
then by the definition of $\ell_{b+i}^{(n+j)}$ there are
no singular strings of length
$\ell_{b+i}^{(n+j)}<h\le \ell_{b+i}^{(n+j+1)}$ in $(\nu,J)^{(n+j)}$.
By symmetry there are no singular strings of length
$\ell_k^{(n-j)}<h\le \ell_{b+i}^{(n+j+1)}$ in $(\nu,J)^{(n-j)}$.
Taking into account the change of vacancy number, there are still no
singular strings of length
$\max\{\ell_k^{(n-j+1)},\ell_k^{(n-j)}+1\}<h\le \ell_{b+i}^{(n+j+1)}$
in the $(n-j)$-th rigged partition  after the addition of all letters
up to $u_{a+i-1}$. Since $s_{a+i-j}^{(n-j)}\le s_{a+i-j-1}^{(n-j-1)}=
\ell_{b+i}^{(n+j+1)}$ by \eqref{prop2} it follows that
$s_{a+i-j}^{(n-j)}\le \max\{\ell_k^{(n-j+1)},\ell_k^{(n-j)}+1\}$.
If $\ell_{b+i-1}^{(n+j-1)}< \ell_{b+i}^{(n+j+1)}$, then
$[s_{a+i-j}^{(n-j+1)},\ell_{b+i-1}^{(n+j-1)}]$ must be a nut-bolt pair.
By the same arguments as above there are no singular strings of
length $\max\{\ell_k^{(n-j+1)},\ell_k^{(n-j)}+1\}<h\le \ell_{b+i-1}^{(n+j-1)}$
in the $(n-j)$-th rigged partition  after the addition of all letters
$u_{a+i-1}$. Since $s_{a+i-j}^{(n-j)}\le s_{a+i-j}^{(n-j+1)}
=\ell_{b+i-1}^{(n+j-1)}$ by \eqref{prop1} this implies
$s_{a+i-j}^{(n-j)}\le \max\{\ell_k^{(n-j+1)},\ell_k^{(n-j)}+1\}$.

Assume that case (2) holds.
Then the conditions of Lemma \ref{lem ess} are satisfied
for $j$ replaced by $j+1$ and by induction
$s_{a+i-j-1}^{(n-j-1)}=\ell_k^{(n-j-1)}+1$. Note that
$\ell_k^{(n-j)}=\ell_k^{(n-j-1)}$ in case (2). Since
$\ell_{b+i-1}^{(n+j-1)}>\ell_k^{(n-j+1)}$ it follows that
$s_{a+i-j}^{(n-j+1)}>\ell_k^{(n-j+1)}\ge \ell_k^{(n-j)}$.
Hence by \eqref{prop2} it follows that
$s_{a+i-j}^{(n-j)}=\ell_k^{(n-j)}+1$.

Assume that case (3) holds. Let us first show that
$s_{a+i-j}^{(n-j)}\le \ell_{k-1}^{(n-j)}+1$.
If $[\ell_{k-1}^{(n-j)},\ell_{b+i-1}^{(n+j)}]$ is not a nut-bolt pair
then by induction $s_{a+i-j-1}^{(n-j-1)}\le \ell_{k-1}^{(n-j)}$
(note that $\ell_{k-1}^{(n-j)}>\ell_{k-1}^{(n-j-1)}$ since otherwise
$[\ell_{k-1}^{(n-j)},\ell_{b+i}^{(n+j)}]$ would be a nut-bolt pair
which contradicts our assumptions). Since $s_{a+i-j}^{(n-j)}\le
s_{a+i-j-1}^{(n-j-1)}$ by \eqref{prop2} it follows that
$s_{a+i-j}^{(n-j)}\le \ell_{k-1}^{(n-j)}$.
If $[\ell_{k-1}^{(n-j)},\ell_{b+i-1}^{(n+j)}]$ is a nut-bolt pair
then (a) $[s_{a+i-j}^{(n-j+1)},\ell_{b+i-1}^{(n+j-1)}]$ or
(b) $[\ell_{k-1}^{(n-j+1)},\ell_{b+i-1}^{(n+j-1)}]$ is a nut-bolt
pair. In case (a) $s_{a+i-j}^{(n-j)}\le s_{a+i-j}^{(n-j+1)}=
\ell_{b+i-1}^{(n+j-1)}\le \ell_{k-1}^{(n-j)}$ so that the assertion
holds. In case (b) the assumptions of the Lemma hold for $i$ replaced by
$i-1$ so that by induction $s_{a+i-j-1}^{(n-j)}=\ell_{k-1}^{(n-j)}+1$.
Since $s_{a+i-j}^{(n-j)}\le s_{a+i-j}^{(n-j+1)}\le s_{a+i-j-1}^{(n-j)}$
by \eqref{prop1} and \eqref{prop2} the assertion holds.

There are no singular strings of length
$\max\{\ell_{k-1}^{(n-j-1)},\ell_k^{(n-j+1)}\}<h \le \ell_{k-1}^{(n-j)}$
in the $(n-j)$-th rigged partition after the addition of all letters up
to $u_{a+i-1}$ because of the change of the vacancy numbers.
Furthermore there are no singular strings of length
$\ell_{b+i}^{(n+j)}<h\le \min\{\ell_{b+i}^{(n+j+1)},\ell_{b+i-1}^{(n+j-1)}\}$
in $(\nu,J)^{(n+j)}$ by the definition of $\ell_{b+i}^{(n+j)}$.
Hence by symmetry and change in vacancy numbers, there are no
singular strings of length $\max\{\ell_k^{(n-j+1)},\ell_k^{(n-j)}+1\}<h
\le \min\{\ell_{b+i}^{(n+j+1)},\ell_{b+i-1}^{(n+j-1)}\}$ in the
$(n-j)$-th rigged partition after the addition of all letters up to
$u_{a+i-1}$. Now assume that $s_{a+i-j}^{(n-j)}\le \ell_{k-1}^{(n-j)}$.
Then these conditions and the fact that
$[s_{a+i-j}^{(n-j+1)},\ell_{b+i-1}^{(n+j-1)}]$ is a nut-bolt
pair if $\ell_{b+i-1}^{(n+j-1)}\le \ell_{b+i}^{(n+j+1)}$ imply
that $s_{a+i-j}^{(n-j)}\le \max\{\ell_k^{(n-j+1)},\ell_k^{(n-j)}+1\}$.

Hence the only case left to consider is
$s_{a+i-j}^{(n-j)}=\ell_{k-1}^{(n-j)}+1$. By the above arguments
this case can only occur when $[\ell_{k-1}^{(n-j)},\ell_{b+i-1}^{(n+j)}]$
and $[\ell_{k-1}^{(n-j+1)},\ell_{b+i-1}^{(n+j-1)}]$ are nut-bolt pairs
with $s_{a+i-j-1}^{(n-j)}=\ell_{k-1}^{(n-j)}+1$.
We will show that $s_{a+i-j}^{(n-j)}=\ell_{k-1}^{(n-j)}+1$ cannot occur
unless $\ell_k^{(n-j)}=\ell_{k-1}^{(n-j)}$ in which case
Lemma \ref{lem ess} obviously holds.
Assume that $\ell_k^{(n-j)}<\ell_{k-1}^{(n-j)}$.
By the definition of $\ell_{b+i-1}^{(n+j)}$ and symmetry there is no
singular string of length $\ell_{k-1}^{(n-j)}+1$ available for
$s_{a+i-j}^{(n-j)}$ unless (i) $\ell_{b+i-1}^{(n+j+1)}=\ell_{b+i-1}^{(n+j)}$
or (ii) $\ell_{b+i-2}^{(n+j-1)}=\ell_{b+i-1}^{(n+j)}$.
Assume case (i) holds. By \eqref{nb 1} either
$[s_{a+i-j-2}^{(n-j-1)},\ell_{b+i-1}^{(n+j+1)}]$ or
$[\ell_{k-2}^{(n-j-1)},\ell_{b+i-1}^{(n+j+1)}]$ has to be
a nut-bolt pair. However, the first case yields a contradiction
since $s_{a+i-j-2}^{(n-j-1)}\ge s_{a+i-j-1}^{(n-j)}$
so that $[\ell_{k-2}^{(n-j-1)},\ell_{b+i-1}^{(n+j+1)}]$ must be
a nut-bolt pair.
Assume that case (ii) holds. Then
$[\ell_{k-2}^{(n-j+1)},\ell_{b+i-2}^{(n+j-1)}]$ must be a nut-bolt pair.
By \eqref{prop2} this implies that
$\ell_{k-2}^{(n-j-1)}=\ell_{k-2}^{(n-j)}=\ell_{k-1}^{(n-j)}$ and
$[\ell_{k-2}^{(n-j-1)},\ell_{b+i-1}^{(n+j+1)}]$ must be a nut-bolt pair.
Hence in both case (i) and case (ii)
$[\ell_{k-2}^{(n-j-1)},\ell_{b+i-1}^{(n+j+1)}]$ must be a nut-bolt pair
and $\ell_{k-2}^{(n-j-1)}=\ell_{k-1}^{(n-j)}$.
Now, there is no singular string of length
$\ell_{k-1}^{(n-j)}+1$ for $s_{a+i-j}^{(n-j)}$ unless
$\ell_{k-2}^{(n-j-1)}=\ell_{k-2}^{(n-j)}$.
But then $[\ell_{k-2}^{(n-j)},\ell_{b+i-2}^{(n+j)}]$ must
be a nut-bolt pair since otherwise
$[\ell_{k-2}^{(n-j)},\ell_{b+i-1}^{(n+j)}]$ is a nut-bolt pair which
contradicts our assumptions. This in turn implies that
$\ell_{b+i-2}^{(n+j-1)}=\ell_{k-1}^{(n-j)}$. By \eqref{nb 1}
either $[s_{a+i-j-1}^{(n-j+1)},\ell_{b+i-2}^{(n+j-1)}]$
or $[\ell_{k-2}^{(n-j+1)},\ell_{b+i-2}^{(n+j-1)}]$ must be a
nut-bolt pair. However, the first case contradicts
$s_{a+i-j}^{(n-j)}=\ell_{k-1}^{(n-j)}+1$ since
$s_{a+i-j}^{(n-j)}\le s_{a+i-j-1}^{(n-j+1)}$.
Hence $[\ell_{k-2}^{(n-j+1)},\ell_{b+i-2}^{(n+j-1)}]$ must be a
nut-bolt pair. By induction $s_{a+i-j-2}^{(n-j)}=
\ell_{k-1}^{(n-j)}+1$.
In summary both $[\ell_{k-2}^{(n-j)},\ell_{b+i-2}^{(n+j)}]$
and $[\ell_{k-2}^{(n-j+1)},\ell_{b+i-2}^{(n+j-1)}]$ are nut-bolt pairs
with $s_{a+i-j-2}^{(n-j)}=\ell_{k-1}^{(n-j)}+1$.
Now repeat the entire argument which shows that
both $[\ell_{k-3}^{(n-j)},\ell_{b+i-3}^{(n+j)}]$
and $[\ell_{k-3}^{(n-j+1)},\ell_{b+i-3}^{(n+j-1)}]$ are nut-bolt pairs
with $s_{a+i-j-3}^{(n-j)}=\ell_{k-1}^{(n-j)}+1$ etc.
Since there are only finitely many $\ell$-slides these conditions
must eventually break down which shows that
$s_{a+i-j}^{(n-j)}=\ell_{k-1}^{(n-j)}+1$ cannot occur when
$\ell_k^{(n-j)}<\ell_{k-1}^{(n-j)}$. This concludes the proof of
Lemma \ref{lem ess}.
\end{proof}

\begin{proof}[Proof of \eqref{claim eq}]
We prove \eqref{claim eq} by induction on $i$.

\subsection*{Induction beginning: $i=0$}
If $v_b=n+1$ then by Proposition \ref{pp:virtual C} also $u_a=n+1$
and $u_{a-1},v_{b-1}<n$. This implies that $\ell_b^{(n)}$ is the
largest singular string in $(\nu,J)^{(n)}$ and $s_a^{(n)}$ is the
largest singular string in the central rigged partition after the
addition of $v$. Since $\ell_b^{(n)}+1$ is singular it follows
that $s_a^{(n)}\ge \ell_b^{(n)}+1$. After the addition of $v$ the
vacancy numbers of the strings of length $h>\ell_b^{(n)}$ in the
central rigged partition decrease by one. However, since all
labels are even no string of length $h>\ell_b^{(n)}+1$ becomes
singular. Hence $s_a^{(n)}=\ell_b^{(n)}+1$.

Now assume that $v_b>n+1$.

If $[s_{a-1}^{(n-1)},\ell_b^{(n+1)}]$ is a nut-bolt pair then
by \eqref{prop2} $s_a^{(n)}\le s_{a-1}^{(n-1)}=\ell_b^{(n+1)}$.
Furthermore $\ell_{b-1}^{(n-1)}>\ell_b^{(n+1)}$ since otherwise
the folding of the $b$-th $\ell$-slide would not cross the $(b-1)$-th
$\ell$-slide in a bolt. This implies that $s_a^{(n+1)}>\ell_b^{(n+1)}$.
Also $s_{a-1}^{(n-1)}>\ell_b^{(n)}$ since otherwise
$[s_{a-1}^{(n-1)},\ell_b^{(n+1)}]$ was not a nut-bolt pair.
Since there is a singular string of length $\ell_b^{(n)}+1$ this implies
that $s_a^{(n)}\ge \ell_b^{(n)}+1$.
After the addition of $v$, the vacancy numbers corresponding to the strings
of length $\ell_b^{(n)}<h\le \ell_b^{(n+1)}$ in the central rigged partition
have decreased by 1. Again, since there are no singular strings
of length $\ell_b^{(n)}<h\le \ell_b^{(n+1)}$ in $(\nu,J)^{(n)}$
by the definition of $\ell_b^{(n)}$ and all labels are even,
it follows that there are no singular strings of length
$\ell_b^{(n)}+1<h\le \ell_b^{(n+1)}$ after the addition of $v$.
Hence $s_a^{(n)}=\ell_b^{(n)}+1$.

If $[\ell_{b-1}^{(n-1)},\ell_b^{(n+1)}]$ is a nut-bolt pair
then by \eqref{s ineq 1} and \eqref{s ineq 2} we have
$s_{a-1}^{(n-1)}>\ell_{b-1}^{(n-1)}$ and $s_a^{(n+1)}>\ell_b^{(n+1)}$.
This implies $s_a^{(n)}\ge \ell_b^{(n)}+1$. After the addition of $v$
the vacancy number in the central rigged partition
corresponding to the strings of length $\ell_b^{(n)}<h\le \ell_b^{(n+1)}$
decreases by one and of length $\ell_b^{(n+1)}<h$ increases by one.
Hence there are no singular strings of length $h>\ell_b^{(n+1)}$ and
by the now familiar arguments also not for $\ell_b^{(n)}+1<h\le \ell_b^{(n+1)}$
since all labels are even. This proves $s_a^{(n)}=\ell_b^{(n)}+1$.

Finally let $[\ell_b^{(n-1)},\ell_b^{(n+1)}]$ be a nut-bolt pair.
Note that this case can only happen if $u_{a-1}=n$.
We will show that $s_{a-1}^{(n-1)}=\ell_b^{(n)}+1$. Since by
\eqref{s ineq 2} $s_a^{(n+1)}>\ell_b^{(n)}$ it then follows that
$s_a^{(n)}=\ell_b^{(n)}+1$.
If $v_b=n+2$ then $v_{b-1}<n$ for $[\ell_b^{(n-1)},\ell_b^{(n+1)}]$ to hold.
The length of the longest singular string in $(\nu,J)^{(n+1)}$
is $\ell_b^{(n)}$ and hence by symmetry also in $(\nu,J)^{(n-1)}$.
This forces $s_{a-1}^{(n-1)}=\ell_b^{(n)}+1$.
Hence assume $v_b>n+2$.
If $[\ell_b^{(n-2)},\ell_b^{(n+2)}]$ is a nut-bolt pair then
$s_{a-1}^{(n-1)}=\ell_b^{(n)}+1$ holds by Lemma \ref{lem ess}.
If $[\ell_{b-1}^{(n-2)},\ell_b^{(n+2)}]$ is a nut-bolt pair
then the vacancy numbers in the $(n-1)$-st rigged partition corresponding
to strings of length $\ell_{b-1}^{(n-2)}<h\le \ell_{b-1}^{(n-1)}$ are
increased by one after the addition of $v$ so that there are no singular
strings of this length.
By Lemma \ref{lem ess} $\ell_{b-1}^{(n-2)}<s_{a-2}^{(n-2)}\le
\ell_{b-1}^{(n-1)}$ (note that the case $\ell_{b-1}^{(n-1)}=\infty$ is
included here). Hence using $s_{a-1}^{(n-1)}\le s_{a-2}^{(n-2)}$
and \eqref{s ineq 1} this implies $\ell_b^{(n-1)}<s_{a-1}^{(n-1)}
\le \ell_{b-1}^{(n-2)}$. By the definition of $\ell_b^{(n+1)}$ there
are no singular strings in $(\nu,J)^{(n+1)}$ of length
$\ell_b^{(n)}<h\le \ell_{b-1}^{(n-2)}$ and hence by symmetry
also not in $(\nu,J)^{(n-1)}$. After the addition of $v$ there is
a singular string of length $\ell_b^{(n-1)}+1$ so that
$s_{a-1}^{(n-1)}=\ell_b^{(n)}+1$.
To conclude assume that $[s_{a-2}^{(n-2)},\ell_b^{(n+2)}]$
is a nut-bolt pair. Then by \eqref{prop2} and \eqref{s ineq 1}
$\ell_b^{(n-1)}<s_{a-1}^{(n-1)}\le s_{a-2}^{(n-2)}$. By the same
arguments as in the previous case there are no singular strings
of length $\ell_{b-1}^{(n-1)}+1<h\le s_{a-2}^{(n-2)}$ in the
$(n-1)$-st rigged partition after the addition of $v$ so that
$s_{a-1}^{(n-1)}=\ell_b^{(n)}+1$ as asserted.

\subsection*{Induction step: $i-1\to i$}
Assume that $[s_{a+i-1}^{(n-1)},\ell_{b+i}^{(n+1)}]$ is a nut-bolt pair.
Then $\ell_{b+i}^{(n+1)}>\ell_{b+i}^{(n)}$,
$s_{a+i}^{(n+1)}\ge \ell_{b+i}^{(n+1)}$ and
$\ell_{b+i-1}^{(n-1)}>\ell_{b+i}^{(n+1)}$. Since there is a singular string
of length $\ell_{b+i}^{(n)}$ after the addition of $v$ and the letters
$u_1,\ldots,u_{a+i-1}$ it follows that
$\ell_{b+i}^{(n)}<s_{a+i}^{(n)}\le s_{a+i-1}^{(n-1)}$ by \eqref{prop2}.
By the definition of $\ell_{b+i}^{(n)}$ there are no singular strings
of length $\ell_{b+i}^{(n)}<h\le \ell_{b+i}^{(n+1)}$ in $(\nu,J)^{(n)}$.
The vacancy numbers of the central rigged partition corresponding
to the strings of length $\ell_{b+i}^{(n)}<h\le \ell_{b+i}^{(n+1)}$
decrease by one. Hence due to the even labels there are no singular
strings of length  $\ell_{b+i}^{(n)}+1<h\le \ell_{b+i}^{(n+1)}$
in the central rigged partition after the addition of all letters up to
$u_{a+i-1}$. This proves \eqref{claim eq}.

Now assume that $[\ell_{b+i-1}^{(n-1)},\ell_{b+i}^{(n+1)}]$
is a nut-bolt pair. By induction $s_{a+i-1}^{(n)}=\ell_{b+i-1}^{(n)}+1$.
The vacancy numbers of the central rigged partition after the addition
of all letter up to $u_{a+i-1}$ for the strings of length
$\ell_{b+i-1}^{(n-1)}<h\le \ell_{b+i-1}^{(n)}$ increase by one and
for the strings of length $\ell_{b+i}^{(n)}<h\le \ell_{b+i-1}^{(n-1)}$
decrease by one. Since there are no singular strings of length
$\ell_{b+i}^{(n)}<h\le \ell_{b+i}^{(n+1)}=\ell_{b+i-1}^{(n-1)}$ in
$(\nu,J)^{(n)}$ by the definition of $\ell_{b+i}^{(n)}$ there are no
singular strings of length $\ell_{b+i}^{(n)}+1<h\le \ell_{b+i-1}^{(n)}$
in the central rigged partition after the addition of all letters up to
$u_{a+i-1}$ due to the even labels. Also $s_{a+i}^{(n)}\le s_{a+i-1}^{(n)}
=\ell_{b+i-1}^{(n)}+1$. Hence \eqref{claim eq} holds unless
$s_{a+i}^{(n)}=s_{a+i-1}^{(n)}$ and $\ell_{b+i}^{(n)}<\ell_{b+i-1}^{(n)}$.
We will show that in the present case this cannot occur.
Assume that $s_{a+i}^{(n)}=s_{a+i-1}^{(n)}$
and $\ell_{b+i}^{(n)}<\ell_{b+i-1}^{(n)}$.
This requires $s_{a+i-1}^{(n-1)}=s_{a+i-1}^{(n)}$ and
$s_{a+i}^{(n+1)}=s_{a+i-1}^{(n)}$ by \eqref{prop2}.
By \eqref{nb 1} either (1) $[s_{a+i-2}^{(n-1)},\ell_{b+i-1}^{(n+1)}]$
or (2) $[\ell_k^{(n-1)},\ell_{b+i-1}^{(n+1)}]$ forms a nut-bolt pair.
Since $[\ell_{b+i-1}^{(n-1)},\ell_{b+i}^{(n+1)}]$ is a nut-bolt pair
by assumption, $k$ in case (2) must be $k=b+i-2$. Hence if case (2) holds
there is no singular string of length $s_{a+i-1}^{(n)}$ in the central rigged
partition after the addition of all letter up to $u_{a+i-1}$
by the change in vacancy number. Hence $s_{a+i}^{(n)}=s_{a+i-1}^{(n)}$
cannot hold.
By the same arguments there is no singular string of length
$s_{a+i-1}^{(n)}$ if $\ell_{b+i-2}^{(n-1)},\ell_{b+i-1}^{(n+1)}
>\ell_{b+i-1}^{(n)}$.
Hence either (a) $\ell_{b+i-2}^{(n-1)}=\ell_{b+i-1}^{(n)}$
or (b) $\ell_{b+i-1}^{(n+1)}=\ell_{b+i-1}^{(n)}$.
Since $s_{a+i-2}^{(n-1)}\ge s_{a+i-1}^{(n-1)}>\ell_{b+i-1}^{(n)}$
case (b) requires case (2) which we already argued yields a contradiction.
Hence assume case (a) and case (1) holds. But then by \eqref{nb 1}
and \eqref{nb 2} $[\ell_{b+i-2}^{(n-1)},s_{a+i-1}^{(n+1)}]$ must be a nut-bolt
pair which contradicts $s_{a+i-1}^{(n+1)}=\ell_{b+i-1}^{(n)}+1$.
This concludes the proof of \eqref{claim eq} when
$[\ell_{b+i-1}^{(n-1)},\ell_{b+i}^{(n+1)}]$ is a nut-bolt pair.

Finally assume that $[\ell_{b+i}^{(n-1)},\ell_{b+i}^{(n+1)}]$
is a nut-bolt pair.
If $\ell_{b+i}^{(n+2)}=\infty$ then $i=1$ and $v_b=u_a=n+1$,
$v_{b+1}=u_{a+1}=n+2$ and $u_{a-1},v_{b-1}<n-1$.
By symmetry $[\ell_b^{(n-1)},\ell_{b+1}^{(n+1)}]$ forms a
nut-bolt pair which contradicts our assumption.
Hence we may assume from now on that $\ell_{b+i}^{(n+2)}<\infty$.
We will show that $s_{a+i-1}^{(n-1)}=\ell_{b+i}^{(n)}+1$.
Then by \eqref{s ineq 1} and \eqref{prop2} it follows that
$s_{a+i}^{(n)}=\ell_{b+i}^{(n)}+1$.
If $[s_{a+i-2}^{(n-2)},\ell_{b+i}^{(n+2)}]$ is a nut-bolt pair
then $\ell_{b+i-1}^{(n-2)}\ge \ell_{b+i}^{(n+2)}$ since otherwise
the $(b+i-1)$-st and the $(b+i)$-th $\ell$-slide
would not cross in a bolt. Hence by the definition of
$\ell_{b+i}^{(n+1)}$ there are no singular strings of length
$\ell_{b+i}^{(n+1)}<h\le \ell_{b+i}^{(n+2)}$ in $(\nu,J)^{(n+1)}$.
By symmetry there are no singular strings of length
$\ell_{b+i}^{(n+1)}<h\le \ell_{b+i}^{(n+2)}$ in $(\nu,J)^{(n-1)}$.
Since $s_{a+i-1}^{(n-1)}\le s_{a+i-2}^{(n-2)}$ by \eqref{prop2}
it follows that $s_{a+i-1}^{(n-1)}=\ell_{b+i}^{(n)}+1$.
If $[\ell_{b+i}^{(n-2)},\ell_{b+i}^{(n+2)}]$ is a nut-bolt pair
then $s_{a+i-1}^{(n-1)}=\ell_{b+i}^{(n)}+1$ follows from
Lemma \ref{lem ess}.
The last case to consider is that
$[\ell_{b+i-1}^{(n-2)},\ell_{b+i}^{(n+2)}]$ is a nut-bolt pair.
Note that due to the change of vacancy numbers there are no singular
strings of length $\ell_{b+i-1}^{(n-2)}<h\le \ell_{b+i-1}^{(n-1)}$
in the $(n-1)$-st rigged partition after the addition of all
letters up to $u_{a+i-1}$. Also by the definition of
$\ell_{b+i}^{(n+1)}$ and symmetry there are no singular strings
of length $\ell_{b+i}^{(n-1)}+1<h\le \ell_{b+i-1}^{(n-2)}$ in the
same rigged partition. Hence, if $s_{a+i-1}^{(n-1)}\le \ell_{b+i-1}^{(n-1)}$
then $s_{a+i-1}^{(n-1)}=\ell_{b+i}^{(n-1)}+1$ as asserted.
We may therefore restrict our attention to the case
$s_{a+i-1}^{(n-1)}=\ell_{b+i-1}^{(n-1)}+1$ with
$\ell_{b+i}^{(n-1)}<\ell_{b+i-1}^{(n-1)}$ (note that
$\ell_{b+i-1}^{(n-1)}=\ell_{b+i}^{(n-1)}$ would imply the assertion).
We will show that these conditions cannot hold simultaneously.
By Lemma \ref{lem ess} $s_{a+i-1}^{(n-1)}=\ell_{b+i-1}^{(n-1)}+1$
can only happen if (1) $[\ell_{b+i-1}^{(n-1)},\ell_{b+i-1}^{(n+1)}]$
is a nut-bolt pair or (2) $\ell_{b+i-1}^{(n-2)}=\ell_{b+i-1}^{(n-1)}$.
Note that case (2) implies case (1) since otherwise
$[\ell_{b+i-1}^{(n-1)},\ell_{b+i}^{(n+1)}]$ is a nut-bolt pair which
contradicts our assumptions. Case (1) implies by induction that
$s_{a+i-2}^{(n-1)}=\ell_{b+i-1}^{(n)}+1$.
In summary $[\ell_{b+i-1}^{(n-1)},\ell_{b+i-1}^{(n+1)}]$ is a
nut-bolt pair with $s_{a+i-2}^{(n-1)}=\ell_{b+i-1}^{(n-1)}+1$.
There are no singular strings of length $\ell_{b+i-1}^{(n-1)}+1$
available in the $(n-1)$-st rigged partition for $s_{a+i-1}^{(n-1)}$
unless (a) $\ell_{b+i-1}^{(n+2)}=\ell_{b+i-1}^{(n-1)}$ or
(b) $\ell_{b+i-2}^{(n)}=\ell_{b+i-1}^{(n-1)}$. In case (a) by \eqref{nb 1}
either $[s_{a+i-3}^{(n-2)},\ell_{b+i-1}^{(n+2)}]$ or
$[\ell_{b+i-2}^{(n-2)},\ell_{b+i-1}^{(n+2)}]$ forms a nut-bolt pair.
The first case contradicts $s_{a+i-3}^{(n-2)}\ge s_{a+i-2}^{(n-1)}$.
In case (b) $\ell_{b+i-2}^{(n-2)}=\ell_{b+i-1}^{(n-1)}$ so that
$[\ell_{b+i-2}^{(n-2)},\ell_{b+i-1}^{(n+2)}]$ must form a nut-bolt pair.
Hence in both case (a) and (b) $[\ell_{b+i-2}^{(n-2)},\ell_{b+i-1}^{(n+2)}]$
forms a nut-bolt pair and $\ell_{b+i-2}^{(n-2)}=\ell_{b+i-1}^{(n-1)}$.
By the change in vacancy numbers there is no singular string of length
$\ell_{b+i-1}^{(n-1)}+1$ available for $s_{a+i-1}^{(n-1)}$ unless
$\ell_{b+i-2}^{(n-1)}=\ell_{b+i-1}^{(n-1)}$.
Then $[\ell_{b+i-2}^{(n-1)},\ell_{b+i-2}^{(n+1)}]$ must be a nut-bolt pair
since otherwise $[\ell_{b+i-2}^{(n-1)},\ell_{b+i-1}^{(n+1)}]$ forms
a nut-bolt pair which contradicts our assumptions.
By induction $s_{a+i-3}^{(n-1)}=\ell_{b+i-1}^{(n-1)}+1$.
Repeating the same arguments shows that
$[\ell_{b+i-3}^{(n-1)},\ell_{b+i-3}^{(n+1)}]$ forms
a nut-bolt pair with $s_{a+i-4}^{(n-1)}=\ell_{b+i-1}^{(n-1)}+1$ etc.
Since there are only finitely many $\ell$-slides these conditions
must be eventually broken which shows that
$s_{a+i-1}^{(n-1)}=\ell_{b+i-1}^{(n-1)}+1$ is not possible
unless $\ell_{b+i}^{(n-1)}=\ell_{b+i-1}^{(n-1)}$.
This concludes the proof of \eqref{claim eq}.
\end{proof}

So far we have shown that the conditions of Theorem \ref{thm rc prop}
are satisfied if $(\nu,J)\in\Image(\phib\circ\Psi_R)$. It
remains to show the converse, that is, if $(\nu,J)$ satisfies the
conditions of Theorem \ref{thm rc prop} then
$(\nu,J)\in\Image(\phib\circ\Psi_R)$.

Let $p=p_L\otimes\cdots\otimes p_1=\phib^{\;-1}(\nu,J)$ where
$p_i=b_i\bt_i\in B_A^{n,2}$ if $R_i$ is a single column of height
$r_i=n$ and $p_i$ is of type $\Cn$, $p_i=b_i\in B_A^{n,1}$ if
$R_i$ is a single column of height $r_i=n$ and $p_i$ of type $\Dt$, or
$p_i=b_i\otimes \bt_i\in B_A^{2n-r_i,1}\otimes B_A^{r_i,1}$ otherwise.
By Theorem \ref{cdual RC} and the fact
that the combinatorial $R$-matrix on paths yields the identity map on
rigged configurations, condition 1 of Theorem \ref{thm rc prop}
ensures that $p$ is contragredient self-dual, that is $\R(b_i
\otimes \bt_i)=b_i^{\cd\ed} \otimes \bt_i^{\cd\ed}$.
For type $\Dt$ this already shows that $(\nu,J)\in\Image(\phib\circ\Psi_R)$
by Proposition \ref{pp:virtual D}.

For type $\At$ act on $(\nu,J)$ with $j_{r,1}^{-1}\circ j_{r+1,1}^{-1}\circ
\cdots \circ j_{r_L-1,1}^{-1}$ with $r$ as small as possible. That
is, first remove singular strings of length 1 from the $k$-th
rigged partition for $r_L\le k\le 2n-r_L$ if possible. Repeat this
for $r_L-1,r_L-2,\ldots$ until no more such singular strings
can be removed for $r-1$.

By induction we may assume that $(\nub,\Jb)=\phib(p_{L-1}
\otimes\cdots \otimes p_1)$ is in $\Image(\phib\circ\Psi_{R'})$
if it satisfies the properties of Theorem \ref{thm rc prop}
where $R'=(R_1,\ldots,R_{L-1})$.
Set $u\otimes v =b_L\otimes \bt_L$ and $r=r_L$  for type $\Cn$
and $u\otimes v= i_{r,1}^{-1}\circ i_{r+1,1}^{-1}\circ \cdots \circ
i_{r_L-1,1}^{-1}(b_L\otimes \bt_L)$ for type $\At$.
Then by Proposition \ref{pp:virtual C} and Theorem \ref{thm:A2 decomp}
it suffices to prove that $(uv)|_{[n]}$ is a tableau,
$|u|_{[n]}|-|v|_{[n]}|=n-r$ and $(\nub,\Jb)$ satisfies point 3 of
Theorem \ref{thm rc prop} for type $\At$ and points 2 and 3 for type $\Cn$.

The rigged configuration $(\nub,\Jb)$ is obtained from
$(\nu,J)$ by the inverse algorithm to the one described in
Section \ref{sec:paths RC}. Let $(\nu,J)_{(u,i)}$ for
$0\le i\le 2n-r$ and $(\nu,J)_{(v,i)}$ for $0\le i\le r$
be rigged configurations where $(\nu,J)_{(u,2n-r)}=(\nu,J)$,
$(\nu,J)_{(u,0)}=(\nu,J)_{(v,r)}$ and $(\nu,J)_{(v,0)}=(\nub,\Jb)$.
The rigged configuration $(\nu,J)_{(u,i-1)}$ is obtained from
$(\nu,J)_{(u,i)}$ in the following way:
Select singular strings of length $\st_i^{(i)}, \st_i^{(i+1)},
\ldots$ in the $i$-th, $(i+1)$-th, \ldots rigged partition in
$(\nu,J)_{(u,i)}$ minimal such that
$\st_i^{(i)}\le \st_i^{(i+1)}\le \ldots$. If no such
singular string exists in the $k$-th rigged partition set
$\st_i^{(k)}=\infty$. Let $u_i$ be minimal such that
$\st_i^{(u_i)}=\infty$. Then remove a box from each
of the selected strings in the $i$-th to $(u_i-1)$-th
rigged partition in $(\nu,J)_{(u,i)}$ keeping them singular
and leaving all other strings unchanged. The result is
$(\nu,J)_{(u,i-1)}$.
Define $\lt_i^{(k)}$, $v_i$ and $(\nu,J)_{(v,i-1)}$ in the
analogous fashion starting with $(\nu,J)_{(v,i)}$.
Then let $u$ be the column tableau of height $2n-r$
with entries $u_1,\ldots,u_{2n-r}$ and $v$ the column
tableau of height $r$ with entries $v_1,\ldots,v_r$.

Similarly to \eqref{claim eq} the following holds.
Let $a$ and $d_u$ be defined such that $\st_a^{(n)},
\st_{a+1}^{(n)},\ldots,\st_{a+d_u-1}^{(n)}$ are precisely all finite
$\st_i^{(n)}$. Define $b$ and $d_v$ similarly in terms of $\lt_i^{(n)}$.
Then similar to the proof of \eqref{claim eq} it can be shown that
$d:=d_u=d_v$ and
\begin{equation*}
\lt_{b+i}^{(n)}=\st_{a+i}^{(n)}-1 \qquad \text{for $0\le i<d$.}
\end{equation*}
(For type $\Cn$, $\st_{a+i}^{(n)}\neq 1$ since by property 2 all
parts in $\nu^{(n)}$ are even. For type $\At$ observe that
$\st_{a+i}^{(n)}\neq 1$ also due to the application of
$j_{r,1}^{-1}\circ j_{r+1,1}^{-1}\circ \cdots \circ j_{r_L-1,1}^{-1}$
to $(\nu,J)$).
This proves in particular that $(\nub,\Jb)$ satisfies point 3
for type $\At$ and points 2 and 3 of Theorem \ref{thm rc prop} for
type $\Cn$. It also shows that $|u|_{[n]}|-|v|_{[n]}|=(n-d)-(r-d)=n-r$
as desired.

Finally, Lemma \ref{lem nb} holds with $s_i^{(k)}$ and $\ell_i^{(k)}$
replaced by $\st_i^{(k)}$ and $\lt_i^{(k)}$, respectively. The proof
is similar to the proof of Lemma \ref{lem nb}. From this follows
that $\lt_{b+i}^{(n+j)}\le \st_{a+i-j}^{(n-j)}$. In particular, this
implies that $\st_{a+i-j}^{(n-j)}=\infty$ if $\lt_{b+i}^{(n+j)}=\infty$.
Let $c$ be minimal such that $\lt_{b+i}^{(n+c)}=\infty$, that is,
$n+c=v_{b+i}$. Then the above shows that
\begin{equation}\label{uc}
u_{a+i-c}\le n-c.
\end{equation}
By the contragredient self-duality $(uv)|_{[n]}$ is a tableau
if and only if $(uv)|_{[n+1,2n]}$ is a tableau.
Now assume that $(uv)|_{[n+1,2n]}$ is not a tableau, that is,
there exists an index $i$ such that $u_{a+i'}\le v_{b+i'}$
for all $0\le i'<i$ and $u_{a+i}> v_{b+i}$. Let $v_{b+i}=n+c$.
Then it follows from duality that $u_{a+i-c}>n-c$. This contradicts
\eqref{uc}. Hence $(uv)|_{[n]}$ must be a tableau as asserted.

\end{document}